\documentclass[12pt,twoside,english]{extarticle}
\usepackage{lmodern}
\usepackage{helvet}
\usepackage[T1]{fontenc}
\usepackage[latin9]{inputenc}
\usepackage{geometry}
\usepackage{color}
\usepackage{babel}
\usepackage{float}
\usepackage{mathrsfs}
\usepackage{amsmath}
\usepackage{amsthm}
\usepackage{amssymb}
\usepackage{graphicx}
\usepackage{setspace}
\usepackage[unicode=true,pdfusetitle,
 bookmarks=true,bookmarksnumbered=false,bookmarksopen=false,
 breaklinks=true,pdfborder={0 0 0},pdfborderstyle={},backref=false,colorlinks=true]
 {hyperref}

\makeatletter

\providecommand{\tabularnewline}{\\}

\numberwithin{equation}{section}
\numberwithin{table}{section}
\numberwithin{figure}{section}
\usepackage[natbibapa]{apacite}
\theoremstyle{definition}
\newtheorem{defn}{\protect\definitionname}[section]
\theoremstyle{plain}
\newtheorem{assumption}{\protect\assumptionname}
\theoremstyle{plain}
\newtheorem{lem}{\protect\lemmaname}[section]
\theoremstyle{definition}
\newtheorem{condition}{\protect\conditionname}
\theoremstyle{plain}
\newtheorem{thm}{\protect\theoremname}[section]
\theoremstyle{plain}
\newtheorem{prop}{\protect\propositionname}[section]
\theoremstyle{plain}
\newtheorem{lyxalgorithm}{\protect\algorithmname}

\usepackage{graphicx}
\usepackage{amsmath}
\usepackage{dcolumn}
\usepackage{listings}
\usepackage{color}
\usepackage{setspace}
\usepackage[normalem]{ulem}  
\definecolor{hellgelb}{rgb}{1,1,0.8}
\definecolor{colKeys}{rgb}{0,0,1}
\definecolor{colIdentifier}{rgb}{0,0,0}
\definecolor{colComments}{rgb}{1,0,0}
\definecolor{colString}{rgb}{0,0.5,0}

\usepackage{lmodern}
\usepackage[T1]{fontenc}
\usepackage{geometry}
\usepackage{fancyhdr}
\usepackage{amsthm}
\usepackage{thmtools}
\usepackage{amsmath}
\usepackage{amssymb}
\usepackage{esint}
\usepackage{multirow}
\usepackage{mathrsfs} 
\usepackage{textgreek}
\usepackage{amsfonts}
\usepackage{amssymb}
\usepackage{mathrsfs} 

\usepackage[USenglish]{isodate}
\usepackage{chngcntr}

\numberwithin{equation}{section}
\numberwithin{table}{section}
\numberwithin{assumption}{section}
\counterwithout{figure}{section}
\counterwithout{table}{section}

\makeatother

  \providecommand{\algorithmname}{Algorithm}
  \providecommand{\assumptionname}{Assumption}

  \providecommand{\definitionname}{Definition}

  \providecommand{\lemmaname}{Lemma}

  \providecommand{\propositionname}{Proposition}

  \providecommand{\theoremname}{Theorem}
 
 \providecommand{\theoremname}{Theorem}
 
\usepackage{thmtools}

\newtheoremstyle{MyTheoremstyle}
  {\topsep} 
  {\topsep} 
  {} 
  {} 
  {\bfseries} 
  {.} 
  {.90em} 
  {} 
\theoremstyle{MyTheoremstyle} 
\theoremstyle{MyTheoremstyle} 
\theoremstyle{MyTheoremstyle} 
\theoremstyle{MyTheoremstyle} 
\theoremstyle{MyTheoremstyle}

\declaretheoremstyle[
    headfont=\bfseries,
    notefont=\normalfont,
    bodyfont=\itshape,
    headpunct=\newline,
    headformat={%
        \makebox{\NAME\ \NUMBER\ }{\NOTE}%
    },
]{theorem}

\newlength{\spacelength}
\declaretheoremstyle[
    headfont=\bfseries,
    notefont=\normalfont,
    bodyfont=\itshape,
    headpunct=\newline,
    headformat={%
        \makebox[0pt][l]{\NAME\ \NUMBER\ }\hskip-\spacelength{\NOTE}%
    },
]{theore}

 \usepackage{fancyhdr}
\usepackage{booktabs}
\usepackage{float}
\usepackage{graphicx}
\usepackage{epstopdf}
\usepackage{morefloats}
\usepackage[referable]{threeparttablex}
\usepackage{footnote}

\usepackage{setspace} 
\usepackage{graphicx,color}

\usepackage{listings}
\RequirePackage[T1]{fontenc}
\RequirePackage{ae,fancyvrb}
\DefineVerbatimEnvironment{Sinput}{Verbatim}{fontshape=sl}
\DefineVerbatimEnvironment{Soutput}{Verbatim}{}
\DefineVerbatimEnvironment{Scode}{Verbatim}{fontshape=sl}

\title{\bf Theory of Evolutionary Spectra for Heteroschesdasticity and Autocorrelation Robust Inference in possibly Misspecified and Nonstationary Models}
\author{
\textsc{\textcolor{MyBlue}{Alessandro Casini}}\thanks{Corresponding author at: Department of Economics and Finance, University of Rome Tor Vergata, Via Columbia 2, Rome, 00133, IT. 
Email: 
\texttt{\textcolor{MyBlue}{{alessandro.casini@uniroma2.it}}}.} 
\\
\small{University of Rome Tor Vergata}
\and
\textsc{\textcolor{MyBlue}{Pierre Perron}}\thanks{Department of Economics, Boston University, 270 Bay State Road, Boston, MA 02215, US. 
Email: 
\texttt{\textcolor{MyBlue}{\mbox{perron@bu.edu}}}.} 
\\
\small{Boston University}
}

\date{\small{\today}}

 \makeatletter

\numberwithin{equation}{section}
\makeatother

\usepackage{babel}
\addto\captionsenglish{%
}

\usepackage{color}
\usepackage{xcolor}

\renewcommand*{\thesection}{\arabic{section}}

\definecolor{MyRed}{rgb}{0.8,0,0}
\definecolor{MyBlue}{rgb}{0,0,0.7}
\definecolor{Green}{rgb}{0,0.5,0}
\definecolor{hellgelb}{rgb}{1,1,0.8}
\definecolor{colKeys}{rgb}{0,0,1}
\definecolor{colIdentifier}{rgb}{0,0,0}
\definecolor{colComments}{rgb}{1,0,0}
\definecolor{colString}{rgb}{0,0.5,0}

\definecolor{MyLightRed}{rgb}{2.2,0.2,0.4} 
\definecolor{MyLightRed2}{rgb}{0.6,0.2,0.3} 
\definecolor{MyLightRed2temp}{rgb}{0.6,0.2,0.3}
\definecolor{MyLightRed3}{rgb}{0.8,0.1,0.1} 
\definecolor{MyRed}{rgb}{0.7,0.0,0}

\definecolor{MyLigthBlue13}{rgb}{0,0.2,0.7}
 \definecolor{MyLigthBlack}{rgb}{0.2,0.25,0.3} 

\hypersetup{%
  bookmarks=true
  pdftitle = {Title Here},
  pdfsubject = {},
  pdfkeywords = {Keywords},
  pdfauthor = {Alessandro Casini},}
\hypersetup{ 
  colorlinks = {true}, 
  citecolor = {MyBlue},
  linkcolor = {MyRed},
  filecolor= {Green},	
  urlcolor = {MyBlue},
  hyperindex = {true},
  linktocpage = {true},
  linkbordercolor ={1 0 0},
 citebordercolor ={0 1 1},
 urlbordercolor ={0 1 1}
}

\usepackage{bookmark}

\usepackage[bottom]{footmisc}
\usepackage{multibib}
\newcites{ReferencesSupp}{References}

\makeatother

\providecommand{\algorithmname}{Algorithm}
\providecommand{\assumptionname}{Assumption}
\providecommand{\conditionname}{Condition}
\providecommand{\definitionname}{Definition}
\providecommand{\lemmaname}{Lemma}
\providecommand{\propositionname}{Proposition}
\providecommand{\theoremname}{Theorem}

\begin{document}
\pagebreak{}

\setcounter{page}{0}

\raggedbottom
\title{\textbf{Change-Point Analysis of Time Series with Evolutionary Spectra}\thanks{We thank an Associate Editor and two referees for helpful comments.
We also thank Federico Belotti, Leopoldo Catania and Stefano Grassi
for generous help with computer programming. }}
\maketitle
\begin{abstract}
This paper develops change-point methods for the spectrum of a locally
stationary time series. We focus on series with a bounded spectral
density that change smoothly under the null hypothesis but exhibits
change-points or becomes less smooth under the alternative. We address
two local problems. The first is the detection of discontinuities
(or breaks) in the spectrum at unknown dates and frequencies. The
second involves abrupt yet continuous changes in the spectrum over
a short time period at an unknown frequency without signifying a break.
Both problems can be cast into changes in the degree of smoothness
of the spectral density over time. We consider estimation and minimax-optimal
testing. We determine the optimal rate for the minimax distinguishable
boundary, i.e., the minimum break magnitude such that we are able
to uniformly control type I and type II errors. We propose a novel
procedure for the estimation of the change-points based on a wild
sequential top-down algorithm and show its consistency under shrinking
shifts and possibly growing number of change-points. 
\end{abstract}
%

\onehalfspacing
\thispagestyle{empty}
\allowdisplaybreaks

\pagebreak{}

\section{Introduction}

Classical change-point theory focuses on detecting and estimating
structural breaks in the mean or regression coefficients. Early contributions
include, among others, \citet{hinkley:71}, \citet{yao:87}, \citet{andrews:93},
\citet{horvath:93} and \citet{bai/perron:98}, who assume the presence
of a single or multiple change-points in the parameters of an otherwise
stationary time series model; see the reviews of \citet{aue/Horvath:13}
and \citet{casini/perron_Oxford_Survey} for more details. More recently
there has been a growing interest about functional and time-varying
parameter models where the latter  are characterized by infinite-dimensional
parameters which change continuously over time {[}see, e.g., \citet{dahlhaus:96},
\citet{neumann/von_sachs:1997}, \citet{hormann/kokoszka:2010}, \citet{dette/preus/vetter:2011},
\citet{zhang/wu:12}, \citet{panaretos/tavakoli:13}, \citet{aue/hormann:2015}
and \citet{vandelft/eichler:2018}{]}. Several authors have extended
the stationarity tests originally introduced by \citet{priesley/rao:1969},
and further developed by, e.g., \citet{dwivedi/rao:2010}, \citet{jentsch/subbarao:2015}
and \citet{bandyopadhyay/carsten/rao}, to these settings. In the
context of local stationarity, \citet{paparoditis:2009} proposed
a test based on comparing a local estimate of the spectral density
to a global estimate and \citet{preuss/vetter/dette:2013} proposed
a test for stationarity using empirical process theory. In the context
of functional time series, tests for stationarity were considered
by \citet{horvath/kokoszka/rice:2014} and \citet{aue/rice/sonmez}
using time-domain methods, and by \citet{aue/vandelft:2020} and \citet{vandelft/Characiejus/dette:2018}
using frequency-domain methods. 

There is wide agreement in empirical work that besides breaks in the
mean the detection of breaks in the variance or the correlation structure
of a time series is of importance. For example, the discrimination
between regimes of high and low asset volatility is of central interest
in finance and the detection of changes in the parameters of an autoregressive
process is important to build superior forecasting procedures. In
addition, discerning the type of the changes (continuous or abrupt
changes) is useful in applied work. On the one hand, if one assumes
local stationarity but the true data-generating process involves structural
breaks then parameter estimates may be severely biased and inference
may be misleading. On the other hand, if one assumes a structural
break model but the parameters actually change gradually then similar
issues may arise. A general approach that allows for both continuous
changes as well as breaks is needed to avoid these issues. Therefore,
the detection of breaks in an otherwise locally stationary time series
is important. 

We develop inference methods about the changes in the degree of smoothness
of the spectrum of a locally stationary time series, and hence, about
change-points in the spectrum as a special case. The key parameter
is the regularity exponent that governs how smooth the path of the
local spectral density is over time. We address two local problems.
The first is the detection of discontinuities (or breaks) in the spectrum
at some unknown date and frequency. The second involves the detection
of abrupt yet continuous changes in the spectrum over a short time
period at an unknown frequency without signifying a break. For example,
the spectrum becomes rougher over a short time period, meaning that
the paths are less smooth as quantified by the regularity exponent.
This can occur for a stationary process whose parameters start to
evolve smoothly according to Lipschitz continuity, or for a locally
stationary process with Lipschitz parameters that change to continuous
but non-differentiable functional parameters. For example, the volatility
of high-frequency stock prices or of other macroeconomic variables
is known to become rougher (i.e., less smooth) without signifying
a structural break after central banks' official announcements, especially
in periods of high market uncertainty. In seismology, earthquakes
are made up of several seismic waves that arrive at different times
and so changes in the smoothness properties of each wave is important
for locating the epicenter and identifying what materials the waves
have passed through. We consider minimax-optimal testing and estimation
for both problems, following \citet{ingster:93}. We determine the
optimal rate for the minimax distinguishable boundary, i.e., the minimum
break magnitude such that we are still able to uniformly control type
I and type II errors. These results are different from the developments
on minimax optimality obtained recently in the statistics literature
for the classical change-point problem where the mean of the series
is piecewise constant {[}see, e.g., \citet{liu/gao/samworth} and
\citet{verzelen/fromont/lerasle/reynaud-bouret:2020}{]}. 

The problem of discriminating discontinuities from a continuous evolution
in a nonparametric framework has received relatively less attention
than the classical change-point problem with a few exceptions {[}\citet{muller:92},
\citet{spokoiny:98}, \citet{muller/stadtmuller:99}, \citet{wu/zhao:07}
and \citet{bibinger/jirak/vetter:16}{]}. These focused on nonparametric
regression and high-frequency volatility, and considered time-domain
methods while we consider frequency-domain methods. This adds a difficulty
in that, e.g., the search for a break or smooth change has to run
over two dimensions, the time and frequency indices. Our test statistics
are the maximum of local two-sample $t$-tests based on the local
smoothed periodogram. We construct statistics that allow the researcher
to test for a change-point in the spectrum at a prespecified frequency
and others that allow to detect a break in the spectrum without prior
knowledge about the frequencies. These test statistics can detect
both discontinuous and smooth changes, and therefore are useful for
both inference problems discussed above. The asymptotic null distribution
follows an extreme value distribution. In order to derive this result,
we first establish several asymptotic results, including bounds for
higher-order cumulants and spectra of locally stationary processes.
These results are complementary to some in \citet{dahlhaus:96}, \citet{paparoditis:2009},
\citet{panaretos/tavakoli:13}, \citet{aue/vandelft:2020} and \citet{casini_hac},
and extend some classical frequency-domain results for stationary
processes {[}e.g., \citet{brillinger:75}{]} to locally stationary
processes.

Change-point problems have also been studied in the frequency-domain
in several fields, though with less generality. \citet{adak:98} investigated
the detection of change-points in piecewise stationary time series
by looking at the difference in the power spectral density for two
adjacent regimes. She compared several distance metrics such as the
Kolmogorov-Smirnov, Cr\'{a}mer-Von Mises and CUSUM-type distance
proposed by \citet{coates/diggle}. \citet{last/shumway:08} focused
on detecting change-points in piecewise locally stationary series.
They exploited some of the results in \citet{Kakizawa/shumway/taniguchi}
and \citet{huang/ombao/stoffer:2004} to propose a Kullback-Liebler
discrimination information but did not derive the null distribution
of the test statistic. \citet{dette/wu/zhou:2019} considered testing
for change-points in the autocorrelation coefficient at some pre-specified
lag $k$ while allowing for local stationarity under the null hypothesis.
They proposed a time-domain method based on a CUSUM-type test. To
the extent that the autocorrelation at a given lag is only one of
the features contained in the spectrum, our setting is more general.
\citet{zhang:16} considered testing for and estimating change-points
in the mean of a piecewise locally stationary time series. This is
a different problem from ours since the spectrum involves the second-order
properties which are complex to study. \citet{preuss/puchstein/dette:2015}
considered the detection and estimation of change-points in the autocovariance
function but they required stationarity under the null hypothesis.
Several authors considered methods based on segmenting the wavelet
spectrum for piecewise stationary time series {[}see, e.g., \citet{barigozzi/cho/fryzlewicz:2018}
and \citeauthor{cho/fryzlewicz:2012} (\citeyear{cho/fryzlewicz:2012};
\citeyear{cho/fryzlewicz:2017}){]}. Outside of the wavelet context
other contributions are \citet{kirch/muhsal/ombao:2015} and \citet{Schrhoder/ombao:2019}.
\citet{vogt/dette:2015} investigated the detection of gradual changes
in a locally stationary time series. Our contribution is different
since we provide a general change-point analysis about the time-varying
spectrum of a time series and establish the relevant asymptotic theory
of the proposed test statistics under both the null and alternative
hypotheses. 

We also address the problem of estimating the change-points, allowing
their number to increase with the sample size and the distance between
change-points to shrink to zero. We propose a procedure based on
a wild sequential top-down algorithm that exploits the idea of bisection
combining it with a wild resampling technique similar to the one  proposed
by \citet{fryzlewicz:14}. We establish the consistency of the procedure
for the number of change-points and their locations. We compare the
rate of convergence with that of standard change-point estimators
under the classical setting {[}e.g., \citet{yao:87}, \citet{bai:94a},
\citet{casini/perron_CR_Single_Break}, \citet{casini/perron_Lap_CR_Single_Inf}
and \citet{casini/perron_SC_BP_Lap}{]}. We verify the performance
of our methods via simulations which show their benefits. The advantage
of using frequency-domain methods to detect change-points is that
they do not require to make assumptions about the data-generating
process under the null hypothesis beyond the fact that the spectrum
is bounded. Furthermore, the method allows for a broader range of
alternative hypotheses than time-domain methods which usually have
good power only against some specific alternatives. For example, tests
for changes in the volatility do not have power for changes in the
dependence and vice versa. Our methods are readily available for use
in many fields such as speech processing, biomedical signal processing,
seismology, failure detection, economics and finance. It can also
be used as a pre-test before constructing the recently introduced
double kernel long-run variance estimator that accounts more flexibly
for nonstationarity {[}cf. \citet{casini_hac}, \citet{casini/perron_PrewhitedHAC}
and \citet{casini/perron_Low_Frequency_Contam_Nonstat:2020}{]}. 

The rest of the paper is organized as follows. Section \ref{Section Hypotheses Testing Problem}
introduces the statistical setting and the hypothesis testing problems.
Section \ref{Section Change-Point-Tests} presents the test statistics
and states their null limit distributions. Section \ref{Section Consistency-and-Minimax}
addresses the consistency of the tests and their minimax optimality.
Section \ref{Section Estimation-of-the Breaks} discusses the estimation
of the change-points while Section \ref{Section Implementation} provides
details for the implementation of the methods. The results of some
Monte Carlo simulations are presented in Section \ref{Section Monte Carlo}.
An empirical application is presented in Section \ref{Section Empirical Application}.
Section \ref{Section Conclusions} reports brief concluding comments.
An online supplement {[}cf. \citet{casini/perron:change-point-spectra_Supp}{]}
contains additional theoretical and empirical results, and all mathematical
proofs. 

\section{\label{Section Hypotheses Testing Problem}Statistical Environment
and the Testing Problems}

Section \ref{Subsec Segmented-Locally-Stationary} introduces the
statistical setting and Section \ref{subsec The hypotheses testin problem}
presents the hypotheses testing problems. We work in the frequency-domain
under the locally stationary framework introduced by \citet{dahlhaus:96}
who formalized the ideas of \citet{priestley:65}. \citet{casini_hac}
extended his framework to allow for discontinuities in the spectrum
which then results in a segmented locally stationary process. This
corresponds to the relevant process under the alternative hypothesis
of breaks in the spectrum. Since local stationarity is a special case
of segmented local stationarity we begin with the latter. We use an
infill asymptotic setting whereby we rescale the original discrete
time horizon $\left[1,\,T\right]$ by dividing each $t$ by $T.$

\subsection{\label{Subsec Segmented-Locally-Stationary}Segmented Locally Stationary
Processes}

Suppose $\left\{ X_{t}\right\} _{t=1}^{T}$ is defined on an abstract
probability space $\left(\Omega,\,\mathscr{F},\,\mathbb{P}\right)$,
where $\Omega$ is the sample space, $\mathscr{F}$ is the $\sigma$-algebra
and $\mathbb{P}$ is a probability measure. Let $i\triangleq\sqrt{-1}$.
We use the notation $\overline{A}$ for the complex conjugate of $A\in\mathbb{C}$.
\begin{defn}
\label{Definition Segmented-Locally-Stationary}A sequence of stochastic
processes $\{X_{t,T}\}_{t=1}^{T}$ is called segmented locally stationary
(SLS)\textbf{ }with $m_{0}+1$ regimes, transfer function $A^{0}$
 and trend $\mu$ if there exists a representation 
\begin{align}
X_{t,T} & =\mu_{j}\left(t/T\right)+\int_{-\pi}^{\pi}\exp\left(i\omega t\right)A_{j,t,T}^{0}\left(\omega\right)d\xi\left(\omega\right),\qquad\qquad\left(t=T_{j-1}^{0}+1,\ldots,\,T_{j}^{0}\right),\label{Eq. Spectral Rep of SLS}
\end{align}
for $j=1,\ldots,\,m_{0}+1$, where by convention $T_{0}^{0}=0$ and
$T_{m_{0}+1}^{0}=T$ ($\mathcal{T}\triangleq\{T_{1}^{0},\,\ldots,\,T_{m_{0}}^{0}\}$),
and the following holds: 

(i) $\xi\left(\omega\right)$ is a stochastic process on $\left[-\pi,\,\pi\right]$
with $\overline{\xi\left(\omega\right)}=\xi\left(-\omega\right)$
and 
\begin{align*}
\mathrm{cum}\left\{ d\xi\left(\omega_{1}\right),\ldots,\,d\xi\left(\omega_{r}\right)\right\}  & =\varphi\left(\sum_{j=1}^{r}\omega_{j}\right)g_{r}\left(\omega_{1},\ldots,\,\omega_{r-1}\right)d\omega_{1}\cdots d\omega_{r},
\end{align*}
 where $\mathrm{cum}\left\{ \cdot\right\} $ is the cumulant  of
$r$th order, $g_{1}=0,\,g_{2}\left(\omega\right)=1$, $\left|g_{r}\left(\omega_{1},\ldots,\,\omega_{r-1}\right)\right|\leq M_{r}<\infty$
  and $\varphi\left(\omega\right)=\sum_{j=-\infty}^{\infty}\delta\left(\omega+2\pi j\right)$
is the period $2\pi$ extension of the Dirac delta function $\delta\left(\cdot\right)$.

(ii) There exists a constant $K>0$ and a piecewise continuous function
$A:\,\left[0,\,1\right]\times\mathbb{R}\rightarrow\mathbb{C}$ such
that, for each $j=1,\ldots,\,m_{0}+1$, there exists a $2\pi$-periodic
function $A_{j}:\,(\lambda_{j-1}^{0},\,\lambda_{j}^{0}]\times\mathbb{R}\rightarrow\mathbb{C}$
with $A_{j}\left(u,\,-\omega\right)=\overline{A_{j}\left(u,\,\omega\right)}$,
$\lambda_{j}^{0}\triangleq T_{j}^{0}/T$ and, for all $T,$
\begin{align}
A\left(u,\,\omega\right)=A_{j}\left(u,\,\omega\right) & \,\mathrm{\,for\,}\,\lambda_{j-1}^{0}<u\leq\lambda_{j}^{0},\label{Eq A(u) =00003D Ai}\\
\sup_{1\leq j\leq m_{0}+1}\sup_{T_{j-1}^{0}<t\leq T_{j}^{0}}\sup_{\omega\in\left[-\pi,\,\pi\right]}\left|A_{j,t,T}^{0}\left(\omega\right)-A_{j}\left(t/T,\,\omega\right)\right| & \leq K\,T^{-1}.\label{Eq. 2.4 Smothenss Assumption on A}
\end{align}

(iii) $\mu_{j}\left(t/T\right)$ is piecewise continuous. 
\end{defn}
The smoothness properties of $A$ in $u$ guarantees that $X_{t,T}$
has a piecewise locally stationary behavior. This means that $X_{t,T}$
is locally stationary in each segment where the notion of local stationarity
is as introduced by \citet{dahlhaus:96}. We refer to \citet{casini_hac}
for several theoretical properties of SLS processes. \citet{zhou:2013}
considered  piecewise locally stationary processes in a time-domain
setting but his notion is less general. In particular, \citet{casini_hac}
also defined (and worked with) the covariance between observations
belonging to different segments whereas previous works considered
only the covariance between observations belonging to the same segment,
thereby using smoothness which simplifies the analysis.

\subsection{\label{subsec The hypotheses testin problem}The Testing Problems}

We focus on time-varying spectra that are bounded, thereby excluding
unit root, and long memory processes. Unit roots or trending processes
can be handled by, for example, taking first differences or using
other de-trending techniques. For some $D<\infty$, we consider the
following class of time-varying spectra under the null hypothesis,
\begin{align}
\boldsymbol{F}\left(\theta,\,D\right) & =\left\{ \left\{ f\left(u,\,\omega\right)\right\} _{u\in\left[0,\,1\right],\,\omega\in\left[-\pi,\,\pi\right]}:\,\sup_{\omega\in\left[-\pi,\,\pi\right],\,}\sup_{u,\,v\in\left[0,\,1\right],\,\left|v-u\right|<h}\left|f\left(u,\,\omega\right)-f\left(v,\,\omega\right)\right|\leq Dh^{\theta}\right\} .\label{eq (2) BJV}
\end{align}
At all continuity points $u$, $f\left(u,\,\omega\right)$ is related
to $A\left(u,\,\omega\right)$ by the relationship $f\left(u,\,\omega\right)=|A\left(u,\,\omega\right)|^{2}$.
The key parameter of the testing problem under the null hypotheses
is $\theta>0$. This is the regularity exponent of $f$ in the time
dimension. For $\theta>1$, $f$ is constant in $u$ and reduces to
the spectral density of a stationary process. For $\theta=1$, $f$
is Lipschitz continuous in $u$. For $\theta<1$, $f$ is $\theta$-H\"older
continuous. Local stationarity corresponds to $\theta>0$ and $f$
being differentiable {[}see \citet{dahlhaus:96b}{]}. The latter is
the setting that we consider under the null hypothesis. To avoid redundancy,
we do not require differentiability directly for the functions in
$\boldsymbol{F}\left(\theta,\,D\right)$ since below we assume that
the transfer function $A\left(u,\,\omega\right)$ is differentiable
in $u$ which in turn implies that $f$ is differentiable in $u$.
Since most of the applied work concerning local stationarity relies
on Lipschitz continuity (i.e., $\theta=1$), our specification of
the null is more general.

We now discuss features that are relevant under the alternative hypothesis.
Our focus is on (i) discontinuities of $f$ in $u$ which correspond
to $\theta=0$ and (ii) decreases in the smoothness of the trajectory
$u\mapsto f(u,\,\omega)$ for each $\omega$ which correspond to a
decrease in $\theta$. We shall refer to a series affected by the
changes of case (ii) as ``becoming more rough or less smooth''.
Both cases refer to the properties of the spectral density and thus
to the second-order properties of $X_{t,T}$. 

Case (i) involves a break in the spectrum, i.e., there exits $\lambda_{b}^{0}\in\left(0,\,1\right)$
such that $\Delta f\left(\lambda_{b}^{0},\,\omega\right)\triangleq(f(\lambda_{b}^{0},\,\omega)-\lim_{u\downarrow\lambda_{b}^{0}}f(u,\,\omega))\neq0$
for some $\omega\in\left[-\pi,\,\pi\right]$. 

Case (ii) involves a fall in the regularity exponent from $\theta>0$
to $\theta'\in\left(0,\,\theta\right)$ after some $\lambda_{b}^{0}$
for some period of time and some $\omega$; i.e., the spectrum becomes
rougher after some $\lambda_{b}^{0}\in\left(0,\,1\right)$ for some
time period before returning to $\theta$-smoothness. The case of
an increase in $\theta$ is technically more complex to handle (see
Section \ref{Section Consistency-and-Minimax} for details). As an
example, consider a locally stationary AR(1),
\begin{align*}
X_{t,T} & =a\left(t/T\right)X_{t-1,T}+\sigma\left(t/T\right)e_{t},\qquad t=1,\ldots,\,T,
\end{align*}
where $a:\,\left[0,\,1\right]\rightarrow\left(-1,\,1\right)$ and
$\sigma:\,\left[0,\,1\right]\rightarrow\mathbb{R}_{+}$ are functional
parameters satisfying a Lipschitz condition and $\left\{ e_{t}\right\} $
is an i.i.d. mean-zero sequence. Additionally, if $\sup_{u\in\left[0\,1\right]}|\sigma\left(u\right)|<\infty$
and the initial condition satisfies some regularity condition, then
$f\left(u,\,\omega\right)$ is uniformly bounded and $\theta=1$.
Problem (ii) refers to either $a\left(\cdot\right)$ or $\sigma\left(\cdot\right)$,
or both, becoming less smooth, i.e., we have a change from $\theta=1$
(since the parameters satisfy a Lipschitz condition) to some $\theta'\in\left(0,\,1\right)$.
In this case, $X_{t,T}$ becomes an AR(1) process with functional
parameters that are still continuous but less smooth and may not be
differentiable. 

Case (i) has received most attention so far in the time series literature
although under much stronger assumptions {[}e.g.,  $f\left(u,\,\omega\right)=f\left(\omega\right)${]}.
Case (ii) is a new testing problem and can be of considerable interest
in several fields even though it requires larger sample sizes than
problem (i). For example, the volatility dynamics of financial or
macroeconomic variables vary over time. \citet{bibinger/jirak/vetter:16}
provided evidence that the volatility of stock prices can change substantially
its path properties after a press release following a meeting of the
Federal Open Market Committee, in particular, its path may become
more rough. Since $f\left(u,\,\omega\right)$ is a smooth function
of the parameters of the data-generating process, if $\sigma\left(u\right)$
becomes more rough, then also $f\left(u,\,\omega\right)$ becomes
more rough as $u$ varies. Case (ii) can also be relevant in seismology
since the study of the path properties of the seismic waves is important
for locating the epicenter of an earthquake.  We show below that
our tests are consistent and have minimax optimality properties for
both cases (i) and (ii). Note that case (ii) is a local problem. In
this paper, we do not consider more global problems where for example
the spectrum is such that a fall in $\theta$ to $\theta'\in\left(0,\,\theta\right)$
occurs on $(\lambda_{b}^{0},\,1]$. This represents a continuous change
in the smoothness of the spectrum that persists  until the end of
the time interval. Different test statistics are needed for this case,
as will be discussed later. 

As discussed by \citet{last/shumway:08}, an important question is
which magnitude of the discontinuity in the time-varying spectrum\textcolor{red}{{}
}can be detected. Or equivalently, how much the  spectrum can change
over a short time  without indicating a break. We introduce the
quantity $b_{T}$, called the detection boundary or simply ``rate'',
which is defined as the minimum break magnitude $\Delta f\left(\lambda_{b}^{0},\,\omega\right)$
such that we are still able to uniformly control the type I and type
II errors as indicated below. To address the minimax-optimal testing
{[}cf. \citet{ingster:93}{]}, we now introduce the testing problems
(i) and (ii) and defer a more general treatment to Section \ref{Section Consistency-and-Minimax}.

\subsubsection*{Testing Problem for Case (i)}

Given the discussion above, for some fractional break point $\lambda_{b}^{0}\in\left(0,\,1\right)$
and frequency $\omega_{0}$, and a decreasing sequence $b_{T}$, we
consider the following class of alternative hypotheses:
\begin{align}
\boldsymbol{F}_{1,\lambda_{b}^{0},\omega_{0}} & \left(\theta,\,b_{T},\,D\right)\label{Eq. (3) BJV - H1}\\
 & =\{\left\{ f\left(u,\,\omega\right)\right\} _{u\in\left[0,\,1\right],\,\omega\in\left[-\pi,\,\pi\right]}:\,\left(f\left(u,\,\omega\right)-\Delta f\left(u,\,\omega\right)\right)_{u\in\left[0,\,1\right]}\in\boldsymbol{F}\left(\theta,\,D\right);\nonumber \\
 & \quad|\Delta f\left(\lambda_{b}^{0},\,\omega_{0}\right)|\geq b_{T}\}.\nonumber 
\end{align}
We can then present first the hypothesis testing problem that we wish
to address: 
\begin{align}
\mathcal{H}_{0} & :\,\left\{ f\left(u,\,\omega\right)\right\} _{u\in\left[0,\,1\right],\,\omega\in\left[-\pi,\,\pi\right]}\in\boldsymbol{F}\left(\theta,\,D\right)\label{eq (4) BJV - Testing Problem}\\
\ensuremath{\mathcal{H}_{1}^{\mathrm{B}}} & :\,\exists\lambda_{b}^{0}\in\left(0,\,1\right)\,\mathrm{and}\,\omega_{0}\in\left[-\pi,\,\pi\right]\,\mathrm{with}\,\left\{ f\left(u,\,\omega\right)\right\} _{u\in\left[0,\,1\right],\,\omega\in\left[-\pi,\,\pi\right]}\in\boldsymbol{F}_{1,\lambda_{b}^{0},\omega_{0}}\left(\theta,\,b_{T},\,D\right).\nonumber 
\end{align}
Observe that $\mathcal{H}_{1}^{\mathrm{B}}$ requires at least one
break but allows for multiple breaks even across different $\omega$.
For the testing problem \eqref{eq (4) BJV - Testing Problem}, we
establish the minimax-optimal rate of convergence of the tests suggested
{[}see Ch. 2 in \citet{ingster/suslina:03} for an introduction{]}.
A conventional definition is the following. For a nonrandomized test
$\psi$ that maps a sample $\left\{ X_{t,T}\right\} $ to zero or
one, we consider the maximal type I error 
\begin{align*}
\alpha_{\psi}\left(\theta\right) & =\sup_{\left\{ f\left(u,\,\omega\right)\right\} _{u\in\left[0,\,1\right],\,\omega\in\left[-\pi,\,\pi\right]}\in\boldsymbol{F}\left(\theta,\,D\right)}\mathbb{P}_{f}\left(\psi=1\right),
\end{align*}
 and the maximal type II error
\begin{align*}
\beta_{\psi}\left(\theta,\,b_{T}\right) & =\sup_{\lambda_{b}^{0}\in\left(0,\,1\right),\,\omega_{0}\in\left[-\pi,\,\pi\right]}\sup_{\left\{ f\left(u,\,\omega\right)\right\} _{u\in\left[0,\,1\right],\,\omega\in\left[-\pi,\,\pi\right]}\in\boldsymbol{F}_{1,\lambda_{b}^{0},\omega_{0}}\left(\theta,\,b_{T},\,D\right)}\mathbb{P}_{f}\left(\psi=0\right),
\end{align*}
and define the total testing error as $\gamma_{\psi}(\theta,\,b_{T})=\alpha_{\psi}(\theta)+\beta_{\psi}(\theta,\,b_{T}).$
The notion of asymptotic minimax-optimality is as follows. We want
to find sequences of tests and rates $b_{T}$ such that $\gamma_{\psi}(\theta,\,b_{T})\rightarrow0$
as $T\rightarrow\infty$. The larger is $b_{T}$ the easier it is
to distinguish between $\mathcal{H}_{0}$ and $\mathcal{H}_{1}^{\mathrm{B}}$
 but we may incur at the same time a larger type II error $\beta_{\psi}\left(\theta,\,b_{T}\right)$.
The optimal value $b_{T}^{\mathrm{opt}}$, named the minimax distinguishable
rate, is the minimum value of $b_{T}>0$ such that $\lim_{T\rightarrow\infty}\inf_{\psi}\gamma_{\psi}\left(\theta,\,b_{T}\right)=0$.
A sequence of tests $\psi_{T}$ that satisfies the latter relation
for all $b_{T}\geq b_{T}^{\mathrm{opt}}$ is called minimax-optimal. 

Minimax-optimality has been considered in other change-point problems.
\citet{loader:96} and \citet{spokoiny:98} considered the nonparametric
estimation of a regression function with break size fixed. \citet{bibinger/jirak/vetter:16}
considered breaks in the volatility of semimartingales under high-frequency
asymptotics while we focus on breaks in the spectral density and thus
we work in the frequency-domain. Another difference from previous
work is that we do not deal with i.i.d. observations; hence, we cannot
use the same approach to derive the minimax bound as in \citet{bibinger/jirak/vetter:16}
because their information-theoretic reductions exploit independence.
We need to rely on approximation theorems {[}cf. \citet{berkes/philipp:79}{]}
to establish that our statistical experiment is asymptotically equivalent
in a strong Le Cam sense to a high dimensional signal detection problem.
This allows us to derive the minimax bound using classical arguments
based on the results in \citet{ingster/suslina:03}, Ch. 8. The relevant
results are stated in Section \ref{Section Consistency-and-Minimax}.

\subsubsection*{Testing Problem for Case (ii)}

We consider alternative hypotheses where $f$ is less smooth than
under $\mathcal{H}_{0}$, including the case of breaks as a special
case. Suppose that under $\mathcal{H}_{0}$ the spectrum $f\left(u,\,\omega\right)$
is differentiable in both arguments and behaves until time $T\lambda_{b}^{0}$
as specified in $\boldsymbol{F}\left(\theta,\,D\right)$ for some
$\theta>0$ and $D<\infty$. After $T\lambda_{b}^{0}$, the regularity
exponent $\theta$ drops to some $\theta'$ with $0<\theta'<\theta$
for some non-trivial period of time. That is, since $\boldsymbol{F}\left(\theta,\,D\right)\subset\boldsymbol{F}\left(\theta',\,D\right)$,
we need that $f$ behaves as $\theta'$-regular for some period of
time such that there exists a $\omega$ with $\left\{ f\left(u,\,\omega\right)\right\} _{u\in\left[0,\,1\right]}\notin\boldsymbol{F}\left(\theta,\,D\right).$
This guarantees that $\mathcal{H}_{0}$ and $\mathcal{H}_{1}^{\mathrm{S}}$
(to be defined below) are well-separated. To this end, define for
some function $g_{u}$ with $u\in\left[0,\,1\right]$, $\Delta_{h}^{\theta'}g_{u}=\left(g_{u+h}-g_{u}\right)/\left|h\right|^{\theta'}$
for $h\in\left[-u,\,1-u\right].$ The set of possible alternatives
is then defined as 
\begin{align*}
\boldsymbol{F}'_{1,\lambda_{b}^{0},\omega_{0}}\left(\theta,\,\theta',\,b_{T},\,D\right) & =\biggl\{\left\{ f\left(u,\,\omega\right)\right\} _{u\in\left[0,\,1\right],\,\omega\in\left[-\pi,\,\pi\right]}\in\boldsymbol{F}\left(\theta',\,D\right):\\
 & \quad\quad\,\inf_{\left|h\right|\leq2m_{T}/T}\Delta_{h}^{\theta'}f\left(\lambda_{b}^{0},\,\omega_{0}\right)\geq b_{T}\quad\mathrm{or}\quad\sup_{\left|h\right|\leq2m_{T}/T}\Delta_{h}^{\theta'}f\left(\lambda_{b}^{0},\,\omega_{0}\right)\leq-b_{T}\biggr\},
\end{align*}
 where $m_{T}\rightarrow\infty$ as $T\rightarrow\infty$ such that
$m_{T}^{-1}T^{\epsilon}\rightarrow0$ with $\epsilon>0.$  This leads
to the following testing problem,
\begin{align}
\mathcal{H}_{0} & :\,\left\{ f\left(u,\,\omega\right)\right\} _{u\in\left[0,\,1\right]}\in\boldsymbol{F}\left(\theta,\,D\right)\label{eq (29) BJV - Testing Problem}\\
\mathcal{H}_{1}^{\mathrm{S}} & :\,\exists\lambda_{b}^{0}\in\left(0,\,1\right)\,\mathrm{and}\,\omega_{0}\in\left[-\pi,\,\pi\right]\,\mathrm{with}\,\left\{ \left\{ f\left(u,\,\omega\right)\right\} _{u\in\left[0,\,1\right],\,\omega\in\left[-\pi,\,\pi\right]}\right\} \in\boldsymbol{F}'_{1,\lambda_{b}^{0},\omega_{0}}\left(\theta,\,\theta',\,b_{T},\,D\right).\nonumber 
\end{align}
Note that $\mathcal{H}_{1}^{\mathrm{S}}$ allows for multiple changes.
$\mathcal{H}_{1}^{\mathrm{B}}$ is a special case of $\mathcal{H}_{1}^{\mathrm{S}}$
since it can be seen as the limiting case of $\mathcal{H}_{1}^{\mathrm{S}}$
as $\theta'\rightarrow0.$ In the context of infinite-dimensional
parameter problems one faces the issue of distinguishability between
the null and the alternative hypotheses. It is evident that one cannot
test $f\in\boldsymbol{F}\left(\theta,\,D\right)$ versus $f\in\boldsymbol{F}\left(\theta',\,D\right)$
for $\theta>\theta'$. First, since $\boldsymbol{F}\left(\theta,\,D\right)\subset\boldsymbol{F}\left(\theta',\,D\right)$,
one has at least to remove the set of functions in $\boldsymbol{F}\left(\theta,\,D\right)$
from those in $\boldsymbol{F}\left(\theta',\,D\right)$. Still, as
discussed by \citet{ingster/suslina:03}, this would not be enough
since the two hypotheses are still too close. That explains why we
focus on spectral densities $f$ that belong to $\boldsymbol{F}'_{1,\lambda_{b}^{0},\omega_{0}}\left(\theta,\,\theta',\,b_{T},\,D\right)$
under $\mathcal{H}_{1}^{\mathrm{S}}$. These are rough enough so
as not to be close to functions in $\boldsymbol{F}\left(\theta,\,D\right)$.
This is captured by the requirement that the difference quotient $\Delta_{h}^{\theta'}f$
exceeds the so-called rate $b_{T}$. As $T\rightarrow\infty$ the
requirement becomes less stringent since $b_{T}\rightarrow0$. See
\citet{hoffman/nick:11} and \citet{bibinger/jirak/vetter:16} for
similar discussions in different contexts.

\section{\label{Section Change-Point-Tests}Tests for Changes in the Spectrum
and Their Limiting Distributions }

Section \ref{subsec: Test-Statistics} introduces the test statistics
while Section \ref{subsec:The-Limiting-Distribution} presents the
results concerning their asymptotic distributions under the null hypothesis.
These results apply also to the case of smooth alternatives which
we discuss formally in Section \ref{Section Consistency-and-Minimax}.

\subsection{\label{subsec: Test-Statistics}The Test Statistics}

We first define the quantities needed to define the tests. Let $h:\,\mathbb{R}\rightarrow\mathbb{R}$
be a data taper with $h\left(x\right)=0$ for $x\notin[0,\,1)$, 
\[
H_{k,T}\left(\omega\right)=\sum_{s=0}^{T-1}h\left(s/T\right)^{k}\exp\left(-i\omega s\right),
\]
and (for $n_{T}$ even), 
\begin{align*}
d_{L,h,T}\left(u,\,\omega\right)\triangleq\sum_{s=0}^{n_{T}-1}h\left(\frac{s}{n_{T}}\right)X_{\left\lfloor Tu\right\rfloor -n_{T}+s+1,T}\exp\left(-i\omega s\right), & \quad I_{L,h,T}\left(u,\,\omega\right)\triangleq\frac{1}{2\pi H_{2,n_{T}}\left(0\right)}\left|d_{L,h,T}\left(u,\,\omega\right)\right|^{2},\\
d_{R,h,T}\left(u,\,\omega\right)\triangleq\sum_{s=0}^{n_{T}-1}h\left(\frac{s}{n_{T}}\right)X_{\left\lfloor Tu\right\rfloor +n_{T}-s,T}\exp\left(-i\omega s\right), & \quad I_{R,h,T}\left(u,\,\omega\right)\triangleq\frac{1}{2\pi H_{2,n_{T}}\left(0\right)}\left|d_{R,h,T}\left(u,\,\omega\right)\right|^{2},
\end{align*}
 where $I_{L,h,T}\left(u,\,\omega\right)$ (resp., $I_{R,h,T}\left(u,\,\omega\right)$)
is the local periodogram over a segment of length $n_{T}\rightarrow\infty$
that uses observations to the left (resp. right) of $\left\lfloor Tu\right\rfloor $.
 $I_{L,h,T}\left(u,\,\omega\right)$ (resp., $I_{R,h,T}\left(u,\,\omega\right)$)
is a near unbiased estimator of $f\left(u-n_{T}/T,\,\omega\right)$
(resp., $f\left(u+n_{T}/T,\,\omega\right)$). We allow for a data
taper since one may want to put more weight on observations that are
closer to $u$. The smoothed local periodogram is defined as
\begin{align*}
f_{L,h,T}\left(u,\,\omega\right) & =\frac{2\pi}{n_{T}}\sum_{s=1}^{n_{T}-1}W_{T}\left(\omega-\frac{2\pi s}{n_{T}}\right)I_{L,h,T}\left(u,\,\frac{2\pi s}{n_{T}}\right),
\end{align*}
with $f_{R,h,T}\left(u,\,\omega\right)$ defined similarly to $f_{L,h,T}\left(u,\,\omega\right)$
but with $I_{R,h,T}\left(u,\,\omega\right)$ in place of $I_{L,h,T}\left(u,\,\omega\right)$,
where $W_{T}\left(\omega\right)$ $(-\infty<\omega<\infty)$ is a
family of weight functions of period $2\pi$, 
\begin{align*}
W_{T}\left(\omega\right) & =\sum_{j=-\infty}^{\infty}b_{W,T}^{-1}W\left(b_{W,T}^{-1}\left(\omega+2\pi j\right)\right),
\end{align*}
with $b_{W,T}$ a bandwidth and $W\left(\beta\right)$ $(-\infty<\beta<\infty)$
a fixed function. We define  
\begin{align*}
\widetilde{f}_{L,r,T}\left(\omega\right)=M_{S,T}^{-1}\sum_{j\in\mathbf{S}_{r}}f_{L,h,T}\left(j/T,\,\omega\right) & \quad\mathrm{and}\quad\widetilde{f}_{R,r,T}\left(\omega\right)=M_{S,T}^{-1}\sum_{j\in\mathbf{S}_{r}}f_{R,h,T}\left(j/T,\,\omega\right),
\end{align*}
where 
\begin{align*}
\mathbf{S}_{r} & =\{rm_{T}-m_{T}/2+\left\lfloor n_{T}/2\right\rfloor +1,\,rm_{T}-m_{T}/2+\left\lfloor n_{T}/2\right\rfloor +1+m_{S,T},\\
 & \qquad\ldots,\,rm_{T}+\left\lfloor n_{T}/2\right\rfloor +1+m_{S,T}M_{S,T}/2\},
\end{align*}
with $m_{S,T}=\left\lfloor m_{T}^{1/2}\right\rfloor $ and $M_{S,T}=\left\lfloor m_{T}/m_{S,T}\right\rfloor $.
$\widetilde{f}_{a,r,T}\left(\omega\right)$ ($a=L,\,R$) denotes the
average local spectral density around time $rm_{T}$ computed using
$f_{a,h,T}\left(j/T,\,\omega\right)$ where $r=1,\ldots,\,M_{T}=\left\lfloor T/m_{T}\right\rfloor -1$.
We do not use all the $m_{T}$ local spectral densities $f_{a,\,h,T}\left(j/T,\,\omega\right)$
($a=L,\,R$) in the block $r$ but only those separated by $m_{S,T}$
points. Thus, $\mathbf{S}_{r}$ is a subset of the indices in the
block $r$. We need to consider a sub-sample of the $f_{a,h,T}\left(j/T,\,\omega\right)$'s
($a=L,\,R$) because there is strong dependence among the adjacent
terms, e.g., $f_{a,h,T}\left(j/T,\,\omega\right)$ and $f_{a,h,T}\left((j+1)/T,\,\omega\right)$
($a=L,\,R$). A large deviation between $\widetilde{f}_{L,r,T}\left(\omega\right)$
and $\widetilde{f}_{R,r+1,T}\left(\omega\right)$ suggests the presence
of a  break in the spectrum close to time $\left(r+1\right)m_{T}$
at frequency $\omega$. Note that the latest observation used in the
construction of $\widetilde{f}_{L,r,T}$ is $X_{rm_{T}+m_{T}/2+\left\lfloor n_{T}/2\right\rfloor +1}$
while the earliest observation used in the construction of $\widetilde{f}_{R,r+1,T}$
is $X_{rm_{T}+m_{T}/2+\left\lfloor n_{T}/2\right\rfloor +2}$. This
shows that there is no overlapping in the time points used in the
construction of $\widetilde{f}_{L,r,T}$ and $\widetilde{f}_{R,r+1,T}$.
The reason behind this is that in order to maximize power $\widetilde{f}_{L,r,T}$
and $\widetilde{f}_{R,r+1,T}$ should not use common observations,
otherwise the effect of the common observations to the difference
in the averages would offset the effect of the change-point.

Let
\begin{align*}
\mathbf{S}_{r,+}\left(j\right) & =\{\{rm_{T}-m_{T}/2+\left\lfloor n_{T}/2\right\rfloor +1,\,rm_{T}-m_{T}/2+\left\lfloor n_{T}/2\right\rfloor +1+\widetilde{m}_{S,T},\\
 & \qquad\ldots,\,rm_{T}+\left\lfloor n_{T}/2\right\rfloor +1+\widetilde{m}_{S,T}\widetilde{M}_{S,T}/2\}/\,\{\ldots,\,rm_{T}-m_{T}/2+1+\widetilde{m}_{S,T}\left(j-1\right)\}\},\\
\mathbf{S}_{r,-}\left(j\right) & =\{\{rm_{T}-m_{T}/2+\left\lfloor n_{T}/2\right\rfloor +1,\,rm_{T}-m_{T}/2+\left\lfloor n_{T}/2\right\rfloor +1+\widetilde{m}_{S,T},\\
 & \qquad\ldots,\,rm_{T}+\left\lfloor n_{T}/2\right\rfloor +1+\widetilde{m}_{S,T}\widetilde{M}_{S,T}/2\}/\,\{\ldots,\,rm_{T}-m_{T}/2+1+\widetilde{m}_{S,T}\left(-j+1\right)\}\},
\end{align*}
where $\widetilde{m}_{S,T}=\left\lfloor m_{T}^{1/3}\right\rfloor $
and $\widetilde{M}_{S,T}=\left\lfloor m_{T}/\widetilde{m}_{S,T}\right\rfloor $.
Define $\widehat{\sigma}_{L,r}^{2}\left(\omega\right)=\sum_{j=-\widetilde{M}_{S,T}+1}^{\widetilde{M}_{S,T}-1}K_{1}\left(b_{1,T}j\right)\widehat{\Gamma}_{r}\left(j\right)$
where
\begin{align*}
\widehat{\Gamma}_{r}\left(j\right) & =\begin{cases}
\widetilde{M}_{S,T}^{-1}\sum_{t\in\mathbf{S}_{r,+}\left(j\right)}\widehat{f}_{L,h,T}\left(t/T,\,\omega\right)\widehat{f}_{L,h,T}\left(\left(t-j\widetilde{m}_{S,T}\right)/T,\,\omega\right), & j\geq0\\
\widetilde{M}_{S,T}^{-1}\sum_{t\in\mathbf{S}_{r,-}\left(j\right)}\widehat{f}_{L,h,T}\left(t/T,\,\omega\right)\widehat{f}_{L,h,T}\left(\left(t+j\widetilde{m}_{S,T}\right)/T,\,\omega\right), & j<0
\end{cases},
\end{align*}
 and $\widehat{f}_{L,h,T}\left(j/T,\,\omega\right)=f_{L,h,T}\left(j/T,\,\omega\right)-\widetilde{f}_{L,r,T}\left(\omega\right)$
for $j\in\mathbf{S}_{r}$. The quantity $\widehat{\sigma}_{L,r}^{2}\left(\omega\right)$
is a local long-run variance estimator where $K_{1}$ is a kernel
and $b_{1,T}$ is the associated bandwidth.

We first present a test statistic for the detection of a change-point
in the spectrum $f\left(\cdot,\,\omega\right)$ for a given frequency
$\omega.$ A second test statistic that we consider detects change-points
in $u\in(0,\,1)$ occurring at any frequency $\omega\in\left[-\pi,\,\pi\right]$.
The latter is arguably more useful in practice because often the practitioner
does not know a priori at which frequency the spectrum is discontinuous.
We begin with the following test statistic,
\begin{align}
\mathrm{S}_{\max,T}\left(\omega\right)\triangleq\max_{r=1,\ldots,\,M_{T}-2}\left|\frac{\widetilde{f}_{L,r,T}\left(\omega\right)-\widetilde{f}_{R,r+1,T}\left(\omega\right)}{\widehat{\sigma}_{L,r}\left(\omega\right)}\right| & ,\qquad\omega\in\left[-\pi,\,\pi\right].\label{Eq. (14) Smax Stat}
\end{align}

Test statistics of the form of \eqref{Eq. (14) Smax Stat} were also
used in the time-domain in the context of nonparametric change-point
analysis under a less general framework {[}cf. \citet{eichinger/kirch:2018},
\citet{bibinger/jirak/vetter:16} and \citet{wu/zhao:07}{]} and forecasting
{[}cf. \citet{casini_CR_Test_Inst_Forecast}{]}. The derivation of
the null distribution uses a (strong) invariance principle for nonstationary
processes {[}see, e.g., \nocite{wu:07} and \citet{wu/zhou:11}{]}.

The test statistic $\mathrm{S}_{\mathrm{max},T}\left(\omega\right)$
aims at detecting a break in the spectrum at some given frequency
$\omega$. An alternative would be to consider a double-sup statistic
which takes the maximum over $\omega\in\left[-\pi,\,\pi\right]$.
Theorem \ref{Theorem 5.2.8 Brillinger} in the supplement shows that
$I_{h,L,T}\left(u,\,\omega_{j}\right)$ and $I_{h,L,T}\left(u,\,\omega_{k}\right)$
are asymptotically independent if $2\omega_{j},\,\omega_{j}\pm\omega_{k}\not\equiv0\,(\mathrm{mod\,}2\pi)$.\footnote{The notation $2\omega_{j},\,\omega_{j}\pm\omega_{k}\not\equiv0\,(\mathrm{mod\,}2\pi)$
means $2\omega_{j}\not\equiv0\,(\mathrm{mod\,}2\pi)$ and $\omega_{j}\pm\omega_{k}\not\equiv0\,(\mathrm{mod\,}2\pi).$ } However, the smoothing over frequencies introduces short-range dependence
over $\omega$. Due to this short-range dependence, we cannot consider
the maximum over all frequencies in $\Pi$ because the statistics
would not be independent. Thus, we specify a framework based on an
infill procedure over the frequency-domain $\left[-\pi,\,\pi\right]$
by assuming that there are $n_{\omega}$ frequencies $\omega_{1},\ldots,\,\omega_{n_{\omega}},$
with $\omega_{1}=-\pi$ and $\omega_{n_{\omega}}=\pi-\epsilon,\,\epsilon>0$,
and $\left|\omega_{j}-\omega_{j+1}\right|=O\left(n_{\omega}^{-1}\right)$
for $j=1,\ldots,\,n_{\omega}-2$. Assume that $n_{\omega}\rightarrow\infty$
as $T\rightarrow\infty$.  Let $\Pi\triangleq\left\{ \omega_{1},\ldots,\,\omega_{n_{\omega}}\right\} $.
The maximum is taken over the following set of frequencies 
\[
\Pi'\triangleq\{\omega_{1},\,\omega_{2+\left\lfloor n_{T}b_{W,T}\right\rfloor },\ldots,\omega_{n_{\omega}-\left\lfloor n_{T}b_{W,T}\right\rfloor -1},\,\omega_{n_{\omega}}\}.
\]
Let $n'_{\omega}=\left\lfloor n_{\omega}/\left(\left\lfloor n_{T}b_{W,T}\right\rfloor +1\right)\right\rfloor $
with $n'_{\omega}\rightarrow\infty$. Note that $\Pi'\subset\Pi$.
This then leads to the double-sup statistic,
\begin{equation}
\mathrm{S}_{\mathrm{Dmax},T}\triangleq\max_{\omega_{k}\in\Pi'}\sqrt{\log\left(M_{T}\right)}(M_{S,T}^{1/2}\mathrm{S}_{\mathrm{max},T}\left(\omega_{k}\right)-\gamma_{M_{T}})-\log\left(n'_{\omega}\right),\label{Eq. (SDmax)}
\end{equation}
where $\gamma_{M_{T}}=\left[4\log\left(M_{T}\right)-2\log\left(\log\left(M_{T}\right)\right)\right]^{1/2}$.
This double-sup form  is a new feature for change-point testing under
the frequency-domain. 

Next, we consider alternative test statistics that are self-normalized
such that one does not need to estimate $\sigma_{L,r}^{2}\left(\omega\right)$.
We consider the following test statistic,
\begin{align}
\mathrm{R}_{\max,T}\left(\omega\right)\triangleq\max_{r=1,\ldots,\,M_{T}-2}\left|\frac{\widetilde{f}_{L,r,T}\left(\omega\right)}{\widetilde{f}_{R,r+1,T}\left(\omega\right)}-1\right| & ,\label{Eq. (14) Smax Stat-1}
\end{align}
where $\omega\in\left[-\pi,\,\pi\right].$ We can define a test statistic
corresponding to $\mathrm{S}_{\mathrm{Dmax},T}$ by 
\[
\mathrm{R}_{\mathrm{Dmax},T}\triangleq\max_{\omega_{k}\in\Pi'}\sqrt{\log\left(M_{T}\right)}(M_{S,T}^{1/2}\mathrm{R}_{\mathrm{max},T}\left(\omega_{k}\right)-\gamma_{M_{T}})-\log\left(n'_{\omega}\right).
\]

\subsection{\label{subsec:The-Limiting-Distribution}The Limiting Distribution
Under the Null Hypothesis}

Let $\mathbf{X}_{t,T}=(X_{t,T}^{\left(a_{1}\right)},\ldots,\,X_{t,T}^{\left(a_{p}\right)})$
with finite $p\geq1$. Denote by $\kappa_{\mathbf{X},t}^{\left(a_{1},\ldots,a_{r}\right)}\left(k_{1},\ldots,\,k_{r-1}\right)$
the time-$t$ cumulant of order $r$ of $(X_{t+k_{1},T}^{\left(a_{1}\right)},\ldots,X_{t+k_{r-1},T}^{\left(a_{r-1}\right)},\,X_{t,T}^{\left(a_{r}\right)})$
with $r\leq p$. 
\begin{assumption}
\label{Assumption Locally Stationary}(i) $\left\{ \mathbf{X}_{t,T}\right\} $
is a mean-zero  locally stationary process (i.e., $m_{0}=0$); (ii)
for all $j=1,\ldots,\,p$, $A^{\left(a_{j}\right)}\left(u,\,\omega\right)$
is \textcolor{red}{ } $2\pi$-periodic in $\omega$ and the periodic
extensions are differentiable in $u$ and $\omega$ with uniformly
bounded derivative $\left(\partial/\partial u\right)\left(\partial/\partial\omega\right)A\left(u,\,\omega\right)$;
(iii) $g_{4}$ is continuous. 
\end{assumption}
Assumption \ref{Assumption Locally Stationary} requires $\left\{ \mathbf{X}_{t,T}\right\} $
to be locally stationary. Without loss of generality, we assume that
$\left\{ \mathbf{X}_{t,T}\right\} $ has zero mean. All results go
through when the mean is non-zero or when using demeaned series. The
differentiability of $A\left(u,\,\omega\right)$ implies that $f\left(u,\,\omega\right)$
is also differentiable. This means that under the null hypothesis
we require $f\left(u,\,\omega\right)$ to be differentiable in $u$
and to have some regularity exponent $\theta>0.$ The differentiability
of $A\left(u,\,\omega\right)$ in $u$ can be relaxed at the expense
of more complex proofs to establish the results in the supplement
on high-order cumulants.  Without differentiability, for any $\theta>0$
the test statistics above follow the same asymptotic distribution
as when differentiability holds, though we do not discuss this case
formally.

We need to impose some conditions on the temporal dependence. Let
$\left\{ e_{t}\right\} _{t\in\mathbb{Z}}$ be a sequence of i.i.d.
random variables and $\left\{ e'_{t}\right\} _{t\in\mathbb{Z}}$ be
an independent copy of $\left\{ e_{t}\right\} _{t\in\mathbb{Z}}.$
Assume $X_{t,T}=H_{T}\left(t/T,\,\mathscr{F}_{t}\right)$ where $\mathscr{F}_{t}\triangleq\left\{ \ldots,\,e_{t-1},\,e_{t}\right\} $
and $H_{T}:\,\left[0,\,1\right]\times\mathbb{R}^{\infty}\mapsto\mathbb{R}$
is a measurable function. We use the  dependence measure introduced
by \citeauthor{wu:05} (\citeyear{wu:05}, \citeyear{wu:07}) for
stationary processes and extended to nonstationary processes by \citet{wu/zhou:11}.
Let $\mathscr{L}^{q}$ denote the space generated by the $q$-norm,
$q>0$. For all $t$, assume $X_{t,T}\in\mathscr{L}^{q}$. For $w\geq0$
define the dependence measure,
\begin{align}
\phi_{w,q} & =\sup_{t}\left\Vert X_{t,T}-X_{t,T,\left\{ w\right\} }\right\Vert _{q}\label{Eq. (2.1) WZ Delta}\\
 & =\sup_{t}\left\Vert H_{T}\left(t/T,\,\mathscr{F}_{t}\right)-H_{T}\left(t/T,\,\mathscr{F}_{t,\left\{ w\right\} }\right)\right\Vert _{q},\nonumber 
\end{align}
where $\mathscr{F}_{t,\left\{ w\right\} }$ is a coupled version of
$\mathscr{F}_{t}$ with $e_{w}$ replaced by an i.i.d. copy $e'_{w}$.
For $\{X_{t,T}\}$ locally stationary, $H_{T}\left(\cdot,\,\mathscr{F}_{t}\right)$
is stochastic $\theta$-H\"older continuous in the sense that there
exists $C_{T,H}<\infty$ such that
\begin{align}
\sup_{0\leq u<u'\leq1} & \frac{\left\Vert H_{T}\left(u,\,\mathscr{F}_{t}\right)-H_{T}\left(u',\,\mathscr{F}_{t}\right)\right\Vert }{|u-u'|^{\theta}}\leq C_{T,H}.\label{Eq. Stochastic Lip cont}
\end{align}
Assume $\Upsilon_{n,q}=\sum_{j=n}^{\infty}\phi_{j,q}<\infty$ for
some $n\in\mathbb{Z}$. Let $\tau_{T}=T^{\vartheta_{1}}\left(\log\left(T\right)\right)^{\vartheta_{2}}$
where $\vartheta_{1}=\left(1/2-1/q+\gamma/q\right)$ $/\left(1/2-1/q+\gamma\right)$
and $\vartheta_{2}=\left(\gamma+\gamma/q\right)/\left(1/2-1/q+\gamma\right)$
for some $\gamma>0$. 
\begin{assumption}
\label{Assumption WZ 2011}For $q\geq2$, $\sum_{n=0}^{\infty}n^{l+q-1}(\sum_{j=n}^{\infty}\phi_{j,q}^{2})^{1/2}<\infty$
where $l\geq0$. 
\end{assumption}
Assumption \ref{Assumption WZ 2011} was also used by \citet{shao/wu:2007ET}
who showed that it is satisfied for many nonlinear time series processes.
Using Assumption \ref{Assumption WZ 2011} we can give a sufficient
condition for the summability of the joint cumulant up to a certain
order. The latter is a common assumption in spectral analysis and
we use it to establish results on the high-order cumulants and spectra
in Section \ref{Section Cumulant and Spectra}. These are used to
obtain the null limiting distribution of the test statistics. 
\begin{lem}
\label{Lemma Theorem 4.1 in Shao and Wu (2007)}Let Assumption \ref{Assumption WZ 2011}
hold and $X_{t,T}\in\mathscr{L}^{r}$, $r\geq2$. We have 
\begin{align}
\sum_{k_{1},\ldots,\,k_{r-1}=-\infty}^{\infty}\left(1+\left|k_{j}\right|^{l}\right)\sup_{1\leq t\leq T}\left|\kappa_{\mathbf{X},t}^{\left(a_{1},\ldots,a_{r}\right)}\left(k_{1},\ldots,\,k_{r-1}\right)\right| & <\infty,\label{Eq. (2.6.6) Condition on cumulant}
\end{align}
 for $l\geq0$, $j=1,\ldots,\,r-1$, and any $r$-tuple $a_{1},\ldots,a_{r}$.
\end{lem}
\citet{shao/wu:2007ET} proved \eqref{Eq. (2.6.6) Condition on cumulant}
for $l=0$ and $\{X_{t,T}\}$ a stationary process.
\begin{assumption}
\label{Taper, nT, W and K1, b1}(i) The data taper $h:\,\mathbb{R}\rightarrow\mathbb{R}$
with $h\left(x\right)=0$ for $x\notin[0,\,1)$ is bounded and of
bounded variation; (ii) The sequence $\left\{ n_{T}\right\} $ satisfies
$n_{T}\rightarrow\infty$ as $T\rightarrow\infty$ with $n_{T}/T\rightarrow0$;
(iii) $W\left(\beta\right)$ $(-\infty<\beta<\infty)$ is real-valued,
even, of bounded variation, and satisfies $\int_{-\infty}^{\infty}W\left(\beta\right)d\beta=1$;
(iv) $b_{1,T}\rightarrow0$ such that $Tb_{1,T}\rightarrow\infty$
and $(\widetilde{M}_{S,T}b_{1,T})^{-1/2}M_{S,T}^{1/2}\log T\rightarrow0$
and $K_{1}\left(\cdot\right):\,\mathbb{R}\rightarrow\left[-1,\,1\right],$
$K_{1}\left(0\right)=1,\,K_{1}\left(x\right)=K_{1}\left(-x\right),\,\forall x\in\mathbb{R}$
and ${\textstyle \int\nolimits _{-\infty}^{\infty}}K_{1}^{2}\left(x\right)dx<\infty$.
 
\end{assumption}
Assumption \ref{Taper, nT, W and K1, b1}-(i,ii) are standard in the
nonparametric estimation literature while Assumption \ref{Taper, nT, W and K1, b1}-(iii)
is also used for spectral density estimation under stationarity {[}e.g.,
\citet{brillinger:75}{]}. The conditions on $b_{1,T}$ in Assumption
\ref{Taper, nT, W and K1, b1}-(iv) are necessary for the consistency
of the long-run variance estimator. The class of kernels allowed
by Assumption \ref{Taper, nT, W and K1, b1}-(iv) includes popular
kernels such as the Truncated, Bartlett, Parzen, Quadratic Spectral
(QS) and Tukey-Hanning kernels. For technical reasons inherent to
the proofs, we need to assume that the spectral density is strictly
positive. Theorem \ref{Theorem 5.6.4 Brillinger} in the supplement
shows that the variance of $f_{L,h,T}\left(u,\,\omega\right)$ depends
on $f\left(u,\,\omega\right)$. Thus, the denominator of the test
statistic depends on $f\left(u,\,\omega\right)$. Assumption \ref{Assumption Minimum f}
requires the latter to be bounded away from zero. In practice, if
one suspects that at some frequencies $f\left(u,\,\omega\right)$
can be close to zero, then one can add a small number $\epsilon_{f}>0$
to the denominator of the test statistic to guarantee numerical stability.
\begin{assumption}
\label{Assumption Minimum f}$f_{-}=\min_{u\in\left[0,\,1\right],\,\omega\in\left[-\pi,\,\pi\right]}f\left(u,\,\omega\right)>0$. 
\end{assumption}
The next assumption ensures that the local spectral density estimates
are asymptotically independent when evaluated at some given frequencies
(see Theorem \ref{Theorem 5.6.4 Brillinger}). It is used to derive
the asymptotic null distribution of the double-sup test statistics
$\mathrm{S}_{\mathrm{Dmax},T}$ and $\mathrm{R}_{\mathrm{Dmax},T}$. 
\begin{assumption}
\label{Assumption Big PI}Assume that $2\omega_{j},\,\omega_{j}\pm\omega_{k}\not\equiv0\,(\mathrm{mod\,}2\pi)$
for $\omega_{j},\,\omega_{k}\in\Pi$.
\end{assumption}
\begin{condition}
\label{Condition n_T h_T}(i) The sequence $\left\{ m_{T}\right\} $
satisfies $m_{T}\rightarrow\infty$ as $T\rightarrow\infty,$ and
\begin{align}
M_{S,T}^{1/2} & m_{T}^{\theta}T^{-\theta}\left(\log\left(M_{T}\right)\right)^{1/2}+\tau_{T}^{2}\log\left(M_{T}\right)M_{S,T}^{-1}\label{eq. Condition auxiliary nT}\\
 & +M_{S,T}n_{T}^{4}\log\left(M_{T}\right)T^{-4}+M_{S,T}\left(\log\left(n_{T}\right)\right)^{2}\log\left(M_{T}\right)n_{T}^{-2}\rightarrow0;\nonumber 
\end{align}
 (ii) $b_{W,T}\rightarrow0$ such that $Tb_{W,T}\rightarrow\infty$,
 $\log\left(M_{T}\right)M_{S,T}b_{W,T}^{4}\rightarrow0$ and $\log\left(M_{T}\right)M_{S,T}(n_{T}b_{W,T})^{-1}\rightarrow0.$
\end{condition}
Part (i) imposes lower and upper bounds on the growth condition of
the sequence $\left\{ m_{T}\right\} $. The upper bound relates to
the smoothness of $A\left(u,\,\omega\right)$, the value of $\theta$
under the null hypothesis, $n_{T}$ and the number of summands $M_{S,T}$
in $\widetilde{f}_{L,r,T}^{*}\left(\omega\right)$. 

Let $\mathscr{V}$ denote a random variable with an extreme value
distribution defined by $\mathbb{P}\left(\mathscr{V}\leq v\right)=\exp($
$-\pi^{-1/2}\exp\left(-v\right)).$ 
\begin{thm}
\label{Theorem Asymptotic H0 Distrbution Smax}Let Assumption \ref{Assumption Locally Stationary}-\ref{Assumption Minimum f}
and Condition \ref{Condition n_T h_T} hold. Under $\mathcal{H}_{0}$,
$\sqrt{\log\left(M_{T}\right)}(M_{S,T}^{1/2}\mathrm{S}_{\mathrm{max},T}\left(\omega\right)$
$-\gamma_{M_{T}})\Rightarrow\mathscr{V}$ for any $\omega\in\left[-\pi,\,\pi\right]$.
\end{thm}
Theorem \ref{Theorem Asymptotic H0 Distrbution Smax} shows that the
asymptotic null distribution follows an extreme value distribution.
The derivation of the null distribution uses a (strong) invariance
principle for nonstationary processes {[}see, e.g., \nocite{wu:07}
and \citet{wu/zhou:11}{]}.  The following theorems shows that the
asymptotic null distribution of the remaining tests $\mathrm{S}_{\mathrm{Dmax},T}$,
$\mathrm{R}_{\mathrm{max},T}\left(\omega\right)$ and $\mathrm{R}_{\mathrm{Dmax},T}$
also follows an extreme value distribution, though the additional
Assumption \ref{Assumption Big PI} and the extra factor $\log\left(n'_{\omega}\right)$
are needed for $\mathrm{S}_{\mathrm{Dmax},T}$ and $\mathrm{R}_{\mathrm{Dmax},T}$.

\begin{thm}
\label{Theorem Asymptotic H0 Distrbution SDmax}Let Assumption \ref{Assumption Locally Stationary}-\ref{Assumption Big PI}
and Condition \ref{Condition n_T h_T} hold. Under $\mathcal{H}_{0}$,
$\mathrm{S}_{\mathrm{Dmax},T}\overset{}{\Rightarrow}\mathscr{V}$.
\end{thm}
\begin{thm}
\label{Theorem Asymptotic H0 Distrbution Rmax and RDmax}Let Assumption
\ref{Assumption Locally Stationary}-\ref{Assumption Minimum f} and
Condition \ref{Condition n_T h_T} hold. Under $\mathcal{H}_{0}$,
$\sqrt{\log\left(M_{T}\right)}(M_{S,T}^{1/2}\mathrm{R}_{\mathrm{max},T}\left(\omega_{k}\right)$
$-\gamma_{M_{T}})\Rightarrow\mathscr{V}$ and, in addition if Assumption
\ref{Assumption Big PI} holds, then $\mathrm{R}_{\mathrm{Dmax},T}\overset{}{\Rightarrow}\mathscr{V}$.
\end{thm}
The limiting distributions in Theorem \ref{Theorem Asymptotic H0 Distrbution Smax}-\ref{Theorem Asymptotic H0 Distrbution Rmax and RDmax}
are pivotal, so that  critical values can be obtained immediately
without the need to rely on simulations. The tests have power against
both breaks (alternative hypothesis (i)) and changes in the smoothness
(alternative hypothesis (ii)). Unfortunately, discerning between the
two types of alternative hypotheses is a quite hard technical problem.
Although knowing whether the rejection of the null hypothesis is due
to a break or a change in the smoothness would be useful, in many
cases knowing that there has been some change in the data-generating
process is sufficient to modify the estimation and/or inference strategy
to account for the change. If the change involves a break, a reasonable
approach would be to apply some sample-splitting method for estimation.
For example, in the context of long-run variance estimation the existence
of a break suggests to modify the double kernel HAC (DK-HAC) estimator
to avoid mixing two different regimes {[}cf. \citet{casini_hac}{]}.
On the other hand, for a change in the smoothness a modification of
a standard nonparametric kernel smoothing could be enough. However,
using a sample-splitting technique even when there is a change in
the smoothness would be robust to the change and would result in better
estimation and inference. In general, treating a change as a break
and applying some sample-splitting method would be technically valid
even when the rejection was due to a change in the smoothness (though
it may not be efficient). 

The property of being robust to different alternative hypotheses is
not specific to our method. This property is shared by most of the
existing structural break tests. For example, \citet{andrews:93}
showed that structural break tests also have some power against some
forms of smoothly-varying parameters. This property was seen as a
positive feature in the structural break literature. Following the
same reasoning, the property of being robust to breaks as well as
changes in the smoothness can actually be seen as a virtue of our
method.

\section{\label{Section Consistency-and-Minimax}Consistency and Minimax Optimal
Rate of Convergence}

In this section, we discuss the consistency and minimax-optimal lower
bound for the testing problem \eqref{eq (29) BJV - Testing Problem}
(i.e., case (ii)). The discussion also covers the testing problem
(i) since $\mathcal{H}_{1}^{\mathrm{B}}$ can be seen as the limiting
case of $\mathcal{H}_{1}^{\mathrm{S}}$ as $\theta'\rightarrow0.$
We assume that $X_{t,T}$ is segmented locally stationary with transfer
function $A\left(u,\,\omega\right)$ satisfying the following smoothness
properties. 
\begin{assumption}
\label{Assumption Smothness of A}(i) $\left\{ X_{t,T}\right\} $
is a mean-zero segmented locally stationary process; (ii) $A\left(u,\,\omega\right)$
is \textcolor{red}{}twice continuously differentiable in $u$ at
all $u\neq\lambda_{j}^{0}$ $(j=1,\ldots,\,m_{0}+1)$ with uniformly
bounded derivatives $\left(\partial/\partial u\right)A\left(u,\,\cdot\right)$
and $\left(\partial^{2}/\partial u^{2}\right)A\left(u,\,\cdot\right)$;
(iii) $A\left(u,\,\omega\right)$ is twice left-differentiable in
$u$ at $u=\lambda_{j}^{0}$ $(j=1,\ldots,\,m_{0}+1)$ with uniformly
bounded derivatives $\left(\partial/\partial_{-}u\right)A\left(u,\,\cdot\right)$
and $\left(\partial^{2}/\partial_{-}u^{2}\right)A\left(u,\,\cdot\right)$.
\end{assumption}
\begin{assumption}
\label{Assumption: (i) smothness in omega; (ii) g4 is continuous}(i)
$A\left(u,\,\omega\right)$ is twice differentiable in $\omega$
with uniformly bounded derivatives $\left(\partial/\partial\omega\right)$
$A\left(\cdot,\,\omega\right)$ and $\left(\partial^{2}/\partial\omega^{2}\right)A\left(\cdot,\,\omega\right)$;
(ii) $g_{4}\left(\omega_{1},\,\omega_{2},\,\omega_{3}\right)$ is
continuous in its arguments. 
\end{assumption}
We now move to the derivation of the minimax lower bound. As explained
before, we restrict attention to a strictly positive spectral density
in the frequency dimension at which the null hypotheses is violated.
That is, $f_{-}\left(\omega_{0}\right)=\inf_{u\in\left[0,\,1\right]}f\left(u,\,\omega_{0}\right)>0$.
Such restriction is not imposed on $f\left(u,\,\omega\right)$ for
$\omega\neq\omega_{0}$.
\begin{thm}
\label{Theorem 4.1 BJV}Let Assumption \ref{Assumption WZ 2011}-\ref{Taper, nT, W and K1, b1},
\ref{Assumption Smothness of A}-\ref{Assumption: (i) smothness in omega; (ii) g4 is continuous}
and $f_{-}\left(\omega_{0}\right)>0$ hold. Consider either set of
hypotheses $\{\mathcal{H}_{0},\,\mathcal{H}_{1}^{\mathrm{B}}\}$ with
$\theta'=0$ or $\{\mathcal{H}_{0},\,\ensuremath{\mathcal{H}_{1}^{\mathrm{S}}}\}$
with $0<\theta'<\theta$. Then, for 
\begin{align*}
b_{T} & \leq\left(T/\log\left(M_{T}\right)\right)^{-\frac{\theta-\theta'}{2\theta+1}}D^{-\frac{2\theta'+1}{2\theta+1}}f_{-}\left(\omega_{0}\right),
\end{align*}
we have $\lim_{T\rightarrow\infty}\inf_{\psi}\gamma_{\psi}\left(\theta,\,b_{T}\right)=1.$ 
\end{thm}
The theorem implies the need for 
\[
b_{T}^{\mathrm{opt}}>\left(T/\log\left(M_{T}\right)\right)^{-\frac{\theta-\theta'}{2\theta+1}}D^{-\frac{2\theta'+1}{2\theta+1}}f_{-}\left(\omega_{0}\right),
\]
 otherwise there cannot exist a minimax-optimal test yielding $\lim_{T\rightarrow\infty}\inf_{\psi}\gamma_{\psi}\left(\theta,\,b_{T}\right)=0.$
 Note that the lower bound does not depend on $\omega$. In Theorem
\ref{Theorem 4.2 BJV 17} we establish a corresponding upper bound.
From the lower and upper bounds we deduce the optimal rate for the
minimax distinguishable boundary. We can also derive tests based on
$b_{T}^{\mathrm{opt}}$. For example, using the test statistic \eqref{Eq. (14) Smax Stat}
for $\{\mathcal{H}_{0},\,\mathcal{H}_{1}^{\mathrm{B}}\}$ we obtain
the following test $\psi^{*}:$ $\psi^{*}(\left\{ X_{t,T}\right\} )=1$
if $\mathrm{S}_{\mathrm{max},T}\left(\omega\right)\geq2D^{*}\sqrt{\log\left(M_{T}^{*}\right)/m_{T}^{*}}$
for $\omega\in\left[-\pi,\,\pi\right]$ where $D^{*}>2$, $m_{T}^{*}=(\sqrt{\log\left(M_{T}^{*}\right)}T^{\theta}/D)^{\frac{2}{2\theta+1}}$
and $M_{T}^{*}=\left\lfloor T/m_{T}^{*}\right\rfloor $.  Hence,
in order to construct such a test we need knowledge of $\theta$ under
$\mathcal{H}_{0}$. We discuss this in Section \ref{Section Implementation}.

Next, we establish the optimal rate for minimax distinguishability.
Note that either alternatives $\mathcal{H}_{1}^{\mathrm{B}}$ or
$\mathcal{H}_{1}^{\mathrm{S}}$ allows for multiple breaks.  The
following results require further restrictions on the relation between
$n_{T}$ and $m_{T}$.
\begin{thm}
\label{Theorem 4.2 BJV 17}Let Assumption \ref{Assumption WZ 2011}-\ref{Taper, nT, W and K1, b1},
\ref{Assumption Smothness of A}-\ref{Assumption: (i) smothness in omega; (ii) g4 is continuous}
hold. Consider either alternative hypotheses $\ensuremath{\mathcal{H}_{1}^{\mathrm{B}}}$
with $\theta'=0$ and $\lambda_{j}^{0}<\lambda_{j+1}^{0}$ for $j=1,\ldots,\,m_{0}$,
or $\ensuremath{\mathcal{H}_{1}^{\mathrm{S}}}$ with $0<\theta'<\theta$.
If 
\begin{align}
\left(\sqrt{\log\left(M_{T}^{*}\right)/m_{T}^{*}}\right)^{-1}\left(\left(m_{T}^{*}/T\right)^{\theta}+\left(n_{T}/T\right)^{2}+\log\left(n_{T}\right)/n_{T}+b_{W,T}^{2}\right) & \rightarrow0,\label{Eq. Cond for Sec. 5}
\end{align}
 and 
\begin{align}
b_{T}^{*} & >\left(4D^{*}\sup_{u\in\left[0,\,1\right]}f\left(u,\,\omega_{0}\right)+2\right)^{-\frac{\theta-\theta'}{2\theta+1}}\left(T/\log\left(M_{T}\right)\right)^{-\frac{\theta-\theta'}{2\theta+1}}D^{-\frac{2\theta'+1}{2\theta+1}},\label{Eq. (34) BJV 17}
\end{align}
then $\lim_{T\rightarrow\infty}\gamma_{\psi^{*}}\left(\theta,\,b_{T}^{*}\right)=0$
and $b_{T}^{\mathrm{opt}}\propto\left(T/\log\left(M_{T}\right)\right)^{-\frac{\theta-\theta'}{2\theta+1}}.$ 
\end{thm}
 The theorem shows that a smooth change in the regularity exponent
$\theta$ cannot be distinguished from a break of magnitude smaller
than $b_{T}^{\mathrm{opt}}$  because the change from $\theta$ to
$\theta'$ has to persist for some time. This is also indicated by
the restriction $\theta'>0.$ The minimax bound is similar to the
one established by \citet{bibinger/jirak/vetter:16} for the volatility
of a It\^o semimartingale. The theorem suggests that knowledge
of the frequency $\omega_{0}$ at which the spectrum changes regularity
is irrelevant for the determination of the bound. However, we conjecture
that if the spectrum exhibits a break or smooth change of the form
discussed above simultaneously across multiple frequencies then the
lower bound may be further decreased as one can pool additional information
from inspection of the spectrum for the set of frequencies subject
to the change. The key assumption would be that the change occurs
at the same time $\lambda_{b}^{0}$ for a given set of frequencies
$\omega.$  This may be of interest for economic and financial time
series since they often exhibit a break simultaneously at high and
low frequencies. We leave this to future research. 

\section{Estimation of the Change-Points\label{Section Estimation-of-the Breaks}}

We now discuss the estimation of the break locations for the case
of discontinuities in the spectrum (i.e., $\mathcal{H}_{1}^{\mathrm{B},m_{0}}$
where $m_{0}$ is the number of breaks, recall Definition \ref{Definition Segmented-Locally-Stationary}).
The same estimator is valid for the locations of the smooth changes
as under $\mathcal{H}_{1}^{\mathrm{S}}$. For the latter case we later
provide intuitive remarks about the consistency result and the conditions
needed for it. We first consider the case of a single break (i.e.,
$\mathcal{H}_{1}^{\mathrm{B},1}$) and then present the results for
the case of multiple breaks (i.e., $\mathcal{H}_{1}^{\mathrm{B},m_{0}}$).

\subsection{Single Break Alternatives $\mathcal{H}_{1}^{\mathrm{B},1}$ }

Let  
\begin{align*}
\mathrm{D}_{r,T}\left(\omega\right)\triangleq M_{S,T}^{-1/2}\left|\sum_{j\in\mathbf{S}_{L,r}}f_{L,h,T}\left(j/T,\,\omega\right)-\sum_{j\in\mathbf{S}_{R,r}}f_{R,h,T}\left(j/T,\,\omega\right)\right|, & \qquad\omega\in\left[-\pi,\,\pi\right].
\end{align*}
where 
\begin{align*}
\mathbf{S}_{L,r} & =\{r-m_{T}+1,\,r-m_{T}+1+m_{S,T},\ldots,\,r-m_{T}+1+m_{S,T}M_{S,T}\},\\
\mathbf{S}_{R,r} & =\{r+1,\,r+1+m_{S,T},\ldots,\,r+1+m_{S,T}M_{S,T}\},
\end{align*}
and $r=2m_{T},\,3m_{T},\ldots$ with $r<\left(M_{T}-1\right)m_{T}-n_{T}$.
Note that the maximum of the statistics $\mathrm{D}_{r,T}\left(\omega\right)$
is a version of $\mathrm{S}_{\max,T}$ that does not involve the normalization.
The change-point estimator is defined as
\[
T\widehat{\lambda}_{b,T}=\underset{r=2m_{T},\,3m_{T},\ldots}{\mathrm{argmax}}\max_{\omega\in\left[-\pi,\,\pi\right]}\mathrm{D}_{r,T}\left(\omega\right).
\]
Recall that we consider the following alternative hypothesis:
\begin{align*}
\mathcal{H}_{1}^{\mathrm{B},1}: & \,\left\{ f\left(T_{b}^{0}/T,\,\omega_{0}\right)-\lim_{s\downarrow T_{b}^{0}}f\left(s/T,\,\omega_{0}\right)=\delta_{T}\neq0,\quad\omega_{0}\in\left[-\pi,\,\pi\right]\right\} .
\end{align*}
 Note that a break does not need to occur simultaneously at all frequencies
for the procedure to work. The break magnitude can be either fixed
or converge to zero as specified by the following assumption.
\begin{assumption}
\label{Assumption small shifts}$\delta_{T}=\delta\neq0$ is fixed
or $\delta_{T}\rightarrow0$ and $\delta_{T}M_{S,T}^{1/2}/\sqrt{\log\left(T\right)}\rightarrow(0,\,\infty]$.
\end{assumption}
\begin{prop}
\label{Prop 4.5 BJV Consitency Change-point}Let Assumption \ref{Assumption WZ 2011}-\ref{Taper, nT, W and K1, b1}-(i-iii),
\ref{Assumption Smothness of A} with $m_{0}=1$ and Condition \ref{Condition n_T h_T}
hold. Under $\mathcal{H}_{1}^{\mathrm{B},1},$ if $\delta_{T}$ satisfies
Assumption \ref{Assumption small shifts}, we have $\widehat{\lambda}_{b,T}-\lambda_{b}^{0}=O_{\mathbb{P}}(m_{S,T}\sqrt{M_{S,T}\log(T)}/(T\delta_{T}))$. 
\end{prop}
We compare the rate of convergence in Proposition \ref{Prop 4.5 BJV Consitency Change-point}
with that of classical change-point estimators in the piecewise constant
mean model. For fixed shifts, the latter rate of convergence is $O_{\mathbb{P}}(T^{-1})$
while for shrinking shifts it is $O_{\mathbb{P}}((T\delta_{T}^{2})^{-1})$
where $\delta_{T}\rightarrow0$ with $\delta_{T}T^{1/2-\vartheta}$
for some $\vartheta\in\left(0,\,1/2\right)$ {[}cf. \citet{yao:87}{]}.\footnote{See also \citet{verzelen/fromont/lerasle/reynaud-bouret:2020} for
recent developments on nimimax optimality for change-point estimation
in the piecewise constant mean model. They considered as a significance
of a change-point a measure that depends on both the break magnitude
and the location of the break. They named it the \textit{energy} of
the change-point. They established the uniform detection threshold
for the energy.} Unlike the classical change-point problem where the mean is piecewise
constant, our problem involves a spectrum that can vary smoothly.
The latter represents a local problem that cannot be addressed by
standard sample-splitting methods. Our method is local in nature and
so it is sub-optimal for the classical change-point problem but it
is valid for the more general case of a piecewise smooth spectrum.
Hence, for fixed shifts, the rate of convergence in our problem is
slower. The smallest break magnitude allowed by Proposition \ref{Prop 4.5 BJV Consitency Change-point}
is $\delta_{T}=O(\sqrt{\log\left(T\right)}/M_{S,T}^{1/2})$. Under
this condition the convergence rate for the classical change-point
estimator is $O_{\mathbb{P}}(M_{S,T}(T\log\left(T\right))^{-1})$
which is faster by a factor $O(m_{S,T}\sqrt{\log\left(T\right)})$
than the one suggested by Proposition \ref{Prop 4.5 BJV Consitency Change-point}.
In addition, in classical change-point setting $\delta_{T}\rightarrow0$
is allowed at a faster rate. This is obvious since in our setting
a small break can be confounded with a smooth local change.

Under the smooth alternative $\mathcal{H}_{1}^{\mathrm{S}}$ the
estimator is consistent when $\theta$-regularity is violated only
once in the sample and also when the violation occurs in a small interval
around $\lambda_{b}^{0}$ which does not exceeds $O(m_{S,T}\sqrt{M_{S,T}\log\left(T\right)}/T\delta_{T})$.
If that interval is longer then this becomes a global problem which
cannot be addressed by the estimation method considered in this section.
This also relates to the discussion in Section \ref{Section Consistency-and-Minimax}
that one cannot perfectly separate functions with $\theta$-smoothness
from functions with $\theta'$-smoothness such that $\theta'<\theta.$ 

\subsection{Multiple Breaks Alternatives $\mathcal{H}_{1}^{\mathrm{B},m_{0}}$ }

Let us assume that there are $m_{0}>1$ break points in $f\left(u,\,\omega\right)$.
Let $0<\lambda_{1}^{0}<\ldots<\lambda_{m_{0}}^{0}<1$. We consider
the following class of alternative hypotheses:
\begin{align*}
\mathcal{H}_{1}^{\mathrm{B},m_{0}}: & \,\left\{ f\left(T_{l}^{0}/T,\,\omega_{l}\right)-\lim_{s\downarrow T_{l}^{0}}f\left(s/T,\,\omega_{l}\right)=\delta_{l,T}\neq0,\quad\omega_{l}\in\left[-\pi,\,\pi\right]\,\mathrm{for\,}1\leq l\leq m_{0}\right\} .
\end{align*}
We provide a consistency result for both $m_{0}$ and the actual locations
of the breaks $\lambda_{l}^{0}$ $(1\leq l\leq m_{0})$. Let $j^{*}$
be the largest integer such that $\left(M_{T}-j^{*}\right)m_{T}\leq\left(M_{T}-1\right)m_{T}-n_{T}$
and $\mathcal{I}\subseteq\left\{ 2m_{T},\,3m_{T},\ldots,\,\left(M_{T}-j^{*}\right)m_{T}\right\} $
denote a generic index set. One can test for a break at some time
index in $\mathcal{I}$ by using the test $\psi^{*}(\{X_{r}\}_{r\in\mathcal{I}})$
based on $\max_{\omega_{k}\in\Pi'}\mathrm{S}_{\mathrm{max},T}\left(\omega_{k}\right)$
and if the test rejects one can estimate the break location using
\begin{align}
T\widehat{\lambda}_{T}\left(\mathcal{I}\right) & =\underset{r\in\mathcal{I}}{\mathrm{argmax}}\max_{\omega\in\left[-\pi,\,\pi\right]}\mathrm{D}_{r,T}\left(\omega\right).\label{eq. (38) BJV - Obj Func Estimation}
\end{align}
We can then update the set $\mathcal{I}$ by excluding a $v_{T}$-neighborhood
of $T\widehat{\lambda}_{T}$ and repeat the above steps. This is a
sequential top-down algorithm exploiting the classical idea of bisection.
However, this procedure may not be efficient. For example, consider
the first step of the algorithm in which we test for the first break;
this is associated with the largest break magnitude ($\delta_{1,T}>\delta_{l,T}$
for all $l=2,\ldots,\,m_{0}$). If the true break date $T_{1}^{0}$
falls in between two indices in $\mathcal{I}$, say $r_{1}$ and $r_{2}=r_{1}+m_{T}$,
then this does not maximize either power or precision of the location
estimate because one would need to compare two adjacent blocks exactly
separated at $T_{1}^{0}$ but $T_{1}^{0}\notin\mathcal{I}$ since
$T_{1}^{0}\in\left(r_{1},\,r_{2}\right).$ Hence, we introduce a
wild sequential top-down algorithm. 

Continuing with the above example, we draw randomly without replacement
$K\geq1$ separation points $r^{\diamond}$ from the interval $\left(r_{1},\,r_{2}\right)$
and for each separation point compute $\mathrm{D}_{r^{\diamond},T}\left(\omega\right)$
where $r^{\diamond}\in\left(r_{1},\,r_{2}\right).$ We take the maximum
value. Then, we update $\mathcal{I}$ by removing $r_{1}$ and adding
$r^{\diamond}$. We repeat this for all indices in $\mathcal{I}$.
Because the $K$ separation points are drawn randomly, there is always
some probability to pick up the separation point that guarantees the
highest power. A natural question is why not take all integers between
$r_{1}$ and $r_{2}$ and compute $\mathrm{D}_{r^{\diamond},T}\left(\omega\right)$
for each. The reason is that in applications involving high frequency
data (e.g., weakly, daily, and so on) that would be highly computationally
intensive especially with multiple breaks as one wishes to change
$m_{T}$ when searching for an additional break. This procedure exploits
idea of bisection and combines it with a wild resampling technique
similar to the one in \citet{fryzlewicz:14}. The latter is characterized
by using binary segmentation and drawing a large number of random
intervals. Here the idea of using draws of random intervals is applied
to the sequential top-down algorithm. 

We are now ready to present the algorithm. Guidance as to a suitable
choice of $K$ will be given below. Let $v_{T}\rightarrow\infty$
with $v_{T}/T\rightarrow0$ and $m_{T}/v_{T}\rightarrow0.$ Consider
the test $\psi(\left\{ X_{t,T}\right\} ,\,\mathcal{I})=1$ if $\mathrm{S}_{\mathrm{Dmax},T}\left(\mathcal{I}\right)\geq2D^{*}\sqrt{\log\left(M_{T}^{*}\right)/m_{T}^{*}}$
where 
\begin{align*}
\mathrm{S}_{\mathrm{Dmax},T}\left(\mathcal{I}\right)\triangleq\max_{r\in\mathcal{I}}\max_{\omega_{k}\in\Pi'}\left|\frac{\sum_{j\in\mathbf{S}_{L,r}}f_{L,h,T}\left(j/T,\,\omega_{k}\right)-\sum_{j\in\mathbf{S}_{R,r}}f_{R,h,T}\left(j/T,\,\omega_{k}\right)}{\widehat{\sigma}_{L,r}\left(\omega_{k}\right)}\right| & ,
\end{align*}
with $D^{*}$, $m_{T}^{*}$ and $M_{T}^{*}$ as defined in Section
\ref{Section Consistency-and-Minimax}.  

\begin{lyxalgorithm}
\label{Algorithm 1}Set $\widehat{\mathcal{I}}=\{2m_{T},\,3m_{T},\ldots,\,\left(M_{T}-j^{*}\right)m_{T}\}$
and $\widehat{\mathcal{T}}=\emptyset$. \\
(1) For  $r\in\widehat{\mathcal{I}}\backslash\{2m_{T}\}$ uniformly
draw (without replacement) $K\in\{1,\ldots,\,m_{T}\}$ points $r_{k}^{\diamond}$
from $\mathbf{I}\left(r\right)=\{r-m_{T}+1,\ldots,r\}$ and compute
$\overline{r}^{\diamond}=\arg\max_{k=1,\ldots,\,K}\max_{\omega\in\left[-\pi,\,\pi\right]}\mathrm{D}_{r_{k}^{\diamond},T}\left(\omega\right)$;
set $\widehat{\mathcal{I}}=(\widehat{\mathcal{I}}\backslash\left\{ r\right\} )\cup\{\overline{r}^{\diamond}\}$.
\\
(2) If $\psi(\left\{ X_{t,T}\right\} ,\,\widehat{\mathcal{I}})=0$
return $\widehat{\mathcal{T}}=\emptyset$. Otherwise proceed with
step (3).\\
(3) Estimate the change-point $T\widehat{\lambda}_{T}(\widehat{\mathcal{I}})$
via \eqref{eq. (38) BJV - Obj Func Estimation} using $\widehat{\mathcal{I}}$.
\\
 (4) Set $\widehat{\mathcal{I}}=\widehat{\mathcal{I}}\backslash\{T\widehat{\lambda}_{T}(\mathcal{\widehat{I}})-v_{T},\ldots,\,T\widehat{\lambda}_{T}(\mathcal{\widehat{I}})+v_{T}\}$
and $\widehat{\mathcal{T}}=\widehat{\mathcal{T}}\cup\{T\widehat{\lambda}_{T}(\mathcal{\widehat{I}})\}$.
Return to step (1). 
\end{lyxalgorithm}
Finally, arrange the estimated change-points $\widehat{\lambda}_{l,T}$
in $\widehat{\mathcal{T}}$ in chronological order and use the symbol
$\left|\mathcal{S}\right|$ for the cardinality of a set $\mathcal{S}$.
To each $\widehat{\lambda}_{l,T}$ the procedure can return the frequency
$\widehat{\omega}_{l}$ at which the break is found. 
\begin{assumption}
\label{Assumption small shift Multiple breaks}$\delta_{l,T}=\delta_{l}\neq0$
is fixed or $\delta_{l,T}\rightarrow0$ with $\inf_{1\leq l\leq m_{0}}\delta_{l,T}\geq2D^{*}M_{S,T}^{-1/2}(\log(T))^{2/3}$.
For $\nu_{T}\rightarrow\infty$ with $\nu_{T}=o\left(T/v_{T}\right),$
it holds that $\inf_{1\leq l\leq m_{0}-1}|\lambda_{l+1}^{0}-\lambda_{l}^{0}|\geq\nu_{T}^{-1}.$
\end{assumption}
Assumption \ref{Assumption small shift Multiple breaks} allows for
shrinking shifts and a possibly growing number of change-points as
long as $m_{0}/\nu_{T}\rightarrow0.$ The following proposition presents
the consistency result for the number of change-points $m_{0}$ and
for the change-point locations $\lambda_{l}^{0}$ $\left(l=1,\ldots,\,m_{0}\right)$,
and the rate of convergence of their estimates. 
\begin{prop}
\label{Prop 4.9 BJV Consitency Change-point}Let Assumption \ref{Assumption WZ 2011}-\ref{Taper, nT, W and K1, b1},
\ref{Assumption Smothness of A} with $m_{0}=1$ and Condition \ref{Condition n_T h_T}
hold. Then, under $\mathcal{H}_{1}^{\mathrm{B},m_{0}}$ we have (i)
$\mathbb{P}(|\widehat{\mathcal{T}}-m_{0}|>\epsilon_{2})\rightarrow0$
for any $\epsilon_{2}>0$ and $\sup_{1\leq l\leq m_{0}}|\widehat{\lambda}_{l,T}-\lambda_{l}^{0}|=o_{\mathbb{P}}\left(1\right)$,
and (ii) $\sup_{1\leq l\leq m_{0}}|\widehat{\lambda}_{l,T}-\lambda_{l}^{0}|=O_{\mathbb{P}}(m_{S,T}\sqrt{M_{S,T}\log\left(T\right)}/(T\inf_{1\leq l\leq m_{0}}\delta_{l,T}))$.
 Furthermore, if $K=O(a_{T}m_{T})$ with $a_{T}\in(0,\,1]$ such
that $a_{T}\rightarrow1,$ then the breaks are detected in decreasing
order of magnitude. 
\end{prop}
The number of draws $K$ may be fixed or increase with the sample
size. However, the algorithm can return the change-point dates in
decreasing order of the break magnitudes only if $K$ is sufficiently
large. Note that at each loop of the algorithm it is not possible
to know to which $\lambda_{l}^{0}$ $\left(l=1,\ldots,\,m_{0}\right)$
the estimate $\widehat{\lambda}_{l,T}$ is consistent for. Only after
all breaks are detected and we rearrange the estimated change-points
in $\widehat{\mathcal{T}}$ in chronological order, we can learn such
information. The same procedure can be applied for the case of multiple
smooth local changes, though the notation becomes cumbersome and so
we omit it.

\section{Implementation\label{Section Implementation}}

In this section we explain how to choose the tuning parameters. The
choice of $m_{T}$, $n_{T}$ and $b_{W,T}$ could be based on a mean-squared
error (MSE) criterion or cross-validation exploiting results derived
for locally stationary series {[}e.g., data-dependent methods for
bandwidths in the context of locally stationary processes were investigated
by, among others, \citet{casini_hac}, \citet{dahlhaus:12}, \citet{Dahlhaus/Giraitis:98}
and \citet{ritcher/dahlhaus:19}{]}. The optimal amount of smoothing
depends on the regularity exponent $\theta$, on the  boundness of
the moments and on the extent of the dependence in $\left\{ X_{t,T}\right\} $.
 Here we choose the order of the bandwidths, neglecting the constants,
by following the restrictions in Condition \ref{Condition n_T h_T}.
In particular, we choose the largest possible values allowed by Condition
\ref{Condition n_T h_T} in order to ensure the highest possible power.
We conduct a sensitivity analysis based on simulations in the supplement.
We relegate to future work a more detailed analysis of data-dependent
methods for this problem with multiple smoothing directions. 

For spectral densities satisfying Lipschitz continuity, $\theta=1$
so that $m_{T}\propto T^{2/3-\epsilon}$ while for $\theta=1/2$ we
have $m_{T}\propto T^{1/2-\epsilon}$ where in both cases $\epsilon>0.$
 In applied work, it is common to work under stationarity ($\theta>1$)
or local stationarity with Lipschitz smoothness ($\theta=1$). Hence,
we use the bandwidths corresponding to $\theta=1$ which works for
both cases. Of course, if one has prior knowledge about the smoothness
of the parameters of the data-generating process, one can choose a
suitable $\theta$. Assuming $q=8$ and $\gamma$ large enough we
have $\tau_{T}\propto T^{1/4}$ and so values that satisfy Condition
\ref{Condition n_T h_T} are $m_{T}=T^{0.66}$, $n_{T}=T^{0.62}$
and $b_{W,T}=n_{T}^{-1/6}$.  The scaling is normalized to 1, as
our simulations show this to provide good finite-sample properties,
see Section \ref{Section Monte Carlo}.

As for the tapering function $h\left(\cdot\right)$ and weight function
$W\left(\cdot\right)$ we use a rectangular taper (i.e., $h\left(u\right)=1$
for all $u$) and a rectangular kernel. The rectangular kernel is
known as the Daniell kernel with parameter $n_{T}b_{W,T}$ (it is
a centered moving average which creates a smoothed value at time $Tu$
by averaging all values between $Tu-n_{T}b_{W,T}$ and $Tu+n_{T}b_{W,T}$).
These are the simplest choices for $h\left(\cdot\right)$ and $W\left(\cdot\right)$.
As for the bandwidth and kernel of the estimator $\widehat{\sigma}_{L,r}\left(\omega\right),$
we follow the results in \citeauthor{casini_comment_andrews91} (\citeyear{casini_comment_andrews91},
\citeyear{casini_hac}) that suggest $b_{1,T}=M_{S,T}^{-1/3}.$ This
corresponds to the MSE-optimal bandwidth when $K_{1}\left(\cdot\right)$
is  the Bartlett kernel. For the choice of the number and values
of the frequencies, the theory does not suggest particular values.
Thus, we tried several values $n_{\omega}=15,\,11,\,7,\,5.$ Our default
choice is $n_{\omega}=7$. A sensitivity analysis suggests that different
choices for $n_{\omega}$ lead to negligible differences in the results.
For the selection of the set of frequencies we use the function linspace
that generates a linearly spaced sequence, e.g., in \textsc{Matlab}
we used the command $\mathrm{\mathsf{linspace(-\pi+1e-3,\,0,\,(\mathit{n}_{\omega}+1)/2)}.}$

The regularity exponent $\theta$ also affects the test $\psi(\left\{ X_{t,T}\right\} ,\,\mathcal{I})$
in Algorithm \ref{Algorithm 1}. It is possible to get an estimate
of $\theta$ under the null as follows. Compute $\mathrm{S}_{\mathrm{Dmax},T}$
where the maximum is taken among the indices of the blocks such that
the null hypothesis is not violated and label it $s_{\mathrm{Dmax}}^{*}$.
 Solve $s_{\mathrm{Dmax}}^{*}=2\sqrt{\log\left(M_{T}^{*}\right)/m_{T}^{*}}$
for $\theta$, where recall that $m_{T}^{*}$ and $M_{T}^{*}$ depend
on $\theta.$ This yields a preliminary estimate of $\theta$ which
can then be used for the test $\psi(\left\{ X_{t,T}\right\} ,\,\mathcal{I})$.
Similarly, $b_{T}^{\mathrm{opt}}$ depends on $\theta$ and $\theta'$.
Using the same approach, for a given $\theta'$ one can solve $s_{\mathrm{Dmax}}^{*}=b_{T}^{\mathrm{opt}}$
for $\theta$ as function of $\theta'.$ If one is interested in the
alternative $\mathcal{H}_{1}^{\mathrm{B}}$, we have $\theta'=0$
and so this immediately yields an estimate for $\theta$. If one is
interested in the alternative $\mathcal{H}_{1}^{\mathrm{S}}$, then
one can try a few values of $\theta'$ in the range $\left(0,\,\theta\right)$.
However, note that in order to use Algorithm \ref{Algorithm 1} only
$\theta$ is needed. The knowledge of $\theta'$ under $\mathcal{H}_{1}^{\mathrm{S}}$
is only needed to obtain $b_{T}^{\mathrm{opt}}.$ 

We set $v_{T}=T^{0.666}$ which satisfies $m_{T}/v_{T}\rightarrow0.$
Our default recommendation is $K=10.$ Our simulations with different
data-generating processes and sample sizes show that this choice strikes
a good balance between the precision of the change-point estimates
and computing time. For $T>1000$, we recommend setting $K=\left\lfloor m_{T}/3\right\rfloor $. 

The test statistics $\mathrm{S}_{\mathrm{max},T}\left(\omega\right)$
and $\mathrm{R}_{\mathrm{max},T}\left(\omega\right)$ depend on $\omega.$
The choice of $\omega$ is, of course, important as it involves different
frequency components and hence different periodicities. If the user
does not have a priori knowledge about the frequency at which the
spectrum has a change-point, our recommendation is to run the tests
for multiple values of $\omega\in\left[0,\,\pi\right]$. Even if the
change-point occurs at some $\omega_{0}$ and one selects a value
of $\omega$ close but not equal to $\omega_{0}$ the tests are still
able reject the null hypothesis given the differentiability of $f\left(u,\,\omega\right)$.
Thus, one can select a few values of $\omega$ evenly spread on $\left[0,\,\pi\right]$.

\section{\label{Section Monte Carlo}Small-Sample Evaluations}

In this section, we conduct a Monte Carlo analysis to evaluate the
properties of the proposed methods. We first discuss the detection
of the change-points and then their localization. We investigate different
types of changes and consider the test statistics $\mathrm{S}_{\max,T}\left(\omega\right),$
$\mathrm{S}_{\mathrm{Dmax},T}$, $\mathrm{R}_{\max,T}\left(\omega\right)$,
$\mathrm{R}_{\mathrm{Dmax},T}$ proposed here and the test statistic
$\widehat{D}$ proposed by \citet{last/shumway:08}. The latter is
included for comparison since it applies to the same problems. We
consider the following data-generating processes where in all models
the innovation $e_{t}$ is a Gaussian white noise $e_{t}\sim\mathrm{i.i.d.}\,\mathscr{N}\left(0,\,1\right)$.
Models M1 involves a stationary AR(1) process $X_{t}=\rho X_{t-1}+e_{t}$
with $\rho=0.3$ and 0.6, while M2 involves a locally stationary AR(1)
$X_{t}=\rho\left(t/T\right)X_{t-1}+e_{t}$ where $\rho\left(t/T\right)=0.4\cos\left(0.8-\cos\left(2t/T\right)\right).$
Note that $\rho\left(t/T\right)$ varies smoothly from 0.1389 to 0.3920.
Model M1 and M2 are used to verify the finite-sample size of the tests.
We verify the power in models M3 and M4 using the specification in
model M1 and M2, respectively, for the first regime and consider two
additional regimes with different specifications. Hence, two breaks
are present. In model M3, 
\begin{align*}
X_{t} & =\begin{cases}
0.3X_{t-1}+e_{t}, & 1\leq t\leq\left\lfloor T\lambda_{1}^{0}\right\rfloor \\
0.6X_{t-1}+0.7e_{t}, & \left\lfloor T\lambda_{1}^{0}\right\rfloor +1\leq t\leq\left\lfloor T\lambda_{2}^{0}\right\rfloor \\
0.6X_{t-1}+e_{t}, & \left\lfloor T\lambda_{2}^{0}\right\rfloor +1\leq t\leq T
\end{cases},
\end{align*}
 while, for model M4 
\begin{align*}
X_{t} & =\begin{cases}
\rho\left(t/T\right)X_{t-1}+0.7e_{t}, & 1\leq t\leq\left\lfloor T\lambda_{1}^{0}\right\rfloor \\
0.8X_{t-1}+e_{t}, & \left\lfloor T\lambda_{1}^{0}\right\rfloor +1\leq t\leq\left\lfloor T\lambda_{2}^{0}\right\rfloor \\
\rho\left(t/T\right)X_{t-1}+0.7e_{t}, & \left\lfloor T\lambda_{2}^{0}\right\rfloor +1\leq t\leq T
\end{cases},
\end{align*}
where $\rho\left(t/T\right)$ is as in model M2. In model M3, the
second regime involves higher serial dependence while in the third
regime the variance doubles relative to the second regime. In model
M4, the second regime involves a stationary autoregressive process
with strong serial dependence while in the third regime $X_{t}$ assumes
the same dynamics as in the first regime. Models M3-M4 feature alternative
hypotheses in the forms of breaks in the spectrum. 

We consider the alternative hypothesis of more rough variation without
signifying a break (i.e., $\mathcal{H}_{1}^{\mathrm{S}}$ defined
in Section \ref{Section Consistency-and-Minimax}) in model M5 given
by $X_{t}=\sigma\left(t/T\right)e_{t}$ where $\sigma^{2}\left(t/T\right)=\max\{1.5,\,\overline{\sigma}^{2}+\cos\left(1+\cos\left(10t/T\right)\right)\}$
with $\overline{\sigma}^{2}=1.$ Note that even though $\sigma^{2}\left(\cdot\right)$
is locally stationary, the degree of smoothness alternates throughout
the sample. It starts from $\sigma^{2}\left(\cdot\right)=1.5$ and
maintains this value for some time, then within a short period it
increases slowly to $\sigma^{2}\left(\cdot\right)=2$ and falls slowly
back to $\sigma^{2}\left(\cdot\right)=1.5$. It keeps this value until
the final part of the sample where it increases slowly to $\sigma^{2}\left(\cdot\right)=2$
for a short period. Thus, $\sigma^{2}\left(\cdot\right)$ alternates
between periods where it is constant (i.e., $\theta>1$) and periods
where it becomes non-constant but less smooth (i.e., $\theta=1$).
Importantly, no break occurs; only a change in the smoothness as specified
in $\ensuremath{\mathcal{H}_{1}^{\mathrm{S}}}$. In unreported simulations
we also considered the case where $\theta$ changes from Lipschitz
continuity (i.e., $\theta=1$) to the continuity-path of Wiener processes
(i.e., $\theta\thickapprox1/2$) with results that are similar to
those reported here. For the test statistic $\widehat{D}$ of \citet{last/shumway:08},
we obtain the critical value by simulations. As suggested by the authors
we compute the finite-sample distribution of $\widehat{D}$ by simulating
a white noise under the null hypotheses with a sample size $T=1000$
and then obtain the critical value. We consider the three sample sizes
$T=250,\,500$ and 1000. The significance level is $\alpha=0.05$.
For the test statistics $\mathrm{S}_{\mathrm{max},T}\left(\omega\right)$
and $\mathrm{R}_{\mathrm{max},T}\left(\omega\right)$, we use as a
default value $\omega=0$ given that the interest is often in low
frequency analysis. We set $\lambda_{1}^{0}=0.33$ and $\lambda_{2}^{0}=0.66$
throughout. The number of simulations is 5,000 for all cases.

The results are reported in Table \ref{Table SM1-M2}-\ref{Table Power M3-M5}.
We first discuss the size of the tests. The tests proposed in this
paper have good empirical size for both models and all sample sizes.
The test statistics $\mathrm{S}_{\mathrm{Dmax},T}$ and $\mathrm{R}_{\mathrm{Dmax},T}$
are slightly undersized for $T=250$ but their empirical size improves
for $T=500$ and 1000. The test statistics $\mathrm{S}_{\mathrm{max},T}$
and $\mathrm{R}_{\mathrm{max},T}$ share accurate empirical sizes
in all cases. In contrast, the test statistic $\widehat{D}$ of \citet{last/shumway:08}
is largely oversized for $T=250$ and 500. For $T=1000$ it works
better but it is still oversized. This means that the finite-sample
distribution of $\widehat{D}$ has high variance and changes substantially
across different sample sizes. Since the simulated critical value
is obtained with a sample size $T=1000$ it works better for this
sample size than for the others for which the size control is poor. 

Turning to the power of the tests, we note that it is not fair to
compare the proposed tests with the test $\widehat{D}$ when $T=250$
and 500 since the latter is largely oversized in those cases. In model
M3, all the proposed tests have good power which increases with the
sample size. The tests $\mathrm{S}_{\mathrm{Dmax},T}$ and $\mathrm{R}_{\mathrm{max},T}$
have the highest power, followed by $\mathrm{S}_{\mathrm{max},T}$
and lastly $\mathrm{R}_{\mathrm{Dmax},T}$. The power differences
are not large except those involving $\mathrm{R}_{\mathrm{Dmax},T}$
for $T=250$ which has substantially lower power. It is important
to note that for $T=1000$ the proposed tests have higher power than
the test $\widehat{D}$ of \citet{last/shumway:08} even though the
latter is oversized. For $T=250$ and 500, where the $\widehat{D}$
test is largely oversized, the proposed tests only have slightly lower
power. This confirms that the proposed tests have very good power.
Similar comments apply to model M4.

Model M5 involves changes in the smoothness without involving a break.
This constitutes a more challenging alternative hypothesis, and as
expected, the power for each test is lower than in models M3-M4. The
test with the highest power is $\mathrm{S}_{\mathrm{Dmax},T}$. For
$T=1000$, the test with the lowest power is $\widehat{D}$ (for $T=250$
the $\widehat{D}$ test has higher power again due to its oversize
problem). Overall, the results show that the proposed tests have accurate
empirical size even for small sample sizes and have good power against
different forms of breaks or changes in ``roughness''. 

Next, we consider the estimation of the number of change-points ($m_{0}$)
and their locations. We consider the following two models, both with
$m_{0}=2$. The model M6 is given by
\begin{align*}
X_{t} & =\begin{cases}
0.7e_{t}, & 1\leq t\leq\left\lfloor T\lambda_{1}^{0}\right\rfloor \\
0.6X_{t-1}+0.7e_{t}, & \left\lfloor T\lambda_{1}^{0}\right\rfloor +1\leq t\leq\left\lfloor T\lambda_{2}^{0}\right\rfloor \\
0.6X_{t-1}+e_{t}, & \left\lfloor T\lambda_{2}^{0}\right\rfloor +1\leq t\leq T
\end{cases},
\end{align*}
while model M7 is the same as model M4. We set $\lambda_{1}^{0}=0.33$
and $\lambda_{2}^{0}=0.66$ and $T=1000$ throughout. Table \ref{Table Estimation}
reports summary statistics for $\widehat{m}-m_{0}$. It displays the
percentage of times with $\widehat{m}=m_{0}$, the median, and the
25\% and 75\% quantile of the distribution of $\widehat{m}$. We only
consider Algorithm \ref{Algorithm 1}. We do not report the results
for the corresponding procedure of \citet{last/shumway:08} because
it is based on $\widehat{D}$ which is oversized and so it finds many
more breaks than $m_{0}$. Table \ref{Table Estimation} shows that
$\widehat{m}=m_{0}$ occurs for about 85\% of the simulations with
model M6 and about 80\% with model M7. This suggests that Algorithm
\ref{Algorithm 1} is quite precise. As expected it performs better
in model M6 since the specification of the alternative is farther
from the null. The quantiles of the empirical distribution also suggest
that the change-point estimates $\widehat{T}_{1}$ and $\widehat{T}_{2}$
are accurate. For example the median is very close to the their respective
true value $T_{1}^{0}=333$ and $T_{2}^{0}=666$. Similar conclusions
arise from different models and sample sizes, in unreported simulations.

\section{\label{Section Empirical Application}Empirical Application}

We demonstrate how to use our change-point methods for studying the
causal effects of monetary policy. A fast growing literature in macroeconomics
uses high-frequency data to identify the effects of monetary policy
on the real economy (i.e., money non-neutrality). The identifying
assumption used by \citet{nakamura/steinsson:2018} is that the volatility
of the daily change in the nominal 2-year Treasury yields, say $\Delta i_{t}$,
is higher during days when the Federal Open Market Committee (FOMC)
meets to make monetary policy announcements relative to regular Tuesdays
and Wednesdays with no announcement. Let $\eta_{t}$ denote a pure
monetary shock and suppose that the policy instrument $\Delta i_{t}$,
which is observed in the data, is governed by both monetary and non-monetary
shocks:
\begin{align*}
\Delta i_{t} & =\mu_{i}+\eta_{t}+\varepsilon_{t},
\end{align*}
where $\varepsilon_{t}$ is a function of all other shocks that affect
$\Delta i_{t}$ and $\mu_{i}$ is a constant. We normalize the impact
of $\eta_{t}$ and $\varepsilon_{t}$ on $\Delta i_{t}$ to one. $\Delta i_{t}$
is a measure of the monetary policy news revealed in the FOMC announcement.
The idea is that changes in the policy instrument during days when
there is a FOMC announcement are dominated by the information about
future monetary policy contained in the announcement. Let $\Delta s_{t}$
denote the change in the outcome variable which is the yield on a
five year zero-coupon Treasury bond. We wish to estimate the effects
of the monetary shock $\eta_{t}$ on the outcome variable $\Delta s_{t}$.
The latter is also affected by both the monetary and non-monetary
shocks: 
\begin{align*}
\Delta s_{t} & =\mu_{s}+\beta\eta_{t}+\alpha\varepsilon_{t},
\end{align*}
where $\mu_{s}$ is a constant and $\alpha$ and $\beta$ are two
parameters. The parameter of interest is $\beta$ which represents
the impact of the pure monetary shock $\eta_{t}$ on $\Delta s_{t}$
relative to its impact on $\Delta i_{t}$. 

The identifying assumption is that the variance of monetary shocks
increases during days of FOMC announcements, while the variances of
the other shocks are unchanged. Let $T_{P}$ denote the number of
days containing a FOMC announcement, and let $T_{C}$ the number of
days with no such announcements. The subscript ``$P$'' in $T_{P}$
refers to the \textquotedblleft policy or treatment\textquotedblright{}
sample while the subscript ``$C$'' in $T_{C}$ refers to the \textquotedblleft control\textquotedblright{}
sample. The days in the control sample are comparable on other dimensions
since they are all Tuesdays and Wednesdays with no FOMC meeting. The
identifying restriction can then be written as 
\begin{align}
\sigma_{\eta,P} & >\sigma_{\eta,C}\qquad\mathrm{and}\qquad\sigma_{\varepsilon,P}=\sigma_{\varepsilon,C},\label{Eq. Vol_P >Vol_C}
\end{align}
 where $\sigma_{a,i}$ is the volatility of variable $a=\eta,\,\varepsilon$
in the sample $i=P,\,C$. One can show that 
\begin{align}
\beta & =\frac{\mathrm{Cov}_{P}\left(\Delta i_{t},\,\Delta s_{t}\right)-\mathrm{Cov}_{C}\left(\Delta i_{t},\,\Delta s_{t}\right)}{\mathrm{Var}_{P}\left(\Delta i_{t}\right)-\mathrm{Var}_{C}\left(\Delta i_{t}\right)},\label{Eq. beta}
\end{align}
 where $\mathrm{Cov}_{i}(\cdot,\,\cdot)$ (resp. $\mathrm{Var}_{i}(\cdot,\,\cdot)$)
denotes the population covariance (resp. variance) in the $i=P,\,C$
sample. The parameter $\beta$ can be identified only if $\mathrm{Var}_{P}\left(\Delta i_{t}\right)\neq\mathrm{Var}_{C}\left(\Delta i_{t}\right)$.
If the volatility of the policy instrument $\Delta i_{t}$ does not
change across treatment and control samples, then $\beta$ is not
identified. If the change in volatility is small, then $\beta$ is
weakly identified. 

Subsequent developments in the literature employed robust weak identification
tests to argue that $\beta$ is not strongly identified when using
daily data. In contrast, $\beta$ can be strongly identified if one
uses ultra high-frequency data based on a 30 minute window around
the announcement time. For the case of daily data, we show that the
volatility of $\Delta i_{t}$ in the control sample varies substantially
over time and that several change-points can be detected. Thus, the
regimes in the control sample where the volatility is high contribute
to an average (over the full control sample) volatility that approaches
the average volatility in the treatment sample, thereby violating
$\mathrm{Var}_{P}\left(\Delta i_{t}\right)\neq\mathrm{Var}_{C}\left(\Delta i_{t}\right)$.
This implies that $\beta$ may be weakly identified and its estimates
may be imprecise which supports the recent evidence in the literature. 

We obtain the data from Emi Nakamura's webpage. The sample of ``treatment''
days is all regularly scheduled FOMC meeting day from 1/1/2000 to
3/19/2014. The sample of ``control'' days is all Tuesdays and Wednesdays
that are not FOMC meeting days from 1/1/2000 to 12/31/2012. In both
the treatment and control samples, the second half of 2008, the first
half of 2009 and a 10 day period after 9/11/2001 are dropped in \citet{nakamura/steinsson:2018}.
We follow the same practice. 

\begin{center}
\begin{figure}[H]
\includegraphics[width=18cm,height=9cm]{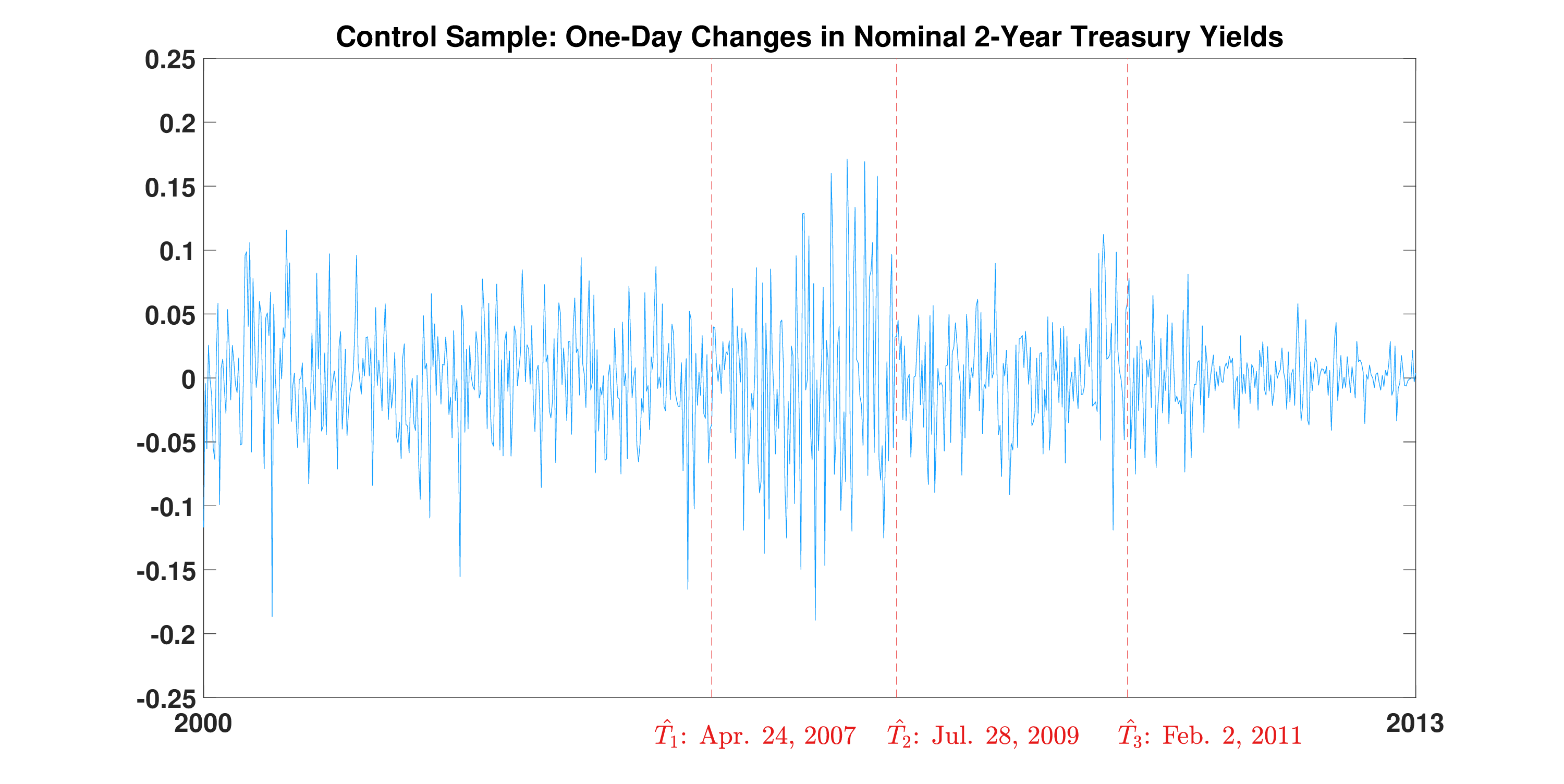}

{\footnotesize{}\caption{{\scriptsize{}\label{Fig1} }{\footnotesize{}Plot of one-day changes
in the nominal Treasury yields ($\Delta i_{t}$) in the control sample.
The sample size is $T_{C}=762$ which corresponds to all Tuesdays
and Wednesdays that are not FOMC meeting days from 1/1/2000 to 12/31/2012.
Following \citet{nakamura/steinsson:2018} we drop the second half
of 2008, the first half of 2009 and a 10 day period after 9/11/2001.
The red dashed lines are change-point dates estimated using Algorithm
\ref{Algorithm 1}.}}
}{\footnotesize\par}
\end{figure}
\end{center}

The plot of $\left\{ \Delta i_{t}\right\} $ for the control sample
is reported in Figure \ref{Fig1}. The series displays substantial
changes in volatility and some changes in persistence. The tests $\mathrm{R}_{\mathrm{max},T}$
and $\mathrm{R}_{\mathrm{Dmax},T}$ strongly reject the null hypothesis
of no change-points in the spectrum. Algorithm \ref{Algorithm 1}
detects three change-points. 

The first change-point date, denoted $\widehat{T}_{1}$, corresponds
to April 24, 2007. Thus, the first regime $[1,\,\widehat{T}_{1}]$
refers to the period prior to the beginning of the 2007-09 financial
crisis. It is evident that a large change in volatility and possibly
a change in persistence occurred. We note that the change in volatility
 does not occur abruptly, it is rather gradual. This means that
the volatility path changes gradually and possibly becomes more rough.
This supports the usefulness of our hypothesis testing framework with
local stationarity under the null hypothesis and with changes in the
smoothness of the parameters that govern the data-generating process
under the alternative hypothesis. The first change-point $\widehat{T}_{1}$
is associated to the volatility path becoming more rough with the
level of the volatility increasing gradually. 

The second change-point date, denoted $\widehat{T}_{2}$, corresponds
to July 28, 2009. Thus, the second regime $[\widehat{T}_{1}+1,\,\widehat{T}_{2}]$
includes roughly the 2007-09 financial crisis. In this regime the
volatility of the series is remarkably high. After $\widehat{T}_{2}$
the series shares a pattern similar to that in the first regime $[1,\,\widehat{T}_{1}]$,
both in terms of persistence and volatility. The second change-point
date $\widehat{T}_{2}$ is associated to an abrupt fall in volatility. 

The third change-point date, denoted $\widehat{T}_{3}$, corresponds
to February 2, 2011. The third regime $[\widehat{T}_{2}+1,\,\widehat{T}_{3}]$
corresponds to a zero lower bound (ZLB) period when the FOMC announced
the conduct of unconventional monetary policies to stimulate the economy
after the crisis. The ZLB refers to a situation in which the short-term
nominal interest rate is at or near zero, limiting the central bank's
capacity to stimulate growth. Unconventional monetary policies refer
to large-scale asset purchases and active use of communication (i.e.,
``forward guidance'') to shape expectations about future monetary
policies, which can mitigate the limitations imposed by the ZLB {[}see,
e.g., \citet{swanson:2021}{]}.

The last regime, $[\widehat{T}_{3}+1,\,T_{C}]$, corresponds to the
period when the economy was witnessing the first effects of the expansive
monetary policy and of the stability of the unconventional monetary
policies that were introduced previously. This is a regime where initially
the economy started the recovery and then reached stable economic
growth. In this regime, the volatility level of the series is the
lowest of the sample. 

Overall, the first change-point date, $\widehat{T}_{1}$, corresponds
to a change in the smoothness while the second change-point date,
$\widehat{T}_{2}$, corresponds to an abrupt break. For the third
change-point date, it is more difficult to tell from the plot whether
this corresponds to an abrupt or smooth break.

We now discuss how the results about the change-points can be useful
for the identification issue based on \eqref{Eq. Vol_P >Vol_C}-\eqref{Eq. beta}.
One should think of $\mathrm{Var}_{C}\left(\Delta i_{t}\right)$ as
the average variance of $\Delta i_{t}$ in the control sample. Our
results show that there is significant time variation in $\mathrm{Var}_{C}\left(\Delta i_{t}\right)$.
In the second regime, $[\widehat{T}_{1}+1,\,\widehat{T}_{2}]$, the
volatility is the highest of the sample. During this period, it is
lower but very close to the average volatility of $\Delta i_{t}$
in the treatment sample, $\mathrm{Var}_{P}\left(\Delta i_{t}\right)$.\footnote{We refer to \citet{nakamura/steinsson:2018} for details about the
series in the treatment sample. We do not test for change-points in
the treatment sample because the sample size in the treatment sample
is relatively small ($T_{P}=74)$ and so we treat it as a single regime.} This contributes to make $\mathrm{Var}_{P}\left(\Delta i_{t}\right)-\mathrm{Var}_{C}\left(\Delta i_{t}\right)$
(over the full sample) closer to zero which then would lead to weak
identification. In fact, \citet{nakamura/steinsson:2018} found the
estimate of $\beta$ to be imprecise and not meaningful from an economic
standpoint. Furthermore, it was quite different from that obtained
using 30 minute data instead of daily data. Our change-point analysis
is useful because it suggests which periods or sub-samples contribute
to this identification problem and which sub-samples can be used to
obtain consistent estimates for $\beta.$ 

\section{\label{Section Conclusions}Conclusions}

We develop a theoretical framework for inference about the smoothness
of the spectral density over time. We provide frequency-domain statistical
tests for the detection of discontinuities in the spectrum of a segmented
locally stationary time series and for changes in the regularity exponent
of the spectral density over time. The null distribution of the test
follows an extreme value distribution. We rely on the theory on minimax-optimal
testing developed by \citet{ingster:93}. We determine the optimal
rate for the minimax distinguishable boundary, i.e., the minimum break
magnitude such that we are still able to uniformly control type I
and type II errors. We propose a novel procedure to estimate the change-points
based on a wild sequential top-down algorithm and show its consistency
under shrinking shifts and possibly growing number of change-points.
The advantage of using frequency-domain methods to detect change-points
is that it does not require to make assumptions about the data-generating
process under the null hypothesis beyond the fact that the spectrum
is differentiable and bounded. Furthermore, the method allows for
a broader range of alternative hypotheses compared to time-domain
methods which usually have power against a limited set of alternatives.
Overall, our simulations and empirical results show the usefulness
of our method.

\newpage{}

\bibliographystyle{elsarticle-harv}
\bibliography{References_JoE}
\addcontentsline{toc}{section}{References}

\section{Appendix}

\subsection{Tables}

\begin{table}[H]
\caption{\label{Table SM1-M2}Empirical small-sample size for models M1-M2}

\begin{centering}
\begin{tabular}{lccc}
\hline 
 & \multicolumn{3}{c}{Model M1}\tabularnewline
$\alpha=0.05$ & $T=250$ & $T=500$ & $T=1000$\tabularnewline
\hline 
\hline 
$\mathrm{S}_{\max,T}\left(0\right)$ & 0.039 & 0.043 & 0.053\tabularnewline
$\mathrm{S}_{\mathrm{Dmax},T}$ & 0.029 & 0.049 & 0.047\tabularnewline
$\mathrm{R}_{\mathrm{max},T}\left(0\right)$ & 0.040 & 0.054 & 0.042\tabularnewline
$\mathrm{R}_{\mathrm{Dmax},T}$ & 0.025 & 0.032 & 0.038\tabularnewline
$\widehat{D}$ statistic & 0.581 & 0.471 & 0.068\tabularnewline
\hline 
 & \multicolumn{3}{c}{Model M2}\tabularnewline
 & $T=250$ & $T=500$ & $T=1000$\tabularnewline
\hline 
\hline 
$\mathrm{S}_{\max,T}\left(0\right)$ & 0.061 & 0.059 & 0.057\tabularnewline
$\mathrm{S}_{\mathrm{Dmax},T}$ & 0.035 & 0.055 & 0.058\tabularnewline
$\mathrm{R}_{\mathrm{max},T}\left(0\right)$ & 0.036 & 0.035 & 0.039\tabularnewline
$\mathrm{R}_{\mathrm{Dmax},T}$ & 0.025 & 0.032 & 0.035\tabularnewline
$\widehat{D}$ statistic & 0.731 & 0.583 & 0.102\tabularnewline
\hline 
\end{tabular}
\par\end{centering}
\end{table}

\begin{table}[H]
\caption{\label{Table Power M3-M5}Empirical small-sample power for models
M3-M5}

\begin{centering}
\begin{tabular}{lccc}
\hline 
 & \multicolumn{3}{c}{Model M3}\tabularnewline
$\alpha=0.05$ & $T=250$ & $T=500$ & $T=1000$\tabularnewline
\hline 
\hline 
$\mathrm{S}_{\max,T}\left(0\right)$ & 0.694 & 0.850 & 0.889\tabularnewline
$\mathrm{S}_{\mathrm{Dmax},T}$ & 0.734 & 0.890 & 0.921\tabularnewline
$\mathrm{R}_{\mathrm{max},T}\left(0\right)$ & 0.768 & 0.940 & 0.973\tabularnewline
$\mathrm{R}_{\mathrm{Dmax},T}$ & 0.456 & 0.752 & 0.874\tabularnewline
$\widehat{D}$ statistic & 0.961 & 0.967 & 0.790\tabularnewline
\hline 
 & \multicolumn{3}{c}{Model M4}\tabularnewline
 & $T=250$ & $T=500$ & $T=1000$\tabularnewline
\hline 
\hline 
$\mathrm{S}_{\max,T}\left(0\right)$ & 0.868 & 0.964 & 0.973\tabularnewline
$\mathrm{S}_{\mathrm{Dmax},T}$ & 0.938 & 0.988 & 0.996\tabularnewline
$\mathrm{R}_{\mathrm{max},T}\left(0\right)$ & 0.927 & 0.997 & 0.999\tabularnewline
$\mathrm{R}_{\mathrm{Dmax},T}$ & 0.775 & 0.983 & 0.998\tabularnewline
$\widehat{D}$ statistic & 1.000 & 1.000 & 1.000\tabularnewline
\hline 
 & \multicolumn{3}{c}{Model M5}\tabularnewline
 & $T=250$ & $T=500$ & $T=1000$\tabularnewline
\hline 
\hline 
$\mathrm{S}_{\max,T}$ & 0.223 & 0.475 & 0.565\tabularnewline
$\mathrm{S}_{\mathrm{Dmax},T}$ & 0.325 & 0.801 & 0.918\tabularnewline
$\mathrm{R}_{\mathrm{max},T}$ & 0.028 & 0.237 & 0.369\tabularnewline
$\mathrm{R}_{\mathrm{Dmax},T}$ & 0.025 & 0.189 & 0.304\tabularnewline
$\widehat{D}$ statistic & 0.834 & 0.695 & 0.172\tabularnewline
\hline 
\end{tabular}
\par\end{centering}
\end{table}

\begin{table}[H]
\caption{\label{Table Estimation}Summary statistics for the empirical distribution
of $\widehat{m}-m_{0}$}

\begin{centering}
\begin{tabular}{ccccc}
Percent time $\widehat{m}=m_{0}$ &  & $Q_{0.25}$ & Median & $Q_{0.75}$\tabularnewline
\hline 
\hline 
\multicolumn{5}{c}{Model M6}\tabularnewline
\hline 
85.50  & $\widehat{T}_{1}$ & 299 & 333 & 352\tabularnewline
 & $\widehat{T}_{2}$ & 632 & 663 & 688\tabularnewline
\multicolumn{5}{c}{Model M7}\tabularnewline
\hline 
80.12  & $\widehat{T}_{1}$ & 317 & 336 & 359\tabularnewline
 & $\widehat{T}_{2}$ & 623 & 655  & 685 \tabularnewline
\hline 
\end{tabular}
\par\end{centering}
\end{table}

\newpage{}

\clearpage 
\pagenumbering{arabic}
\renewcommand*{\thepage}{A-\arabic{page}}
\appendix

\clearpage{}

\newpage{}

\pagebreak{}

\section*{}
\addcontentsline{toc}{part}{Supplemental Material}
\begin{center}
\Large{{Supplemental Material} to} 
\end{center}

\begin{center}
\title{\textbf{\Large{Change-Point Analysis of Time Series with Evolutionary Spectra}}} 
\maketitle
\end{center}
\medskip{} 
\medskip{} 
\medskip{} 
\thispagestyle{empty}

\begin{center}
$\qquad$ \textsc{\textcolor{MyBlue}{Alessandro Casini}} $\qquad$ \textsc{\textcolor{MyBlue}{Pierre Perron}}\\
\small{University of Rome Tor Vergata} $\quad$ \small{Boston University} 
\\
\medskip{}
\medskip{} 
\medskip{} 
\medskip{} 
\date{\small{\today}} 
\medskip{} 
\medskip{} 
\medskip{} 
\end{center}
\begin{abstract}
{\footnotesize{}This supplemental material is structured as follows.
Section \ref{Section Cumulant and Spectra} develops asymptotic results
about high-order cumulants and spectra for locally stationary series
that are needed in the proofs of the main results and are also of
independent interest. Section \ref{Section Mathematical-Appendix}
presents the Mathematical Appendix, which includes the proofs of
the results of the paper and of Section \ref{Section Cumulant and Spectra}.
Section \ref{Section Sensitivity-Analysis} includes a sensitivity
analysis for the tuning parameter choices. In Section \ref{Section Additional-Monte-Carlo}
we present additional simulation results.}{\footnotesize\par}
\end{abstract}
\setcounter{page}{0}
\setcounter{section}{0}
\setcounter{table}{0}
\renewcommand*{\theHsection}{\the\value{section}}

\newpage{}

\begin{singlespace} 
\noindent 
\small

\allowdisplaybreaks


\renewcommand{\thepage}{S-\arabic{page}}   
\renewcommand{\thesection}{S.\Alph{section}}   
\renewcommand{\thetable}{S.\arabic{table}}    
\renewcommand{\theequation}{S.\arabic{equation}}




\section{\label{Section Cumulant and Spectra}Results About High-Order Cumulants
and Spectra of Locally Stationary Series}

This section establishes asymptotic results about high-order cumulants
and spectra for locally stationary series. These are used to derive
the limiting distributions of the test statistics  introduced in
Section \ref{Section Change-Point-Tests}. They are also of independent
interest in the literature related to locally stationary and nonstationary
processes more generally. We consider the tapered finite Fourier
transform, the local and the smoothed local periodogram. Let
\begin{align*}
d_{h,T}\left(u,\,\omega\right) & \triangleq\sum_{s=0}^{n_{T}-1}h\left(\frac{s}{n_{T}}\right)X_{\left\lfloor Tu\right\rfloor -n_{T}/2+s+1,T}\exp\left(-i\omega s\right),\\
I_{h,T}\left(u,\,\omega\right) & \triangleq\frac{1}{2\pi H_{2,n_{T}}\left(0\right)}\left|d_{h,T}\left(u,\,\omega\right)\right|^{2},
\end{align*}
where $I_{h,T}\left(u,\,\omega\right)$ is the periodogram over a
segment of length $n_{T}$ with midpoint $\left\lfloor Tu\right\rfloor $.
 The smoothed local periodogram is defined as
\begin{align*}
f_{h,T}\left(u,\,\omega\right) & =\frac{2\pi}{n_{T}}\sum_{s=1}^{n_{T}-1}W_{T}\left(\omega-\frac{2\pi s}{n_{T}}\right)I_{h,T}\left(u,\,\frac{2\pi s}{n_{T}}\right),
\end{align*}
 where $W_{T}\left(\omega\right)$ and $b_{W,T}$ are defined in Section
\ref{Section Change-Point-Tests}. Note that $d_{L,h,T}\left(u,\,\omega\right),$
$I_{L,h,T}\left(u,\,\omega\right)$ and $f_{L,h,T}\left(u,\,\omega\right)$
considered in Section \ref{Section Change-Point-Tests} are asymptotically
equivalent to $d_{h,T}\left(u,\,\omega\right),$ $I_{h,T}\left(u,\,\omega\right)$
and $f_{h,T}\left(u,\,\omega\right)$, respectively. If \eqref{Eq. (2.6.6) Condition on cumulant}
holds for $l=0$, then we can define the $r$th order cumulant spectrum
at the rescale time $u\in\left(0,\,1\right)$, 
\begin{align}
f_{\mathbf{X}}^{\left(a_{1},\ldots,a_{r}\right)}\left(u,\,\omega_{1},\ldots,\,\omega_{r-1}\right) & =\left(2\pi\right)^{r-1}\sum_{k_{1},\ldots,\,k_{r-1}=-\infty}^{\infty}\kappa_{\mathbf{X},Tu}^{\left(a_{1},\ldots,a_{r}\right)}\left(k_{1},\ldots,\,k_{r-1}\right)\exp\left(-i\sum_{j=1}^{k-1}\omega_{j}k_{j}\right),\label{Eq. Higher-order Spectra}
\end{align}
for any $r$ tuple $a_{1},\ldots,a_{r}$ with $r=2,\,3,\ldots$

\subsection{Local Finite Fourier Transform}

We first present the asymptotic expression for the joint cumulants
of the finite Fourier transform. Next, we use this result to obtain
the limit distribution of the transform. This result is subsequently
used to derive the second-order properties of the local periodogram
and smoothed local periodogram in the next subsections. Corresponding
results for a stationary series can be found in \citet{brillinger:75}
and references therein. Let $\mathbf{d}_{h,T}\left(u,\,\omega\right)=[d_{h,T}^{\left(a_{j}\right)}\left(u,\,\omega\right)]\,\,(j=1,\ldots,\,r)$,
\begin{align*}
H_{n_{T}}^{\left(a_{1},\ldots,\,a_{r}\right)}\left(\omega\right) & =\sum_{s=0}^{n_{T}-1}\left(\prod_{j=1}^{r}h_{a_{j}}\left(s/n_{T}\right)\right)\exp\left(-i\omega s\right),\quad\mathrm{and}\\
H^{\left(a_{1},\ldots,\,a_{r}\right)}\left(\omega\right) & =\int\left(\prod_{j=1}^{r}h_{a_{j}}\left(t\right)\right)\exp\left(-i\omega t\right)dt.
\end{align*}
Let $\mathscr{N}_{p}^{\mathrm{C}}\left(\mathbf{c},\,\Sigma\right)$
denote the complex normal distribution for some $p$-dimensional vector
$\mathbf{c}$ and $p\times p$ Hermitian positive semidefinite matrix
$\Sigma.$ 
\begin{thm}
\label{Theorem 4.3.1-2 Brillinger}Let Assumption \ref{Assumption Locally Stationary},
\ref{Assumption WZ 2011} with $l=0$ and Assumption \ref{Taper, nT, W and K1, b1}-(ii)
hold. Let $h_{a_{j}}\left(x\right)$ satisfy Assumption \ref{Taper, nT, W and K1, b1}-(i)
for all $j=1,\ldots,\,p$. We have
\begin{align*}
\mathrm{cum} & \left(d_{h,T}^{\left(a_{1}\right)}\left(u,\,\omega_{1}\right),\ldots,\,d_{h,T}^{\left(a_{r}\right)}\left(u,\,\omega_{r}\right)\right)\\
 & =\left(2\pi\right)^{r-1}H_{n_{T}}^{\left(a_{1},\ldots,\,a_{r}\right)}\left(\sum_{j=1}^{r}\omega_{j}\right)f_{\mathbf{X}}^{\left(a_{1},\ldots,a_{r}\right)}\left(u,\,\omega_{1},\ldots,\,\omega_{r-1}\right)+\varepsilon_{T},
\end{align*}
 where $\varepsilon_{T}=o\left(n_{T}\right)$ uniformly in $\omega_{j}$
$\left(j=1,\ldots,\,r\right)$. If Assumption \ref{Assumption WZ 2011}
holds with $l=1,$ then $\varepsilon_{T}=O\left(n_{T}/T\right)$ uniformly
in $\omega_{j}$ $\left(j=1,\ldots,\,r\right)$. Furthermore, 
\begin{align*}
f_{\mathbf{X}}^{\left(a_{1},\ldots,a_{r}\right)}\left(u,\,\omega_{1},\ldots,\,\omega_{r-1}\right) & =A^{\left(a_{1}\right)}\left(\left\lfloor Tu\right\rfloor ,\,\omega_{1}\right)\cdots A^{\left(a_{p}\right)}\left(\left\lfloor Tu\right\rfloor ,\,\omega_{p}\right)g_{p}\left(\omega_{1},\ldots,\,\omega_{p-1}\right),
\end{align*}
i.e., the spectrum that corresponds to the spectral representation
\eqref{Eq. Spectral Rep of SLS} with $m_{0}=0$.
\end{thm}
\begin{thm}
\label{Theorem 4.4.1-2 Brillinger}Let Assumption \ref{Assumption Locally Stationary},
\ref{Assumption WZ 2011} with $l=0$ and Assumption \ref{Taper, nT, W and K1, b1}-(ii)
hold. Let $h_{a_{j}}\left(x\right)$ satisfy Assumption \ref{Taper, nT, W and K1, b1}-(i)
for all $j=1,\ldots,\,p$. We have: (i) If $2\omega_{j},\,\omega_{j}\pm\omega_{k}\not\equiv0\,(\mathrm{mod\,}2\pi)$
for $1\leq j<k\leq J_{\omega}$ with $1\leq J_{\omega}<\infty$, $\mathbf{d}_{h,T}\left(u,\,\omega_{j}\right)$
$(j=1,\ldots,\,J_{\omega})$ are asymptotically independent $\mathscr{N}_{p}^{\mathrm{C}}(0,2\pi n_{T}[H^{\left(a_{l},a_{r}\right)}\left(0\right)$
$f_{\mathbf{X}}^{\left(a_{l},a_{r}\right)}\left(u,\,\omega_{j}\right)])$
$\left(l,\,r=1,\ldots,\,p\right)$ variables; (ii) If $\omega=0,\,\pm\pi,\,\pm2\pi,\,\pm3\pi,$
$\ldots,$ $\mathbf{d}_{h,T}\left(u,\,\omega\right)$ is asymptotically
$\mathscr{N}_{p}(0,2\pi n_{T}[H^{\left(a_{l},a_{r}\right)}\left(0\right)f^{\left(a_{l},a_{r}\right)}\left(u,\,\omega\right)])$
$\left(l,\,r=1,\ldots,\,p\right)$ independently from the previous
variates.
\end{thm}
When the series is stationary, (i.e., $f_{\mathbf{X}}^{\left(a_{1},\ldots,a_{r}\right)}\left(u,\,\omega_{1},\ldots,\,\omega_{r-1}\right)$
does not depend on $u$), the results in Theorem \ref{Theorem 4.3.1-2 Brillinger}-\ref{Theorem 4.4.1-2 Brillinger}
reduce to the well-known results on the cumulants and asymptotic distribution
of the Fourier transform when the latter is constructed using a segment
of length $n_{T}.$ 

\subsection{Local Periodogram}

We now study several properties of the tapered local periodogram.
We begin with the finite-sample bias and variance. We then present
results about its asymptotic distribution which allow us to conclude
that the local periodogram evaluated at distinct ordinates results
in estimates that are asymptotically independent thereby mirroring
the stationary case. This result is exploited when deriving the limit
distribution of the test statistics that do not require knowledge
of the frequency at which the change-point occurs. 
\begin{thm}
\label{Theorem 5.2.3 Brillinger}Let Assumption \ref{Assumption Locally Stationary},
\ref{Assumption WZ 2011} with $l=0$ and Assumption \ref{Taper, nT, W and K1, b1}-(i,ii)
hold. We have for $-\infty<\omega<\infty$, 
\begin{align}
\mathbb{E}\left(I_{h,T}\left(u,\,\omega\right)\right) & =\left(\int_{-\pi}^{\pi}\left|H_{n_{T}}\left(\alpha\right)\right|^{2}d\alpha\right)^{-1}\int_{-\pi}^{\pi}\left|H_{n_{T}}\left(\alpha\right)\right|^{2}f_{\mathbf{X}}\left(u,\,\omega-\alpha\right)d\alpha+O\left(\frac{\log\left(n_{T}\right)}{n_{T}}\right)\label{Eq. Bias Local periodogram}\\
 & =f_{\mathbf{X}}\left(u,\,\omega\right)+\frac{1}{2}\left(\frac{n_{T}}{T}\right)^{2}\left(\int_{0}^{1}h^{2}\left(x\right)dx\right)^{-1}\int_{0}^{1}x^{2}h^{2}\left(x\right)dx\frac{\partial^{2}}{\partial u}f_{\mathbf{X}}\left(u,\,\omega\right)\nonumber \\
 & \quad+o\left(\left(\frac{n_{T}}{T}\right)^{2}\right)+O\left(\frac{\log\left(n_{T}\right)}{n_{T}}\right).\nonumber 
\end{align}
\end{thm}
The first equality  shows that the expected value of $I_{h,T}\left(u,\,\omega\right)$
is a weighted average of the local spectral density at rescaled time
$u$ with weights concentrated in a neighborhood of $\omega$ and
relative weights determined by the taper. The second equality shows
that $I_{h,T}\left(u,\,\omega\right)$ is asymptotically unbiased
for $f_{\mathbf{X}}\left(u,\,\omega\right)$ and provides a bound
on the asymptotic bias. 
\begin{thm}
\label{Theorem 5.2.8 Brillinger}Let Assumption \ref{Assumption Locally Stationary},
\ref{Assumption WZ 2011} with $l=0$ and Assumption \ref{Taper, nT, W and K1, b1}-(i,ii)
hold. We have: (i) For $-\infty<\omega_{j},\,\omega_{k}<\infty$,
\begin{align}
\mathrm{Cov} & \left\{ I_{h,T}\left(u,\,\omega_{j}\right),\,I_{h,T}\left(u,\,\omega_{k}\right)\right\} \label{Eq. (5.2.23) Brillinger}\\
 & =\left|H_{2,n_{T}}\left(0\right)\right|^{-2}\left(\left|H_{2,n_{T}}\left(\omega_{j}-\omega_{k}\right)\right|^{2}+\left|H_{2,n_{T}}\left(\omega_{j}+\omega_{k}\right)\right|^{2}\right)f_{\mathbf{X}}\left(u,\,\omega_{j}\right)^{2}+O\left(n_{T}^{-1}\right),\nonumber 
\end{align}
where $O(n_{T}^{-1})$ is uniform in $\omega_{j}$ and $\omega_{k}$;
(ii) If $2\omega_{j},\,\omega_{j}\pm\omega_{k}\not\equiv0\,(\mathrm{mod\,}2\pi)$
with $1\leq j<k\leq J_{\omega}$, the variables $I_{h,T}\left(u,\,\omega_{j}\right)$
$\left(j=1,\ldots,\,J_{\omega}\right)$ are asymptotically independent
$f_{\mathbf{X}}\left(u,\,\omega_{j}\right)\chi_{2}^{2}/2$ variates.
Also, if $\omega=\pm\pi,\,\pm3\pi,\ldots,$ $I_{h,T}\left(u,\,\omega\right)$
is asymptotically $f_{\mathbf{X}}\left(u,\,\omega\right)\chi_{1}^{2}$,
independent of the previous variates. 
\end{thm}
The variance expression for the tapered local periodogram follows
 as a special case of \eqref{Eq. (5.2.23) Brillinger}.

\subsection{Smoothed Local Periodogram}

We now extend Theorem \ref{Theorem 5.2.3 Brillinger}-\ref{Theorem 5.2.8 Brillinger}
to the smoothed local periodogram. Since our test statistics are based
on it, these results are directly employed to derive their limiting
null distributions.   
\begin{thm}
\label{Theorem 5.6.1}Let Assumption \ref{Assumption Locally Stationary},
\ref{Assumption WZ 2011} with $l=0$ and Assumption \ref{Taper, nT, W and K1, b1}
hold. Let $b_{W,T}\rightarrow0$ as $T\rightarrow\infty$ with $b_{W,T}n_{T}\rightarrow\infty$.
Then,
\begin{align}
\mathbb{E}\left(f_{h,T}\left(u,\,\omega\right)\right) & =\int_{-\infty}^{\infty}W\left(\beta\right)f_{\mathbf{X}}\left(u,\,\omega-b_{W,T}\beta\right)d\beta+O\left(\left(n_{T}b_{W,T}\right)^{-1}\right)+O\left(\log\left(n_{T}\right)n_{T}^{-1}\right)\label{Eq. (5.6.7)}\\
 & =f_{\mathbf{X}}\left(u,\,\omega\right)+\frac{1}{2}\left(\frac{n_{T}}{T}\right)^{2}\left(\int_{0}^{1}h^{2}\left(x\right)dx\right)^{-1}\int_{0}^{1}x^{2}h^{2}\left(x\right)dx\frac{\partial^{2}}{\partial u}f\mathbf{_{X}}\left(u,\,\omega\right)\nonumber \\
 & \quad+\frac{1}{2}b_{W,T}^{2}\int_{0}^{1}x^{2}W\left(x\right)dx\frac{\partial^{2}}{\partial\omega}f_{\mathbf{X}}\left(u,\,\omega\right)+O\left(\left(n_{T}/T\right)^{-2}\right)+O\left(\log\left(n_{T}\right)n_{T}^{-1}\right)+o\left(b_{W,T}^{2}\right).\nonumber 
\end{align}
 The error terms are uniform in $\omega$. 
\end{thm}

\begin{thm}
\label{Theorem 5.6.4 Brillinger}Let Assumption \ref{Assumption Locally Stationary},
\ref{Assumption WZ 2011} with $l=0$ and Assumption \ref{Taper, nT, W and K1, b1}-(i,ii)
hold. Let $b_{W,T}\rightarrow0$ as $T\rightarrow\infty$ with $b_{W,T}n_{T}\rightarrow\infty$.
Then, $f_{h,T}\left(u,\,\omega_{1}\right),\ldots,\,f_{h,T}\left(u,\,\omega_{J_{\omega}}\right)$
are asymptotically jointly normal satisfying 
\begin{align}
\lim_{T\rightarrow\infty} & n_{T}b_{W,T}\mathrm{Cov}\left(f_{h,T}\left(u,\,\omega_{j}\right),\,f_{h,T}\left(u,\,\omega_{k}\right)\right)\label{Eq. (5.6.19)}\\
 & =2\pi\left[\eta\left\{ \omega_{j}-\omega_{k}\right\} +\eta\left\{ \omega_{j}+\omega_{k}\right\} \right]\int h\left(t\right)^{4}dt\,\left[\int h\left(t\right)^{2}dt\right]^{-2}\int W\left(\alpha\right)^{2}d\alpha\,f_{\mathbf{X}}\left(u,\,\omega_{j}\right)^{2}.\nonumber 
\end{align}
\end{thm}
The variance expression for $f_{h,T}\left(u,\,\omega_{j}\right)$
follows immediately as a special case of \eqref{Eq. (5.6.19)}. Consistency
of the spectral density estimates of a stationary time series was
obtained by \citet{grenander/rosenblatt:57} and \citet{parzen:57}.
Asymptotic normality was considered by \citet{rosenblatt:59}, \citet{brillinger/rosenblatt:67},
\citet{hannan:70} and \citet{anderson:71}. Theorem \ref{Theorem 5.6.4 Brillinger}
presents corresponding results for the locally stationary case which
highlight the effect of the time-smoothing in addition to the smoothing
over the frequency-domain. \citet{panaretos/tavakoli:13} established
similar results for functional stationary processes while \citet{aue/vandelft:2020}
established some results for functional locally stationary processes
using a different notion of local stationarity. \citet{paparoditis:2009}
established similar results for linear locally stationary processes.

\section{\label{Section Mathematical-Appendix}Mathematical Appendix}

\subsection{\label{Subsec Preliminary-Lemmas}Preliminary Lemmas}

Let $L_{T}:\,\mathbb{R}\rightarrow\mathbb{R}$, $T\in\mathbb{R}_{+}$
be the $2\pi$-periodic extension of 
\begin{align*}
L_{T}\left(\omega\right) & \triangleq\begin{cases}
T, & \left|\omega\right|\leq1/T,\\
1/\left|\omega\right|, & 1/T\leq\left|\omega\right|\leq\pi.
\end{cases}
\end{align*}
 For a complex-valued function $w$ define $H_{n_{T}}\left(w\left(\cdot\right),\,\omega\right)=\sum_{s=0}^{n_{T}-1}w\left(s\right)\exp\left(-i\omega s\right)$,
and, for the taper $h\left(x\right),$ $H_{k,n_{T}}\left(\omega\right)=H_{n_{T}}\left(h^{k}\left(\frac{\cdot}{n_{T}}\right),\,\omega\right),$
and $H_{n_{T}}\left(\omega\right)=H_{1,n_{T}}\left(\omega\right).$ 
\begin{lem}
\label{Lemma A.4 Dahl97}Let  $\varPi\triangleq(-\pi,\,\pi]$. With
a constant $K$ independent of $T$ the following properties hold:
(i) $L_{T}\left(\omega\right)$ is monotone increasing in $T$ and
decreasing in $\omega\in\left[0,\,\pi\right]$; (ii) $\int_{\varPi}L_{T}\left(\alpha\right)d\alpha\leq K\ln T$
for $T>1$; (iii) $\int_{\varPi}L_{T}\left(\alpha\right)^{k}d\alpha\leq KT^{k-1}$
for $k\geq2$.
\end{lem}
\noindent\textit{Proof of Lemma }\ref{Lemma A.4 Dahl97}\textit{.
}See Lemma A.4 in \citeReferencesSupp{dahlhaus:96}. $\square$
\begin{lem}
\label{Lemma A.5 Dahl97}Suppose $h\left(\cdot\right)$ satisfies
Assumption \ref{Taper, nT, W and K1, b1} and $\vartheta:\,\left[0,\,1\right]\rightarrow\mathbb{R}$
is differentiable with bounded derivative. Then we have for $0\leq t\leq n_{T}$,
\begin{align*}
H_{n_{T}}\left(\vartheta\left(\frac{\cdot}{T}\right)h\left(\frac{\cdot}{n_{T}}\right),\,\omega\right) & =\vartheta\left(\frac{t}{T}\right)H_{n_{T}}\left(\omega\right)+O\left(\sup_{x}\left|d\vartheta\left(x\right)/dx\right|\frac{n_{T}}{T}L_{n_{T}}\left(\omega\right)\right)\\
 & =O\left(\sup_{x\leq n_{T}/T}\left|\vartheta\left(x\right)\right|L_{n_{T}}\left(\omega\right)+\sup_{x}\left|d\vartheta\left(x\right)/dx\right|\frac{n_{T}}{T}L_{n_{T}}\left(\omega\right)\right).
\end{align*}
The same holds, if $\vartheta\left(\cdot/T\right)$ is replaced on
the left side by numbers $\vartheta_{s,T}$ with $\sup_{s}\left|\vartheta_{s,T}-\vartheta\left(s/T\right)\right|=O\left(T^{-1}\right)$. 
\end{lem}
\noindent\textit{Proof of Lemma }\ref{Lemma A.5 Dahl97}. \citeReferencesSupp{dahlhaus:96}
proved this result under differentiability of $h\left(\cdot\right)$.
By Abel's transformation {[}cf. Exercise 1.7.13 in \citeReferencesSupp{brillinger:75}{]},
\begin{align}
H_{n_{T}}\left(\vartheta\left(\frac{\cdot}{T}\right)h\left(\frac{\cdot}{n_{T}}\right),\,\omega\right)-\vartheta\left(\frac{t}{T}\right)H_{n_{T}}\left(\omega\right) & =\sum_{s=0}^{n_{T}-1}\left[\vartheta\left(\frac{s}{T}\right)-\vartheta\left(\frac{t}{T}\right)\right]h\left(\frac{s}{n_{T}}\right)\exp\left(-i\omega s\right)\nonumber \\
 & =-\sum_{s=0}^{n_{T}-1}\left[\vartheta\left(\frac{s}{T}\right)-\vartheta\left(\frac{s-1}{T}\right)\right]H_{s}\left(h\left(\frac{\cdot}{n_{T}}\right),\,\omega\right)\nonumber \\
 & \quad+\left[\vartheta\left(\frac{n_{T}-1}{T}\right)-\vartheta\left(\frac{t}{T}\right)\right]H_{n_{T}}\left(h\left(\frac{\cdot}{n_{T}}\right),\,\omega\right).\label{Eq. in A.5 In Dal97}
\end{align}
By repeated applications of Abel's transformation,
\begin{align*}
H_{s}\left(h\left(\frac{\cdot}{n_{T}}\right),\,\omega\right) & =\sum_{t=0}^{s-1}h\left(\frac{t}{n_{T}}\right)\exp\left(-i\omega t\right)\\
 & =\sum_{t=0}^{s-1}\left(h\left(\frac{t}{n_{T}}\right)-h\left(\frac{t-1}{n_{T}}\right)\right)H_{t}\left(1,\,\omega\right)\\
 & \quad+h\left(\frac{n_{T}-1}{n_{T}}\right)H_{n_{T}}\left(1,\,\omega\right)\\
 & =\sum_{t=0}^{s-1}\left(h\left(\frac{t}{n_{T}}\right)-h\left(\frac{t-1}{n_{T}}\right)\right)H_{t}\left(1,\,\omega\right)+0,
\end{align*}
where we have used $h\left(\left(n_{T}-1\right)/n_{T}\right)-h\left(1\right)=O(n_{T}^{-1})$
and $h\left(x\right)=0$ for $x\notin[0,\,1)$. Since $h\left(\cdot\right)$
is of bounded variation, if $\left|\omega\right|\leq1/n_{T}$ we have
\begin{align*}
\sum_{t=0}^{s-1}\left|\left(h\left(\frac{t}{n_{T}}\right)-h\left(\frac{t-1}{n_{T}}\right)\right)\right|\left|H_{t}\left(1,\,\omega\right)\right| & \leq\sum_{t=0}^{s-1}t\left|\left(h\left(\frac{t}{n_{T}}\right)-h\left(\frac{t-1}{n_{T}}\right)\right)\right|\\
 & \leq\left(s-1\right)\sum_{t=0}^{s-1}\left|\left(h\left(\frac{t}{n_{T}}\right)-h\left(\frac{t-1}{n_{T}}\right)\right)\right|\\
 & \leq C\left(s-1\right),
\end{align*}
whereas if $1/n_{T}\leq\left|\omega\right|\leq\pi$ we have, 
\begin{align*}
\sum_{t=0}^{s-1}\left|\left(h\left(\frac{t}{n_{T}}\right)-h\left(\frac{t-1}{n_{T}}\right)\right)\right|\left|H_{t}\left(1,\,\omega\right)\right| & \leq C\frac{1}{\left|\omega\right|}\sum_{t=0}^{s-1}\left|\left(h\left(\frac{t}{n_{T}}\right)-h\left(\frac{t-1}{n_{T}}\right)\right)\right|\\
 & \leq C\frac{1}{\left|\omega\right|}.
\end{align*}
 Thus, $H_{s}(h(\cdot/n_{T}),\,\omega)\leq L_{s}\left(\omega\right)\leq L_{n_{T}}\left(\omega\right)$
where the last inequality follows by Lemma \ref{Lemma A.4 Dahl97}-(i).
It follows from \eqref{Eq. in A.5 In Dal97} that, 
\begin{align*}
H_{n_{T}} & \left(\vartheta\left(\frac{\cdot}{T}\right)h\left(\frac{\cdot}{n_{T}}\right),\,\omega\right)-\vartheta\left(\frac{t}{T}\right)H_{n_{T}}\left(\omega\right)\\
 & =O\left(\sup_{x\leq n_{T}/T}\left|\vartheta\left(x\right)\right|L_{n_{T}}\left(\omega\right)+\sup_{x}\left|d\vartheta\left(x\right)/dx\right|\frac{n_{T}}{T}L_{n_{T}}\left(\omega\right)\right).\,\square
\end{align*}

\begin{lem}
\label{Lemma P4.1 Brillinger}Assume that $h^{\left(a_{j}\right)}\left(x\right)$
satisfies Assumption \ref{Taper, nT, W and K1, b1}-(i) for all $j=1,\ldots,\,p$,
then we have for some $C$ with $0<C<\infty$, 
\begin{align*}
\biggl| & \sum_{s=0}^{n_{T}-1}h_{T}^{\left(a_{1}\right)}\left(s+k_{1}\right)\cdots h_{T}^{\left(a_{p-1}\right)}\left(s+k_{p-1}\right)h_{T}^{\left(a_{1}\right)}\left(s\right)\exp\left(-i\omega s\right)-H_{T}^{\left(a_{1},\cdots,\,a_{p}\right)}\left(\omega\right)\biggl|\\
 & \leq C\left(\left|k_{1}\right|+\ldots+\left|k_{p-1}\right|\right).
\end{align*}
\end{lem}
\noindent\textit{Proof of Lemma }\ref{Lemma P4.1 Brillinger}\textit{.
}See Lemma P4.1 in \citeReferencesSupp{brillinger:75}. $\square$
\begin{lem}
\label{Lemma P4.5 Brillinger}Let $\left\{ Y_{T}\right\} $ be a sequence
of $p$ vector-valued random variables, with (possibly) complex components,
and such that all cumulants of the variate $(Y_{T}^{\left(a_{1}\right)},\,\overline{Y}_{T}^{\left(a_{1}\right)},\ldots,\,Y_{T}^{\left(a_{p}\right)},\,\overline{Y}_{T}^{\left(a_{p}\right)})$
exist and tend to the corresponding cumulants of a variate $(Y^{\left(a_{1}\right)},\,\overline{Y}^{\left(a_{1}\right)},\ldots,\,Y^{\left(a_{p}\right)},\,\overline{Y}^{\left(a_{p}\right)})$
that is determined by its moments. Then $Y_{T}$ tends in distribution
to a variate having components $Y^{\left(a_{1}\right)},\ldots,\,Y^{\left(a_{p}\right)}$. 
\end{lem}
\noindent\textit{Proof of Lemma }\ref{Lemma P4.5 Brillinger}\textit{.
}It follows from Lemma P4.5 in \citeReferencesSupp{brillinger:75}.
$\square$

\subsection{\label{subsec:Proofs-of-Cumulant}Proofs of the Results of Section
\ref{Section Cumulant and Spectra}}

\subsubsection{Proof of Theorem \ref{Theorem 4.3.1-2 Brillinger}}

For $\left\lfloor Tu\right\rfloor -n_{T}/2+1\leq t_{1},\ldots,\,t_{p}\leq\left\lfloor Tu\right\rfloor +n_{T}/2-1$,
\begin{align*}
\mathrm{cum} & \left(X_{t_{1},T},\ldots,\,X_{t_{p},T}\right)\\
 & =\int_{-\pi}^{\pi}\cdots\int_{-\pi}^{\pi}\exp\left(it_{1}\omega_{1}+\cdots+it_{p}\omega_{p}\right)\\
 & \quad\times A_{t_{1},T}^{0}\left(\omega_{1}\right)\cdots A_{t_{p},T}^{0}\left(\omega_{p}\right)\eta\left(\sum_{j=1}^{p}\omega_{j}\right)g_{p}\left(\omega_{1},\ldots,\,\omega_{p-1}\right)d\omega_{1}\cdots d\omega_{p}.
\end{align*}
 We can replace $A_{t_{j},T}^{0}\left(\omega_{j}\right)$ by $A\left(t_{j}/T,\,\omega_{j}\right)$
using \eqref{Eq. 2.4 Smothenss Assumption on A}, and then replace
$A\left(t_{j}/T,\,\omega_{j}\right)$ by $A\left(\left\lfloor Tu\right\rfloor ,\,\omega_{j}\right)$
using the smoothness of $A\left(u,\,\cdot\right)$. Altogether, this
gives an error $O\left(n_{T}/T\right)$. Let $t_{1}=t_{p}+k_{1},\ldots,\,t_{p-1}=t_{p}+k_{p-1}.$
We have
\begin{align}
\mathrm{cum} & \left(X_{t_{1},T},\ldots,\,X_{t_{p},T}\right)\nonumber \\
 & =\int_{-\pi}^{\pi}\cdots\int_{-\pi}^{\pi}\exp\left(i\left(\left(\omega_{1}+\cdots+\omega_{p-1}\right)t_{p}+\omega_{1}k_{1}+\cdots+\omega_{p-1}k_{p-1}+t_{p}\omega_{p}\right)\right)\nonumber \\
 & \quad\times A\left(\left\lfloor Tu\right\rfloor ,\,\omega_{1}\right)\cdots A\left(\left\lfloor Tu\right\rfloor ,\,\omega_{p}\right)\eta\left(\sum_{j=1}^{p}\omega_{j}\right)g_{p}\left(\omega_{1},\ldots,\,\omega_{p-1}\right)d\omega_{1}\cdots d\omega_{p}+O\left(n_{T}/T\right)\nonumber \\
 & =\int_{-\pi}^{\pi}\cdots\int_{-\pi}^{\pi}\exp\left(i\left(\left(\omega_{1}+\cdots+\omega_{p-1}+\omega_{p}\right)t_{p}+\omega_{1}k_{1}+\cdots+\omega_{p-1}k_{p-1}\right)\right)\nonumber \\
 & \quad\times A\left(\left\lfloor Tu\right\rfloor ,\,\omega_{1}\right)\cdots A\left(\left\lfloor Tu\right\rfloor ,\,\omega_{p}\right)\eta\left(\sum_{j=1}^{p}\omega_{j}\right)g_{p}\left(\omega_{1},\ldots,\,\omega_{p-1}\right)d\omega_{1}\cdots d\omega_{p}+O\left(n_{T}/T\right)\nonumber \\
 & \triangleq\kappa_{Tu,t_{p}}\left(k_{1}\ldots,\,k_{p-1}\right)+O\left(n_{T}/T\right).\label{Eq. cum kappa}
\end{align}
 This shows that $\mathrm{cum}\left(X_{t_{1},T},\ldots,\,X_{t_{p},T}\right)$
depends on $t_{p}$ only through $\exp\left(i\left(\omega_{1}+\cdots+\omega_{p-1}+\omega_{p}\right)t_{p}\right)$.
The cumulant of interest in Theorem \ref{Theorem 4.3.1-2 Brillinger}
has the following form, 
\begin{align*}
\mathrm{cum} & \left(d_{h,T}^{\left(a_{1}\right)}\left(u,\,\omega_{1}\right),\ldots,\,d_{h,T}^{\left(a_{p}\right)}\left(u,\,\omega_{p}\right)\right)\\
 & =\int_{-\pi}^{\pi}\cdots\int_{-\pi}^{\pi}H_{n_{T}}\left(A_{\left\lfloor Tu\right\rfloor -n_{T}/2+1+\cdot,T}^{0,(a_{1})}\left(\gamma_{1}\right)h_{a_{1}}\left(\frac{\cdot}{n_{T}}\right),\,\omega_{1}-\gamma_{1}\right)\\
 & \times H_{n_{T}}\left(A_{\left\lfloor Tu\right\rfloor -n_{T}/2+1+\cdot,T}^{0,(a_{2})}\left(\gamma_{2}\right)h_{a_{2}}\left(\frac{\cdot}{n_{T}}\right),\,\omega_{2}-\gamma_{2}\right)\\
 & \times\cdots\\
 & \times H_{n_{T}}\left(A_{\left\lfloor Tu\right\rfloor -n_{T}/2+1+\cdot,T}^{0,(a_{p})}\left(\gamma_{p}\right)h_{a_{p}}\left(\frac{\cdot}{n_{T}}\right),\,\omega_{p}-\gamma_{p}\right)\\
 & \times\exp\left\{ i\left(\left(\gamma_{1}+\cdots+\gamma_{p}\right)\left\lfloor Tu\right\rfloor \right)\right\} \eta\left(\sum_{j=1}^{p}\gamma_{j}\right)g_{p}\left(\gamma_{1},\ldots,\,\gamma_{p-1}\right)d\gamma_{1}\cdots d\gamma_{p}+o\left(1\right).
\end{align*}
 By Lemma \ref{Lemma A.5 Dahl97}, the latter is equal to 
\begin{align}
\int_{-\pi}^{\pi}\cdots\int_{-\pi}^{\pi} & A^{\left(a_{1}\right)}\left(u,\,\gamma_{1}\right)\cdots A^{\left(a_{p}\right)}\left(u,\,\gamma_{p}\right)\label{Eq. (A.3) Dal97 before}\\
 & \times H_{n_{T}}^{\left(a_{1}\right)}\left(\omega_{1}-\gamma_{1}\right)\cdots H_{n_{T}}^{\left(a_{p}\right)}\left(\omega_{p}-\gamma_{p}\right)\nonumber \\
 & \times\exp\left(i\left(\left(\gamma_{1}+\cdots+\gamma_{p}\right)\left\lfloor Tu\right\rfloor \right)\right)\eta\left(\sum_{j=1}^{p}\gamma_{j}\right)g_{p}\left(\gamma_{1},\ldots,\,\gamma_{p-1}\right)d\gamma_{1}\cdots d\gamma_{p},\nonumber 
\end{align}
plus a remainder term $R_{u}$ with 
\begin{align}
\left|R_{u}\right| & \leq C\frac{n_{T}}{T}\int_{-\pi}^{\pi}\cdots\int_{-\pi}^{\pi}L_{n_{T}}\left(\omega_{1}-\gamma_{1}\right)\cdots L_{n_{T}}\left(\omega_{p}-\gamma_{p}\right)\exp\left(i\left(\left(\gamma_{1}+\cdots+\gamma_{p}\right)\left\lfloor Tu\right\rfloor \right)\right)\label{Eq. (A.4) Dahl97}\\
 & \quad\times\eta\left(\sum_{j=1}^{p}\gamma_{j}\right)g_{p}\left(\gamma_{1},\ldots,\,\gamma_{p-1}\right)d\gamma_{1}\cdots d\gamma_{p}\nonumber \\
 & \leq C\frac{n_{T}}{T}\int_{-\pi}^{\pi}\cdots\int_{-\pi}^{\pi}L_{n_{T}}\left(\omega_{1}-\gamma_{1}\right)\cdots L_{n_{T}}\left(\omega_{p}-\gamma_{p}\right)d\gamma_{1}\cdots d\gamma_{p}\nonumber \\
 & \leq C\frac{n_{T}}{T}\left(\ln n_{T}\right)^{p},\nonumber 
\end{align}
 where we have used $g_{p}\left(\gamma_{1},\ldots,\,\gamma_{p-1}\right)\leq\mathrm{const}_{p},$
the fact that $\int_{-\pi}^{\pi}\exp\left\{ i\left(\gamma\left\lfloor Tu\right\rfloor \right)\right\} d\gamma=2\sin\left(\pi\left\lfloor Tu\right\rfloor \right)/\left\lfloor Tu\right\rfloor $,
and the third inequality follows from Lemma \ref{Lemma A.4 Dahl97}-(ii).

Next, note that the function $H_{n_{T}}\left(\omega\right)$ will
have substantial magnitude only for $\omega$ near some multiple of
$2\pi$. Thus, by continuity of $A\left(\cdot,\,\omega\right)$, $g_{p}$,
and of the exponential function we have that \eqref{Eq. (A.3) Dal97 before}
is equal to
\begin{align}
\int_{-\pi}^{\pi}\cdots\int_{-\pi}^{\pi} & A^{\left(a_{1}\right)}\left(u,\,\omega_{1}\right)\cdots A^{\left(a_{p}\right)}\left(u,\,\omega_{p}\right)\label{Eq. (A.3) Dalh97}\\
 & \times H_{n_{T}}^{\left(a_{1}\right)}\left(\omega_{1}-\gamma_{1}\right)\cdots H_{n_{T}}^{\left(a_{p}\right)}\left(\omega_{p}-\gamma_{p}\right)\nonumber \\
 & \times\exp\left(i\left(\left(\omega_{1}+\cdots+\omega_{p}\right)\left\lfloor Tu\right\rfloor \right)\right)\eta\left(\sum_{j=1}^{p}\omega_{j}\right)g_{p}\left(\omega_{1},\ldots,\,\omega_{p-1}\right)d\gamma_{1}\cdots d\gamma_{p}.\nonumber 
\end{align}
 By Lemma \ref{Lemma P4.1 Brillinger}, 
\begin{align*}
| & \sum_{s=0}^{n_{T}-1}h_{a_{1}}\left(\frac{s+k_{1}}{n_{T}}\right)\cdots h_{a_{p-1}}\left(\frac{s+k_{p-1}}{n_{T}}\right)h_{a_{p}}\left(\frac{s}{n_{T}}\right)\exp\left(i\sum_{j=1}^{p}\omega_{j}s\right)-H_{T}^{\left(a_{1},\cdots,\,a_{p}\right)}\left(\sum_{j=1}^{p}\omega_{j}\right)|\\
 & \leq C\left(\left|k_{1}\right|+\ldots+\left|k_{p-1}\right|\right).
\end{align*}
Thus, \eqref{Eq. (A.3) Dalh97} is equal to
\begin{align}
\sum_{k_{1}=-n_{T}}^{n_{T}} & \cdots\sum_{k_{p-1}=-n_{T}}^{n_{T}}\exp\left(-i\sum_{j=1}^{p-1}\omega_{j}k_{j}\right)\label{Eq. exp k H}\\
 & \times\left(\kappa_{Tu,t_{p}}^{\left(a_{1},\ldots,\,a_{p}\right)}\left(k_{1}\ldots,\,k_{p-1}\right)H_{T}^{\left(a_{1},\ldots,\,a_{p}\right)}\left(\sum_{j=1}^{p}\omega_{j}\right)+O\left(n_{T}/T\right)\right)+\varepsilon_{T},\nonumber 
\end{align}
 where $\kappa_{Tu,t_{p}}^{\left(a_{1},\ldots,\,a_{p}\right)}\left(k_{1}\ldots,\,k_{p-1}\right)=\mathrm{cum}(X_{t_{1},T}^{\left(a_{1}\right)},\ldots,\,X_{t_{p},T}^{\left(a_{p}\right)})+\left(n_{T}/T\right)$
and
\begin{align*}
\left|\varepsilon_{T}\right| & \leq C\sum_{k_{1}=-n_{T}}^{n_{T}}\cdots\sum_{k_{p-1}=-n_{T}}^{n_{T}}\kappa_{Tu,t_{p}}^{\left(a_{1},\ldots,\,a_{p}\right)}\left(k_{1}\ldots,\,k_{p-1}\right)\left(\left|k_{1}\right|+\cdots+\left|k_{p}\right|\right)<\infty.
\end{align*}
 Note that $\left|\varepsilon_{T}\right|/n_{T}\rightarrow0$ since
$\left(\left|k_{1}\right|+\cdots+\left|k_{p}\right|\right)/n_{T}\rightarrow0.$
Thus, $\varepsilon_{T}=o\left(n_{T}\right)$ uniformly in $\omega_{j}$
$\left(j=1,\ldots,\,p\right)$. Altogether we have 
\begin{align*}
\mathrm{cum} & \left(d_{h,T}^{\left(a_{1}\right)}\left(u,\,\omega_{1}\right),\ldots,\,d_{h,T}^{\left(a_{p}\right)}\left(u,\,\omega_{p}\right)\right)\\
 & =\left(2\pi\right)^{r-1}H_{n_{T}}^{\left(a_{1},\ldots,\,a_{p}\right)}\left(\sum_{j=1}^{p}\omega_{j}\right)f_{\mathbf{X}}^{\left(a_{1},\ldots,a_{p}\right)}\left(u,\,\omega_{1},\ldots,\,\omega_{p-1}\right)+\varepsilon_{T},
\end{align*}
 where $f_{\mathbf{X}}^{\left(a_{1},\ldots,a_{p}\right)}\left(u,\,\omega_{1},\ldots,\,\omega_{p-1}\right)$
is given in \eqref{Eq. Spectral Rep of SLS}. The proof for the $r$th
cumulant of $d_{h,T}^{\left(a_{j}\right)}\left(u,\,\omega_{1}\right)$
$\left(j=1,\ldots,\,r\right)$ with $r<p$ is the same as for the
$p$th cumulant. 

Note that from \eqref{Eq. exp k H} we have 
\begin{align*}
\sum_{k_{1}=-n_{T}}^{n_{T}} & \cdots\sum_{k_{p-1}=-n_{T}}^{n_{T}}\exp\left(-i\sum_{j=1}^{p-1}\omega_{j}k_{j}\right)\kappa_{Tu,t_{p}}\left(k_{1}\ldots,\,k_{p-1}\right)\\
 & =\sum_{k_{1}=-n_{T}}^{n_{T}}\cdots\sum_{k_{p-1}=-n_{T}}^{n_{T}}\int_{-\pi}^{\pi}\cdots\int_{-\pi}^{\pi}\exp(i(\left(\gamma_{1}+\cdots+\gamma_{p-1}+\gamma_{p}\right)t_{p}\\
 & \quad+\left(\omega_{1}-\gamma_{1}\right)k_{1}+\cdots+\left(\omega_{p-1}-\gamma_{p-1}\right)k_{p-1}))\\
 & \quad\times A^{\left(a_{1}\right)}\left(\left\lfloor Tu\right\rfloor ,\,\gamma_{1}\right)\cdots A^{\left(a_{p}\right)}\left(\left\lfloor Tu\right\rfloor ,\,\gamma_{p}\right)\eta\left(\sum_{j=1}^{p}\gamma_{j}\right)g_{p}\left(\gamma_{1},\ldots,\,\gamma_{p-1}\right)d\gamma_{1}\cdots d\gamma_{p}.
\end{align*}
 Since $\sum_{j=1}^{p}\gamma_{j}\equiv0\,(\mathrm{mod\,}2\pi)$, $\gamma_{p}$
is normalized and so the latter is equivalent to 
\begin{align*}
\sum_{k_{1}=-n_{T}}^{n_{T}} & \cdots\sum_{k_{p-1}=-n_{T}}^{n_{T}}\int_{-\pi}^{\pi}\cdots\int_{-\pi}^{\pi}\exp\left(i\left(\gamma_{1}k_{1}+\cdots+\gamma_{p-1}k_{p-1}\right)\right)\\
 & \times A^{\left(a_{1}\right)}\left(\left\lfloor Tu\right\rfloor ,\,\omega_{1}\right)\cdots A^{\left(a_{p}\right)}\left(\left\lfloor Tu\right\rfloor ,\,\omega_{p}\right)g_{p}\left(\omega_{1},\ldots,\,\omega_{p-1}\right)d\gamma_{1}\cdots d\gamma_{p-1},
\end{align*}
where we have used the continuity of $A\left(\cdot,\,\omega\right)$
and $g_{p}$. Then, 
\begin{equation}
A^{\left(a_{1}\right)}\left(\left\lfloor Tu\right\rfloor ,\,\omega_{1}\right)\cdots A^{\left(a_{p}\right)}\left(\left\lfloor Tu\right\rfloor ,\,\omega_{p}\right)g_{p}\left(\omega_{1},\ldots,\,\omega_{p-1}\right)=f_{\mathbf{X}}^{\left(a_{1},\ldots,\,a_{p}\right)}\left(u,\,\omega_{1},\ldots,\,\omega_{p-1}\right)\label{Eq. (Spectrum of Spectral Rep)}
\end{equation}
 is the spectrum that corresponds to the spectral representation \eqref{Eq. Spectral Rep of SLS}
with $m_{0}=0$.  In view of the following identities {[}see e.g.,
Exercise 1.7.5-(c,d) in \citeReferencesSupp{brillinger:75}{]},
\begin{align*}
\sum_{k=-n_{T}}^{n_{T}}\exp\left(-i\omega k\right) & =\frac{\sin\left(n_{T}+1/2\right)\omega}{\sin\omega/2},\qquad\int_{-\pi}^{\pi}\frac{\sin\left(n_{T}+1/2\right)\omega}{\sin\omega/2}d\omega=2\pi,
\end{align*}
 we have 
\begin{align*}
\mathrm{cum} & \left(d_{T}^{\left(a_{1}\right)}\left(u,\,\omega_{1}\right),\ldots,\,d_{T}^{\left(a_{p}\right)}\left(u,\,\omega_{p}\right)\right)\\
 & =\left(2\pi\right)^{p-1}H_{T}^{\left(a_{1},\ldots,\,a_{p}\right)}\left(\sum_{j=1}^{p}\omega_{j}\right)f^{\left(a_{1},\ldots,\,a_{p}\right)}\left(u,\,\omega_{1},\ldots,\,\omega_{p-1}\right)+\varepsilon_{T},
\end{align*}
which verifies \eqref{Eq. (Spectrum of Spectral Rep)}. $\square$ 

\subsubsection{Proof of Theorem \ref{Theorem 4.4.1-2 Brillinger}}

We have,
\begin{align*}
\mathbb{E}\left(\mathbf{d}_{h,T}\left(u,\,\omega\right)\right) & =\sum_{s=0}^{n_{T}-1}\exp\left(-i\omega s\right)\mathbb{E}\left(\mathbf{X}_{\left\lfloor Tu\right\rfloor -n_{T}/2+s+1,T}\right)\\
 & =0.
\end{align*}
By Theorem \ref{Theorem 4.3.1-2 Brillinger} we deduce 
\begin{align}
n_{T}^{-1} & \mathrm{Cov}\left(d_{h,T}^{\left(a_{l}\right)}\left(u,\,\pm\omega_{j}\right),\,d_{h,T}^{\left(a_{r}\right)}\left(u,\,\pm\omega_{k}\right)\right)\label{Eq. (Cov (dl, dr))}\\
 & =n_{T}^{-1}2\pi H_{n_{T}}^{\left(a_{l},\,a_{r}\right)}\left(\pm\omega_{j}\mp\omega_{k}\right)f_{\mathbf{X}}^{\left(a_{l},\,a_{r}\right)}\left(u,\,\pm\omega_{j}\left(n_{T}\right)\right)+o\left(1\right)+O\left(n_{T}^{-1}\right).\nonumber 
\end{align}
Note that {[}see, e.g., Lemma P4.6 in \citeReferencesSupp{brillinger:75}{]},
\begin{align}
\left|H_{n_{T}}^{\left(a_{1},\ldots,\,a_{p}\right)}\left(\omega\right)\right| & \leq\frac{C}{\left|\sin\left(\omega/2\right)\right|},\label{Eq. (Lemma P4.6, Brillinger)}
\end{align}
where $C$ is a constant with $0<C<\infty$. If $\omega_{j}\pm\omega_{k}\not\equiv0\,(\mathrm{mod\,}2\pi)$
the first term on the right-hand side of \eqref{Eq. (Cov (dl, dr))}
tends to zero using \eqref{Eq. (Lemma P4.6, Brillinger)}.  If $\pm\omega_{j}\mp\omega_{k}\equiv0\,(\mathrm{mod\,}2\pi)$
the right-hand side of \eqref{Eq. (Cov (dl, dr))} tends to 
\begin{align*}
2\pi H_{T}^{\left(a_{l},\,a_{r}\right)}\left(0\right)f_{\mathbf{X}}^{\left(a_{l},a_{r}\right)}\left(u,\,\pm\omega_{j}\right) & =2\pi\left(\int h^{\left(a_{l}\right)}\left(t\right)h^{\left(a_{r}\right)}\left(t\right)dt\right)f_{\mathbf{X}}^{\left(a_{l},a_{r}\right)}\left(u,\,\pm\omega_{j}\right).
\end{align*}
This shows that the second-order cumulants behave as indicated by
the theorem. By Theorem \ref{Theorem 4.3.1-2 Brillinger} for $r>2$,
\begin{align*}
n_{T}^{-r/2} & \mathrm{cum}\left(d_{h,T}^{\left(a_{1}\right)}\left(u,\,\pm\omega_{j_{1}}\right),\ldots,\,d_{h,T}^{\left(a_{r}\right)}\left(u,\,\pm\omega_{j_{r}}\right)\right)\\
 & =n_{T}^{-r/2}\left(2\pi\right)^{r-1}H_{n_{T}}^{\left(a_{1},\ldots,\,a_{r}\right)}\left(\pm\omega_{j_{1}}\pm\cdots\pm\omega_{j_{r}}\right)f_{\mathbf{X}}^{\left(a_{1},\ldots,\,a_{r}\right)}\left(u,\,\pm\omega_{j_{1}},\ldots,\,\pm\omega_{j_{r-1}}\right)+o\left(n_{T}^{1-r/2}\right).
\end{align*}
The latter tends to 0 as $n_{T}\rightarrow\infty$ if $r>2$ because
$H_{n_{T}}^{\left(a_{1},\ldots,\,a_{r}\right)}\left(\omega\right)=O\left(n_{T}\right)$.
Thus, also the cumulants of order higher than two behave as indicated
by the theorem. This implies that the cumulants of the considered
variables and the conjugates of those variables tend to the cumulants
of Gaussian random variable. Since the distribution of the latter
is fully determined by its moments, the theorem follows from Lemma
\ref{Lemma P4.5 Brillinger}. The second part of the theorem follows
from the fact that $\sin\left(\omega\right)=0$ for $\omega=0,\,\pm\pi,\,\pm2\pi,\,\pm2\pi,\ldots$
$\square$

\subsubsection{Proof of Theorem \ref{Theorem 5.2.3 Brillinger}}

The proof of the second equality in \eqref{Eq. Bias Local periodogram}
is similar to \citeReferencesSupp{dahlhaus:96} who proved the result
under stronger assumptions on the data taper. Using the spectral representation
\eqref{Eq. Spectral Rep of SLS},
\begin{align*}
\mathrm{cum} & \left(d_{h,T}\left(u,\,\omega\right),\,d_{h,T}\left(u,\,-\omega\right)\right)\\
 & =\sum_{t=0}^{n_{T}-1}\sum_{s=0}^{n_{T}-1}h\left(\frac{t}{T}\right)h\left(\frac{s}{T}\right)\int_{-\pi}^{\pi}\exp\left(-i\left(\omega-\eta\right)\left(s-t\right)\right)A_{\left\lfloor Tu\right\rfloor -n_{T}/2+t}^{0}\left(\eta\right)A_{\left\lfloor Tu\right\rfloor -n_{T}/2+s}^{0}\left(-\eta\right)d\eta.
\end{align*}
We use Abel's transformation to replace $A_{\left\lfloor Tu\right\rfloor -n_{T}/2+t}^{0}\left(\eta\right)$
by $A\left(u,\,\omega\right),$
\begin{align*}
| & \sum_{t=0}^{n_{T}-1}h\left(\frac{t}{n_{T}}\right)\left(A_{\left\lfloor Tu\right\rfloor -n_{T}/2+t}^{0}\left(\eta\right)-A\left(u,\,\omega\right)\right)\exp\left(-i\left(\omega-\eta\right)t\right)|\\
 & =|\sum_{t=0}^{n_{T}-1}\left(A_{\left\lfloor Tu\right\rfloor -n_{T}/2+t}^{0}\left(\eta\right)-A_{\left\lfloor Tu\right\rfloor -n_{T}/2+t-1}^{0}\left(\eta\right)\right)H_{t}\left(h\left(\frac{\cdot}{n_{T}},\,\omega-\eta\right)\right)|\\
 & \quad+|\left(A_{\left\lfloor Tu\right\rfloor -n_{T}/2+n_{T}-1}^{0}\left(\eta\right)-A\left(u,\,\omega\right)\right)H_{n_{T}}\left(h\left(\frac{\cdot}{n_{T}}\right),\,\omega-\eta\right)|\\
 & \leq O\left(\frac{n_{T}}{T}\right)L_{n_{T}}\left(\omega-\eta\right)+\left(O\left(\frac{n_{T}}{T}\right)+O\left(\left|\omega-\eta\right|\right)\right)L_{n_{T}}\left(\omega-\eta\right),
\end{align*}
 where the inequality follows from using Lemma \ref{Lemma A.5 Dahl97},
\begin{align}
\left|H_{t}\left(h\left(\frac{\cdot}{n_{T}},\,\omega-\eta\right)\right)\right| & \leq L_{t}\left(\omega-\eta\right)\leq L_{n_{T}}\left(\omega-\eta\right).\label{Eq. (H_t)<L_t<L_n}
\end{align}
Since we are dividing by $\sum_{s=0}^{n_{T}-1}h\left(s/n_{T}\right)^{2}\sim n_{T}$
we get, 
\begin{align*}
n_{T}^{-1} & \left|\sum_{t=0}^{n_{T}-1}h\left(\frac{t}{n_{T}}\right)\left(A_{\left\lfloor Tu\right\rfloor -n_{T}/2+t}^{0}\left(\eta\right)-A\left(u+\frac{t-n_{T}/2}{T},\,\omega\right)\right)\exp\left(-i\left(\omega-\eta\right)t\right)\right|\\
 & \leq O\left(\frac{1}{T}\right)L_{n_{T}}\left(\omega-\eta\right)+\left(O\left(\frac{1}{T}\right)+n_{T}^{-1}O\left(\left|\omega-\eta\right|\right)\right)L_{n_{T}}\left(\omega-\eta\right)\\
 & \leq C<\infty
\end{align*}
 where we have used the fact that $L_{n_{T}}\left(\omega-\eta\right)\leq n_{T}$
and 
\begin{align*}
\left|\omega-\eta\right|L_{n_{T}}\left(\omega-\eta\right) & =\begin{cases}
\left|\omega-\eta\right|n_{T}, & \left|\omega-\eta\right|\leq1/n_{T}\\
1, & 1/n_{T}\leq\left|\omega-\eta\right|\leq\pi
\end{cases}.
\end{align*}
Using Lemma \ref{Lemma A.5 Dahl97} and \eqref{Eq. (H_t)<L_t<L_n},
we have
\begin{align*}
n_{T}^{-1} & \left|\sum_{s=0}^{n_{T}-1}h\left(\frac{s}{T}\right)\exp\left(i\left(\omega-\eta\right)s\right)A_{\left\lfloor Tu\right\rfloor -n_{T}/2+s}^{0}\left(-\eta\right)d\eta\right|\\
 & =n_{T}^{-1}\left|A\left(\left(\left\lfloor Tu\right\rfloor -n_{T}/2\right)/T,\,-\eta\right)H_{n_{T}}\left(-\omega+\eta\right)\right|+O\left(T^{-1}\right)\\
 & =n_{T}^{-1}O\left(\sup_{u\in\left[0,\,1\right]}A\left(u,\,-\eta\right)\right)L_{n_{T}}\left(-\omega+\eta\right)+O\left(T^{-1}\right).
\end{align*}
Thus, after integration over $\eta$ we obtain that the error in replacing
$A_{\left\lfloor Tu\right\rfloor -n_{T}/2+t}^{0}\left(\eta\right)$
by $A\left(u,\,\omega\right)$ is $O\left(\left(\log n_{T}\right)/n_{T}\right)$.
Next, we replace $A_{\left\lfloor Tu\right\rfloor -n_{T}/2+s}^{0}\left(-\eta\right)$
by $A\left(u,\,\omega\right)$ and integrate over $\eta$ using the
relation 
\begin{align*}
A\left(u,\,\omega\right)A\left(u,\,-\omega\right) & =\left|A\left(u,\,\omega\right)\right|^{2}=f_{\mathbf{X}}\left(u,\,\omega\right).
\end{align*}
In view of 
\begin{align}
\int_{-\pi}^{\pi}\left|H_{n_{T}}\left(\alpha\right)\right|^{2}d\alpha & =2\pi\sum_{t=0}^{n_{T}-1}\left(\frac{t}{n_{T}}\right)^{2},\label{Eq. (HnT) =00003D integral}
\end{align}
 we then have 
\begin{alignat}{1}
\mathbb{E}\left(I_{h,T}\left(u,\,\omega\right)\right) & =\frac{1}{2\pi H_{2,n_{T}}\left(0\right)}\sum_{t=0}^{n_{T-1}}\sum_{s=0}^{n_{T-1}}h\left(\frac{t}{T}\right)\left(\frac{s}{T}\right)\int_{-\pi}^{\pi}\exp\left(-i\left(\omega-\alpha\right)\left(s-t\right)\right)f_{\mathbf{X}}\left(u,\,\alpha\right)d\alpha+O\left(\frac{\log n_{T}}{n_{T}}\right)\label{eq (1) bias}\\
 & =\frac{1}{\int_{-\pi}^{\pi}\left|H_{n_{T}}\left(\alpha\right)\right|^{2}d\alpha}\int_{-\pi}^{\pi}\left|H_{n_{T}}\left(\omega-\alpha\right)\right|^{2}f_{\mathbf{X}}\left(u,\,\alpha\right)d\alpha+O\left(\frac{\log n_{T}}{n_{T}}\right)\nonumber \\
 & =\frac{1}{\int_{-\pi}^{\pi}\left|H_{n_{T}}\left(\alpha\right)\right|^{2}d\alpha}\int_{-\pi}^{\pi}\left|H_{n_{T}}\left(\alpha\right)\right|^{2}f_{\mathbf{X}}\left(u,\,\omega-\alpha\right)d\alpha+O\left(\frac{\log n_{T}}{n_{T}}\right).\nonumber 
\end{alignat}
This shows the first equality of \eqref{Eq. Bias Local periodogram}.
For the second equality replace $A_{\left\lfloor Tu\right\rfloor -n_{T}/2+t}^{0}\left(\eta\right)$
by $A(u+\left(t-n_{T}/2\right)/T,\,\omega)$ and $A_{\left\lfloor Tu\right\rfloor -n_{T}/2+t}^{0}\left(-\eta\right)$
by $A(u+\left(t-n_{T}/2\right)/T,\,-\omega)$ so that \eqref{eq (1) bias}
holds with $f_{\mathbf{X}}(u+\left(t-n_{T}/2\right)/T,\,\omega)$
in place of $f_{\mathbf{X}}\left(u,\,\alpha\right)$. Then take a
second-order Taylor expansion of $f_{\mathbf{X}}$ around around $u$
to obtain 
\begin{align*}
\mathbb{E}\left(I_{h,T}\left(u,\,\omega\right)\right) & =\frac{1}{2\pi H_{2,n_{T}}\left(0\right)}\sum_{t=0}^{n_{T-1}}h\left(\frac{t}{T}\right)^{2}f_{\mathbf{X}}\left(u+\frac{t-n_{T}/2}{T},\,\omega\right)+O\left(\frac{\log n_{T}}{n_{T}}\right)\\
 & =f_{\mathbf{X}}\left(u,\,\omega\right)+\frac{1}{2}\left(\frac{n_{T}}{T}\right)^{2}\int_{0}^{1}x^{2}h^{2}\left(x\right)dx\frac{\partial^{2}}{\partial u^{2}}f_{\mathbf{X}}\left(u,\,\omega\right)\\
 & \quad+o\left(\left(\frac{n_{T}}{T}\right)^{2}\right)+O\left(\frac{\log n_{T}}{n_{T}}\right).\,\square
\end{align*}

\subsubsection{Proof of Theorem \ref{Theorem 5.2.8 Brillinger}}

By Theorem 2.3.1-(ix) in \citeReferencesSupp{brillinger:75}, $\mathrm{Cov}\left(Y_{j},\,Y_{k}\right)=\mathrm{cum}(Y_{j},\,\overline{Y}_{k})$
for possibly complex variables $Y_{j}$ and $Y_{k}.$ Thus,
\begin{align*}
\mathrm{Cov} & \left(d_{h,T}\left(u,\,\omega_{j}\right)d_{h,T}\left(u,\,-\omega_{j}\right),\,d_{h,T}\left(u,\,\omega_{k}\right)d_{h,T}\left(u,\,-\omega_{k}\right)\right)\\
 & =\mathrm{cum}\left(d_{h,T}\left(u,\,\omega_{j}\right)d_{h,T}\left(u,\,-\omega_{j}\right),\,d_{h,T}\left(u,\,\omega_{k}\right)d_{h,T}\left(u,\,-\omega_{k}\right)\right).
\end{align*}
 By the product theorem for cumulants {[}cf. \citeReferencesSupp{brillinger:75},
Theorem 2.3.2{]}, we have to sum over all indecompasable partitions
$\left\{ P_{1},\ldots,\,P_{m}\right\} $ with $\left|P_{i}\right|\mathrm{=card}(P_{i})\geq2$
of the two-way table,
\begin{center}
\begin{tabular}{|cc|}
$a_{j,1}$ & $a_{j,2}$\tabularnewline
$a_{k,1}$ & $a_{k,2}$\tabularnewline
\end{tabular}, 
\par\end{center}

\noindent where $a_{j,1}$ and $a_{j,2}$ stand for the positions
of $d_{h,T}\left(u,\,\omega_{j}\right)$ and $d_{h,T}\left(u,\,-\omega_{j}\right)$,
respectively.  This results in, 
\begin{align*}
\mathrm{cum} & \left(d_{h,T}\left(u,\,\omega_{j}\right)d_{h,T}\left(u,\,-\omega_{j}\right),\,d_{h,T}\left(u,\,\omega_{k}\right)d_{h,T}\left(u,\,-\omega_{k}\right)\right)\\
 & =\mathrm{cum}\left(d_{h,T}\left(-\omega_{j}\right),\,d_{h,T}\left(-\omega_{j}\right),\,d_{h,T}\left(\omega_{k}\right),\,d_{h,T}\left(-\omega_{k}\right)\right)\\
 & \quad+\mathrm{cum}\left(d_{h,T}\left(\omega_{j}\right)\right)\mathrm{cum}\left(d_{h,T}\left(-\omega_{j}\right),\,d_{h,T}\left(\omega_{k}\right),\,d_{h,T}\left(-\omega_{k}\right)\right)\\
 & \quad+\mathrm{three\,similar\,terms}\\
 & \quad+\mathrm{cum}\left(d_{h,T}\left(\omega_{j}\right)\right)\mathrm{cum}\left(d_{h,T}\left(\omega_{k}\right)\right)\mathrm{cum}\left(d_{h,T}\left(-\omega_{j}\right),\,d_{h,T}\left(-\omega_{k}\right)\right)\\
 & \quad+\mathrm{three\,similar\,terms}\\
 & \quad+\mathrm{cum}\left(d_{h,T}\left(\omega_{j}\right),\,d_{h,T}\left(-\omega_{k}\right)\right)\mathrm{cum}\left(d_{h,T}\left(-\omega_{j}\right),\,d_{h,T}\left(-\omega_{k}\right)\right)\\
 & \quad+\mathrm{cum}\left(d_{h,T}\left(\omega_{j}\right),\,d_{h,T}\left(-\omega_{k}\right)\right)\mathrm{cum}\left(d_{h,T}\left(-\omega_{j}\right),\,d_{h,T}\left(\omega_{k}\right)\right).
\end{align*}
Then, by Theorem \ref{Theorem 4.3.1-2 Brillinger}, 
\begin{align}
\mathrm{cum} & \left(d_{h,T}\left(u,\,\omega_{j}\right)d_{h,T}\left(u,\,-\omega_{j}\right),\,d_{h,T}\left(u,\,\omega_{k}\right)d_{h,T}\left(u,\,-\omega_{k}\right)\right)\label{Eq (Cum (dj, dj, dk, dk)}\\
 & =\left(2\pi\right)^{3}H_{4,n_{T}}\left(0\right)f_{\mathbf{X}}\left(u,\,\omega_{j},\,-\omega_{j},\,\omega_{k}\right)+O\left(1\right)\nonumber \\
 & \quad+\left[2\pi H_{2,n_{T}}\left(\omega_{j}+\omega_{k}\right)f_{\mathbf{X}}\left(u,\,\omega_{j}\right)+O\left(1\right)\right]\left[2\pi H_{2,n_{T}}\left(-\omega_{j}-\omega_{k}\right)f_{\mathbf{X}}\left(u,\,\omega_{j}\right)+O\left(1\right)\right]\nonumber \\
 & \quad+\left[2\pi H_{2,n_{T}}\left(\omega_{j}-\omega_{k}\right)f_{\mathbf{X}}\left(u,\,\omega_{j}\right)+O\left(1\right)\right]\left[2\pi H_{2,n_{T}}\left(-\omega_{j}+\omega_{k}\right)f_{\mathbf{X}}\left(u,\,\omega_{j}\right)+O\left(1\right)\right].\nonumber 
\end{align}
Given 
\begin{align*}
H_{2,n_{T}}\left(0\right) & =\sum_{t=0}^{n_{T}-1}h^{2}\left(t/T\right)\sim n_{T}\int h^{2}\left(\alpha\right)d\alpha
\end{align*}
 and 
\begin{align*}
H_{2,n_{T}}\left(\omega_{j}-\omega_{k}\right)H_{2,n_{T}}\left(-\omega_{j}+\omega_{k}\right) & =\left|H_{2,n_{T}}\left(\omega_{j}-\omega_{k}\right)\right|,^{2}
\end{align*}
the result of the theorem follows because
\begin{align*}
n_{T}^{-2}\left(2\pi\right)^{3}H_{4,n_{T}}\left(0\right)f_{\mathbf{X}}\left(u,\,\omega_{j},\,-\omega_{j},\,\omega_{k}\right) & =O\left(n_{T}^{-1}\right),
\end{align*}
 and because the $O\left(1\right)$ terms on the right-hand side of
\eqref{Eq (Cum (dj, dj, dk, dk)} become negligible when multiplied
by $H_{2,n_{T}}^{-2}\left(0\right)$. 

Next, we prove the second result of the theorem. Recall that $\mathbf{z}\sim\mathscr{N}_{p}^{\mathrm{C}}\left(\mu_{z},\,\Sigma_{z}\right)$
means that the $2p$ vector 
\begin{align*}
\begin{bmatrix}\mathrm{Re\,\mathbf{z}}\\
\mathrm{Im}\,\mathbf{z}
\end{bmatrix}
\end{align*}
 is distributed as
\begin{align*}
\mathscr{N}_{2p}\left(\begin{bmatrix}\mathrm{Re\,}\mu_{z}\\
\mathrm{Im}\,\mu_{z}
\end{bmatrix},\,\frac{1}{2}\begin{bmatrix}\mathrm{Re\,}\Sigma_{z} & -\mathrm{Im\,}\Sigma_{z}\\
-\mathrm{Im\,}\Sigma_{z} & \mathrm{Re\,}\Sigma_{z}
\end{bmatrix}\right) & ,
\end{align*}
 where $\Sigma_{z}$ is a $p\times p$ hermitian positive semidefinite
matrix. By Theorem \ref{Theorem 4.4.1-2 Brillinger} we know that
$\mathrm{Re}\,\mathbf{d}_{h,T}\left(\omega_{j}\right)$ and $\mathrm{Im}\,\mathbf{d}_{h,T}\left(\omega_{j}\right)$
are asymptotically independent $\mathscr{N}\left(0,\,\pi n_{T}f_{\mathbf{X}}\left(u,\,\omega_{j}\right)\right)$
variates. Hence, by the Mann-Wald Theorem,
\begin{align*}
I_{h,T}\left(u,\,\omega_{j}\left(n_{T}\right)\right) & =\left(2\pi n_{T}\right)^{-1}\left\{ \left(\mathrm{Re\,}d_{h,T}\left(u,\,\omega_{j}\left(n_{T}\right)\right)\right)^{2}+\left(\mathrm{Im}\,d_{h,T}\left(\omega_{j}\left(n_{T}\right)\right)\right)^{2}\right\} 
\end{align*}
is asymptotically distributed as $f_{\mathbf{X}}\left(u,\,\omega_{j}\right)\chi_{2}^{2}/2$
if $2\omega_{j}\not\equiv0\,(\mathrm{mod\,}2\pi)$. This proves part
(i). For part (ii), if $\omega=\pm\pi,\,\pm3\pi,\ldots,$ then $I_{h,T}\left(u,\,\omega\right)$
is asymptotically distributed as $f_{\mathbf{X}}\left(u,\,\omega\right)\chi_{1}^{2}$,
independently from the previous variates. $\square$

\subsubsection{Proof of Theorem \ref{Theorem 5.6.1}}

Using Theorem \ref{Theorem 5.2.3 Brillinger}, we have
\begin{align*}
\mathbb{E}\left(f_{h,T}\left(u,\,\omega\right)\right) & =\frac{2\pi}{n_{T}}\sum_{s=0}^{n_{T}-1}W_{T}\left(\omega-\frac{2\pi s}{n_{T}}\right)\mathbb{E}\left(I_{h,T}\left(u,\,\frac{2\pi s}{n_{T}}\right)\right)\\
 & =\frac{2\pi}{n_{T}}\sum_{s=0}^{n_{T}-1}W_{T}\left(\omega-\frac{2\pi s}{n_{T}}\right)f_{\mathbf{X}}\left(u,\,\frac{2\pi s}{T}\right)+O\left(n_{T}T^{-1}\right)+O\left(\log\left(n_{T}\right)n_{T}^{-1}\right).
\end{align*}
The first term on the right-hand side is 
\begin{align*}
\frac{2\pi}{n_{T}} & \sum_{s=0}^{n_{T}-1}W_{T}\left(\omega-\frac{2\pi s}{n_{T}}\right)f_{\mathbf{X}}\left(u,\,\frac{2\pi s}{n_{T}}\right)\\
 & =\int_{0}^{2\pi}W_{T}\left(\omega-\alpha\right)f_{\mathbf{X}}\left(u,\,\alpha\right)d\alpha+O\left(\left(n_{T}b_{T}\right)^{-1}\right)+O\left(\log\left(n_{T}\right)n_{T}^{-1}\right)\\
 & =\int_{0}^{2\pi}\sum_{j=-\infty}^{\infty}b_{T}^{-1}W\left(b_{T}^{-1}\left(\omega-\alpha+2\pi j\right)\right)f_{\mathbf{X}}\left(u,\,\alpha\right)d\alpha+O\left(\left(n_{T}b_{T}\right)^{-1}\right)+O\left(\log\left(n_{T}\right)n_{T}^{-1}\right)\\
 & =\int_{-\infty}^{\infty}W\left(\beta\right)f_{\mathbf{X}}\left(u,\,\omega-\beta b_{T}\right)d\beta+O\left(\left(n_{T}b_{T}\right)^{-1}\right)+O\left(\log\left(n_{T}\right)n_{T}^{-1}\right),
\end{align*}
 where the last equality follows from the change in variable $\beta=b_{T}^{-1}\left(\omega-\alpha\right)$.
This yields the first equality of \eqref{Eq. (5.6.7)}. The second
equality follows from the first and Theorem \ref{Theorem 5.2.3 Brillinger}
along with a Taylor expansion. $\square$ 

\subsubsection{Proof of Theorem \ref{Theorem 5.6.4 Brillinger}}

Let
\begin{align*}
c_{T}\left(u,\,k\right) & =H_{2,T}\left(0\right)^{-1}\sum_{s=0}^{n_{T}-1}h\left(\frac{s+k}{T}\right)h\left(\frac{s}{T}\right)X_{\left\lfloor Tu\right\rfloor -n_{T}/2+s+k+1,T}X_{\left\lfloor Tu\right\rfloor -n_{T}/2+s+1,T}.
\end{align*}
We can rewrite $I_{h,T}\left(u,\,\omega\right)$ using $c_{T}\left(u,\,k\right)$
as follows, 
\begin{align*}
I_{h,T}\left(u,\,\omega\right)=\left(2\pi\right)^{-1}\sum_{k=-\infty}^{\infty}\exp\left(-i\omega k\right)c_{T}\left(u,\,k\right) & .
\end{align*}
 Note that
\begin{align*}
f_{h,T}\left(u,\,\omega\right) & =\int_{0}^{2\pi}W_{2,T}\left(\omega-\alpha\right)I_{h,T}\left(u,\,\alpha\right)d\alpha+O\left(\left(n_{T}b_{W,T}\right)^{-1}\right),
\end{align*}
 where $W_{2,T}\left(\omega\right)=\sum_{k=-\infty}^{\infty}w\left(b_{W,T}k\right)\exp\left(-i\omega k\right)$
and $w\left(k\right)=\int_{-\infty}^{\infty}W_{2,T}\left(\alpha\right)\exp\left(i\alpha k\right)d\alpha$
for $k\in\mathbb{R}$. From Theorem \ref{Theorem 5.2.8 Brillinger},
\begin{align*}
\mathrm{Cov} & \left(f_{h,T}\left(u,\,\omega_{j}\right),\,f_{h,T}\left(u,\,\omega_{k}\right)\right)\\
 & =\int_{0}^{2\pi}\int_{0}^{2\pi}W_{2,T}\left(\omega_{j}-\alpha\right)W_{2,T}\left(\omega_{k}-\beta\right)\mathrm{Cov}\left(I_{h,T}\left(u,\,\alpha\right),\,I_{h,T}\left(u,\,\beta\right)\right)d\alpha d\beta\\
 & =H_{2,n_{T}}\left(0\right)^{-1}H_{2,n_{T}}\left(0\right)^{-1}\int_{0}^{2\pi}\int_{0}^{2\pi}W_{2,T}\left(\omega_{j}-\alpha\right)W_{2,T}\left(\omega_{k}-\beta\right)\\
 & \quad\times\{\left|H_{2,n_{T}}\left(\alpha-\beta\right)\right|^{2}+\left|H_{2,n_{T}}\left(\alpha+\beta\right)\right|^{2}\}\left|f\left(u,\,\alpha\right)\right|^{2}d\alpha d\beta+O\left(n_{T}^{-1}\right).
\end{align*}
 We now show that 
\begin{align}
\int_{0}^{2\pi}W_{2,T}\left(\omega_{k}-\beta\right) & \left|H_{2,n_{T}}\left(\alpha-\beta\right)\right|^{2}d\beta\label{eq (**) Brillinger Ch. 7}\\
 & =2\pi W_{2,T}\left(\omega_{k}-\alpha\right)\sum_{s=0}^{n_{T}-1}h^{4}\left(s\right)+O\left(b_{W,T}^{-2}\right),\nonumber 
\end{align}
 uniformly in $\alpha.$  We can expand \eqref{eq (**) Brillinger Ch. 7}
as follows, 
\begin{align*}
\sum_{t=0}^{n_{T}-1}\sum_{s=0}^{n_{T}-1} & h^{2}\left(t/n_{T}\right)h^{2}\left(s/n_{T}\right)\int_{0}^{2\pi}W_{2,T}\left(\omega_{k}-\beta\right)\times\exp\left\{ -i\left(\alpha-\beta\right)t+i\left(\alpha-\beta\right)s\right\} d\beta\\
 & =\sum_{t=0}^{n_{T}-1}\sum_{s=0}^{n_{T}-1}h^{2}\left(t\right)h^{2}\left(s\right)\int_{0}^{2\pi}\sum_{k=-\infty}^{\infty}w\left(b_{W,T}k\right)\exp\left(-i\left(\omega_{k}-\beta\right)k\right)\\
 & \quad\times\exp\left\{ -i\left(\alpha-\beta\right)t+i\left(\alpha-\beta\right)s\right\} d\beta\\
 & =\sum_{t=0}^{n_{T}-1}\sum_{s=0}^{n_{T}-1}h^{2}\left(t/n_{T}\right)h^{2}\left(s/n_{T}\right)w\left(b_{W,T}\left(t-s\right)\right)\exp\left(i\left(\omega_{k}-\alpha\right)\left(t-s\right)\right)\\
 & =\sum_{k=-\infty}^{\infty}w\left(b_{W,T}k\right)\exp\left(i\left(\omega_{k}-\alpha\right)k\right)\sum_{s=0}^{n_{T}-1}h^{2}\left(\left(s+k\right)n_{T}\right)h^{2}\left(s/n_{T}\right)\\
 & =2\pi W_{2,T}\left(\omega_{k}-\alpha\right)\sum_{s=0}^{n_{T}-1}h^{4}\left(s/n_{T}\right)+R_{T},
\end{align*}
where we have applied Lemma \ref{Lemma P4.1 Brillinger} to $\exp\left(i\left(\omega_{k}-\alpha\right)k\right)\sum_{s=0}^{n_{T}-1}h^{2}\left(s+k\right)h^{2}\left(s\right)$
to yield,
\begin{align*}
\left|\exp\left(i\left(\omega_{k}-\alpha\right)k\right)\sum_{s=0}^{n_{T}-1}h^{2}\left(s+k\right)h^{2}\left(s\right)-\exp\left(i\left(\omega_{k}-\alpha\right)k\right)\sum_{s=0}^{n_{T}-1}h^{4}\left(s/n_{T}\right)\right| & \leq C\left|k\right|,
\end{align*}
and 
\begin{align*}
\left|R_{T}\right| & \leq C\sum_{k=-\infty}^{\infty}\left|w\left(b_{W,T}k\right)\right|\left|k\right|\sim Cb_{W,T}^{-2}\int\left|x\right|\left|w\left(x\right)\right|dx,
\end{align*}
for $0<C<\infty$. The latter result follows because 
\begin{align*}
C\sum_{k=-\infty}^{\infty} & \left|w\left(b_{T}k\right)\right|\left|k\right|\\
 & =Cb_{W,T}^{-2}b_{W,T}\sum_{k=-\infty}^{\infty}\left|w\left(b_{W,T}k\right)\right|\left|b_{W,T}k\right|\\
 & =Cb_{T}^{-2}\int\left|x\right|\left|w\left(x\right)\right|dx,
\end{align*}
 for a finite $0<C<\infty.$ A similar result holds for the second
term involving $\left|H_{2,n_{T}}\left(\alpha+\beta\right)\right|^{2}$.
Overall, we have
\begin{align*}
\mathrm{Cov} & \left(f_{h,T}\left(u,\,\omega_{j}\right),\,f_{h,T}\left(u,\,\omega_{k}\right)\right)\\
 & =2\pi H_{2,T}\left(0\right)^{-2}\sum_{s=0}^{n_{T}-1}h\left(s/n_{T}\right)^{4}\int_{0}^{2\pi}\{W_{2,T}\left(\omega_{j}-\alpha\right)W_{2,T}\left(\omega_{k}-\alpha\right)\left|f\left(u,\,\alpha\right)\right|^{2}\\
 & \quad+W_{2,T}\left(\omega_{j}-\alpha\right)W_{2,T}\left(\omega_{k}+\alpha\right)\left|f\left(u,\,\alpha\right)\right|^{2}\}d\alpha+O\left(b_{W,T}^{-2}n_{T}^{-2}\right)+O\left(n_{T}^{-1}\right).
\end{align*}
Equation \eqref{Eq. (5.6.19)} follows from
\begin{align*}
n_{T}b_{W,T} & \mathrm{Cov}\left(f_{h,T}\left(u,\,\omega_{j}\right),\,f_{h,T}\left(u,\,\omega_{k}\right)\right)\\
 & =b_{W,T}2\pi n_{T}H_{2,T}\left(0\right)^{-1}n_{T}H_{2,T}\left(0\right)^{-1}n_{T}^{-1}\sum_{t=0}^{n_{T}-1}h\left(t/n_{T}\right)^{4}\\
 & \quad\times\int_{0}^{2\pi}\biggl\{\sum_{l=-\infty}^{\infty}b_{W,T}^{-1}W\left(b_{W,T}^{-1}\left(\omega_{j}-\alpha+2\pi l\right)\right)\\
 & \quad\times\sum_{l=-\infty}^{\infty}b_{W,T}^{-1}W\left(b_{W,T}^{-1}\left(\omega_{k}-\alpha+2\pi l\right)\right)\left|f\left(u,\,\alpha\right)\right|^{2}\\
 & \quad+\sum_{l=-\infty}^{\infty}b_{W,T}^{-1}W\left(b_{W,T}^{-1}\left(\omega_{j}-\alpha+2\pi l\right)\right)\\
 & \quad\times\sum_{l=-\infty}^{\infty}b_{W,T}^{-1}W\left(b_{W,T}^{-1}\left(\omega_{k}+\alpha+2\pi l\right)\right)\left|f\left(u,\,\alpha\right)\right|^{2}\}d\alpha+O\left(\left(n_{T}b_{W,T}\right)^{-1}\right)+O\left(b_{W,T}\right)\\
 & =2\pi\left(\int h^{2}\left(t\right)dt\right)^{-2}\int h^{4}\left(t\right)dt\\
 & \quad\int_{0}^{2\pi}\left[\eta\left\{ \omega_{j}-\omega_{k}\right\} \left|f\left(u,\,\omega_{j}\right)\right|^{2}+\eta\left\{ \omega_{j}+\omega_{k}\right\} \left|f\left(u,\,\omega_{j}\right)\right|^{2}\right]\int_{-\infty}^{\infty}W^{2}\left(\alpha\right)d\alpha\\
 & \quad+O\left(\left(n_{T}b_{W,T}\right)^{-1}\right)+O\left(b_{W,T}\right).
\end{align*}
Finally, we consider the magnitude of the joint cumulants of order
$r.$  We have 
\begin{align}
\mathrm{cum} & \left(f_{h,T}\left(u,\,\omega_{1}\right),\ldots,\,\,f_{h,T}\left(u,\,\omega_{r}\right)\right)\label{Eq. (*) p. 442 Brillinger}\\
 & =2\pi\left\{ H_{2,n_{T}}\left(0\right)\right\} ^{-r}\nonumber \\
 & \times\sum_{t_{1}=0}^{n_{T}-1}\cdots\sum_{t_{2r}=0}^{n_{T}-1}w\left(b_{T}\left(t_{1}-t_{2}\right)\right)\cdots w\left(b_{T}\left(t_{2r-1}-t_{2r}\right)\right)\nonumber \\
 & \times\exp\left(-i\omega_{1}\left(t_{1}-t_{2}\right)-\ldots-i\omega_{r}\left(t_{2r-1}-t_{2r}\right)\right)h_{n_{T}}\left(t_{1}\right)\cdots h_{n_{T}}\left(t_{2r}\right)\nonumber \\
 & \times\mathrm{cum}\,(X_{\left\lfloor Tu\right\rfloor -n_{T}/2+t_{1}+1,T}X_{\left\lfloor Tu\right\rfloor -n_{T}/2+t_{2}+1,T},\ldots.,\nonumber \\
 & \quad X_{\left\lfloor Tu\right\rfloor -n_{T}/2+t_{2r-1}+1,T}X_{\left\lfloor Tu\right\rfloor -n_{T}/2+t_{2r}+1,T}).\nonumber 
\end{align}
 Note that
\begin{align*}
\mathrm{cum} & \left(X_{\left\lfloor Tu\right\rfloor -n_{T}/2+t_{1}+1,T}X_{\left\lfloor Tu\right\rfloor -n_{T}/2+t_{2}+1,T},\ldots,\,X_{\left\lfloor Tu\right\rfloor -n_{T}/2+t_{2r-1}+1,T}X_{\left\lfloor Tu\right\rfloor -n_{T}/2+t_{2r}+1,T}\right)\\
 & =\sum_{\mathbf{v}}c_{X\cdots X}\left(u;\,t_{j},\,j\in v_{1}\right)\cdots c_{X\cdots X}\left(u;\,t_{j},\,j\in v_{\mathrm{p}}\right),
\end{align*}
 where $c_{X\cdots X}\left(u;\,t_{j},\,j\in v_{1}\right)$ is the
time-$Tu$ cumulant involving the variables $X_{t_{j}}$ for $j\in v_{1}$
and where the summation is over all indecomposable partitions $\mathbf{v}=\left(v_{1},\ldots,\,v_{\mathrm{P}}\right)$
of the table
\begin{center}
\begin{tabular}{|cc|}
1 & 2\tabularnewline
3 & 4\tabularnewline
$\vdots$ & $\vdots$\tabularnewline
$2r-1$ & $2r$\tabularnewline
\end{tabular}.
\par\end{center}

As the partition is indecomposable, in each set $v_{\mathrm{p}}$
of the partition we may find an element $t_{\mathrm{p}}^{*}$ such
that none of $t_{j}-t_{\mathrm{p}}^{*}$, $j\in v_{\mathrm{p}}$ $\left(\mathrm{p}=1,\ldots,\,\mathrm{P}\right)$
is $t_{2l-1}-t_{2l},$ $l=1,\,2,\ldots,\,r$. Define $2r-\mathrm{P}$
new variables $k_{1},\ldots,\,k_{2r-\mathrm{P}}$ as the nonzero $t_{j}-t_{\mathrm{p}}^{*}$.
Eq. \eqref{Eq. (*) p. 442 Brillinger} is now bounded by
\begin{align*}
C^{r}n_{T}^{-r} & \sum_{\mathbf{v}}\sum_{t_{1}^{*}}\cdots\sum_{t_{\mathrm{P}}^{*}}\sum_{k_{1}}\cdots\sum_{k_{2r-\mathrm{P}}}\biggl|w\left(b_{W,T}\left(k_{\alpha_{1}}+t_{\beta_{1}}^{*}-k_{\alpha_{1}}-t_{\beta_{2}}^{*}\right)\right)\\
 & \quad\cdots\times w\left(b_{W,T}\left(k_{\alpha_{2r-1}}+t_{\beta_{2r-1}}^{*}-k_{\alpha_{2r}}-t_{\beta_{2r}}^{*}\right)\right)\biggr|\\
 & \quad\times\left|h\left(t_{1}^{*}/n_{T}\right)\right|^{2r}\left|c_{X\cdots X}\left(u;\,k_{1},\ldots\right)\cdots c_{X\cdots X}\left(u;\,\ldots,\,k_{2r-\mathrm{P}}\right)\right|,
\end{align*}
 for some finite $C$, where $\alpha_{1},\ldots,\,\alpha_{2r}$ are
selected from $1,\ldots,\,2r$ and $\beta_{1},\ldots,\,\beta_{2r}$
from $1,\ldots,\,\mathrm{P}$. By Lemma 2.3.1 in \citeReferencesSupp{brillinger:75},
there are $\mathrm{P}-1$ linearly independent differences among the
$t_{\beta_{1}^{*}}-t_{\beta_{2}^{*}},\ldots,\,t_{\beta_{2r-1}^{*}}-t_{\beta_{2r}^{*}}$.
Suppose these are $t_{\beta_{1}^{*}}-t_{\beta_{2}^{*}},\ldots,\,t_{\beta_{2r-2}^{*}}-t_{\beta_{2r-1}^{*}}.$
Making the change of variables
\begin{align*}
s_{1} & =k_{\alpha_{1}}+t_{\beta_{1}}^{*}-k_{\alpha_{1}}-t_{\beta_{2}}^{*}\\
 & \;\vdots\\
s_{\mathrm{P}-1} & =k_{\alpha_{2\mathrm{P}-3}}+t_{\beta_{2\mathrm{P}-3}}^{*}-k_{\alpha_{2\mathrm{P}-2}}-t_{\beta_{2\mathrm{P}-2}}^{*},
\end{align*}
the cumulant \eqref{Eq. (*) p. 442 Brillinger} is bounded by
\begin{align*}
C^{r}n_{T}^{-r}\sum_{\mathbf{v}} & \sum_{t_{1}^{*}}\sum_{s_{1}}\cdots\sum_{s_{\mathrm{P}-1}}\sum_{k_{1}}\cdots\sum_{k_{2r-\mathrm{P}}}\left|w\left(b_{W,T}s_{1}\right)\cdots w\left(b_{W,T}s_{P-1}\right)\right|\\
 & \quad\left|h\left(t_{1}^{*}/n_{T}\right)\right|^{2r}\left|c_{X\cdots X}\left(u;\,k_{1},\ldots\right)\cdots c_{X\cdots X}\left(u;\,\ldots,\,k_{2r-\mathrm{P}}\right)\right|\\
 & \leq C^{r}n_{T}^{-r+1}b_{W,T}^{-\left(\mathrm{P}-1\right)}\sum_{\mathbf{v}}C_{n_{2,1}}\cdots C_{n_{2,\mathrm{P}}}\\
 & =O\left(n_{T}^{-r+1}b_{W,T}^{-\left(\mathrm{P}-1\right)}\right),
\end{align*}
where $\mathrm{P}\leq r$ and $C_{n_{2},j}=\sup_{u\in\left[0,\,1\right]}\sum_{t_{1},\ldots,\,t_{n_{2},j}}\left|c_{X\cdots X}\left(u;\,t_{1},\ldots,\,t_{n_{2},j}\right)\right|$
with $n_{2,j}$ denoting the number of elements in the $j$th set
of the partition $\mathbf{v}$. It follows that for $r>2$, 
\begin{align*}
\mathrm{cum}\left(\left(n_{T}b_{W,T}\right)^{1/2}f_{h,T}\left(u,\,\omega_{1}\right),\ldots,\,\left(n_{T}b_{W,T}\right)^{1/2}f_{h,T}\left(u,\,\omega_{r}\right)\right) & \rightarrow0.
\end{align*}
Thus, the variates $f_{h,T}\left(u,\,\omega_{1}\right),\ldots,\,f_{h,T}\left(u,\,\omega_{r}\right)$
are asymptotically normal with the moment structure given in the theorem.
$\square$

\subsection{\label{subsec:Proof-of- Tests and Limiting distribution}Proof of
the Results of Section \ref{Section Change-Point-Tests}}

\subsubsection{Preliminary Lemmas}

Let $I_{T}^{*}\left(j/T,\,\omega\right)=I_{L,h,T}\left(j/T,\,\omega\right)-\mathbb{E}(I_{L,h,T}(j/T,\,\omega))$.
For $w\geq0$ consider the dependence measure,
\begin{align}
\phi_{I,w,q} & =\sup_{j\in\left\{ \mathbf{S}_{r};\,r=1,\ldots,\,M_{T}-2\right\} }\left\Vert I_{h,T}^{*}\left(j/T,\,\omega\right)-I_{h,T,\left\{ w\right\} }^{*}\left(j/T,\,\omega\right)\right\Vert _{q},\label{Eq. (2.1) WZ Delta-1}
\end{align}
where  $X_{t,T}$ in $I_{h,T,\left\{ w\right\} }^{*}\left(j/T,\,\omega\right)$
is replaced by $X_{t,T,\left\{ w\right\} }$.  Let $\Upsilon_{n,q}=\sum_{j=n}^{\infty}\phi_{I,j,q}.$ 
\begin{lem}
\label{Lemma: GMC for I}Let Assumption \ref{Assumption Locally Stationary}-\ref{Assumption WZ 2011}
hold. We have $I_{h,T}^{*}\left(j/T,\,\omega\right)\in\mathscr{L}^{q}$
and, for $q>2$, $\Upsilon_{n,q}=O\left(n^{-\gamma}\right)$ for some
$\gamma>0$. 
\end{lem}
\noindent\textit{Proof of Lemma }\ref{Lemma: GMC for I}\textit{.}
We have\textit{ 
\begin{align*}
\phi_{I,w,q} & =\sup_{j\in\left\{ \mathbf{S}_{r};\,r=1,\ldots,\,M_{T}-2\right\} }\Biggl\Vert\frac{1}{2\pi H_{2,n_{T}}\left(0\right)}\sum_{s=0}^{n_{T}-1}\sum_{t=0}^{n_{T}-1}h\left(\frac{s}{n_{T}}\right)h\left(\frac{t}{n_{T}}\right)\\
 & \quad\left(X_{j-n_{T}+s+1,T}X_{j-n_{T}+t+1,T}-X_{j-n_{T}+s+1,T,\left\{ w\right\} }X_{j-n_{T}+t+1,T,\left\{ w\right\} }\right)\exp\left(-i\omega\left(s-t\right)\right)\Biggl\Vert_{q}.
\end{align*}
}Note that 
\begin{align}
\Biggl\Vert\frac{1}{2\pi H_{2,n_{T}}\left(0\right)}\sum_{s=0}^{n_{T}-1}\sum_{t=0}^{n_{T}-1} & h\left(\frac{s}{n_{T}}\right)h\left(\frac{t}{n_{T}}\right)\left(X_{j-n_{T}+s+1,T}\left(X_{j-n_{T}+t+1,T}-X_{j-n_{T}+t+1,T,\left\{ w\right\} }\right)\right)\Biggr\Vert_{q}\nonumber \\
 & \leq\frac{1}{2\pi H_{2,n_{T}}\left(0\right)}\sum_{s=0}^{n_{T}-1}\sum_{t=0}^{n_{T}-1}\int_{-\pi}^{\pi}\left|h\left(\frac{s}{n_{T}}\right)h\left(\frac{t}{n_{T}}\right)\exp\left(-i\omega\left(s-t\right)\right)\right|\nonumber \\
 & \quad\times\left\Vert X_{j-n_{T}+s+1,T}\right\Vert _{q}\left\Vert X_{j-n_{T}+t+1,T}-X_{j-n_{T}+t+1,T,\left\{ w\right\} }\right\Vert _{q}\nonumber \\
 & \leq C\frac{1}{2\pi H_{2,n_{T}}\left(0\right)}\sum_{s=0}^{n_{T}-1}\sum_{t=0}^{n_{T}-1}\int_{-\pi}^{\pi}\left|h\left(\frac{s}{n_{T}}\right)h\left(\frac{t}{n_{T}}\right)\exp\left(-i\omega\left(s-t\right)\right)\right|\nonumber \\
 & \quad\times\sup_{j}\left\Vert X_{j-n_{T}+t+1,T}-X_{j-n_{T}+t+1,T,\left\{ w\right\} }\right\Vert _{q}\nonumber \\
 & \leq C\frac{1}{2\pi H_{2,n_{T}}\left(0\right)}\sum_{s=0}^{n_{T}-1}\sum_{t=0}^{n_{T}-1}\int_{-\pi}^{\pi}\left|h\left(\frac{s}{n_{T}}\right)h\left(\frac{t}{n_{T}}\right)\exp\left(-i\omega\left(s-t\right)\right)\right|\phi_{w,q}.\label{Eq. A.1 in Lemma GMC}
\end{align}
By Lemma A.7 in \citeReferencesSupp{dahlhaus:96}, we have $\sum_{s=0}^{n_{T}-1}|h\left(s/n_{T}\right)\exp\left(-i\omega s\right)|\leq Cn_{T}^{-1}L_{n_{T}}\left(\omega\right)^{2}$
for some $C<\infty$. Using Lemma \ref{Lemma A.4 Dahl97}-(iii),
the right-hand side of \eqref{Eq. A.1 in Lemma GMC} is less than
or equal to
\begin{align*}
 & C\frac{K^{2}n_{T}^{-2}n_{T}^{3}}{2\pi H_{2,n_{T}}\left(0\right)}\phi_{w,q}\leq C_{2}\phi_{w,q},
\end{align*}
 for some $C_{2}<\infty$. Overall, we obtain $\phi_{I,w,q}\leq C_{2}\phi_{w,q}$
and so $\Upsilon_{n,q}\leq C_{2}\sum_{j=n}^{\infty}\phi_{j,q}.$
Using Assumption \ref{Assumption WZ 2011} we have $\Upsilon_{n,q}=C_{2}O\left(n^{-q+1}\right)$
since $(\sum_{j=n}^{\infty}\phi_{j,q}^{2})^{1/2}\leq\sum_{j=n}^{\infty}\phi_{j,q}$
and $\sum_{n=0}^{\infty}n^{q-1}\phi_{n,q}<\infty$. $\square$ 

\medskip{}

The result in Lemma \ref{Lemma: GMC for I} also holds for $I_{T}^{*}\left(j/T,\,\omega\right)$
constructed using $I_{R,h,T}\left(j/T,\,\omega\right)$ in place of
$I_{L,h,T}\left(j/T,\,\omega\right)$. Let $\tau_{T}=T^{\vartheta_{1}}\left(\log\left(T\right)\right)^{\vartheta_{2}}$
where $\vartheta_{1}=\left(1/2-1/q+\gamma/q\right)$ $/\left(1/2-1/q+\gamma\right)$
and $\vartheta_{2}=\left(\gamma+\gamma/q\right)/$ $\left(1/2-1/q+\gamma\right)$
for some $\gamma>0$. 

The proofs below can be simplified by noting that $\widehat{\sigma}_{L,r}^{2}\left(\omega\right)$
is a consistent estimate of $\sigma_{L,r}^{2}\left(\omega\right)=\mathrm{Var}(\sqrt{M_{S,T}}\widetilde{f}_{L,r,T}^{*}\left(\omega\right))$
where 
\begin{align*}
\widetilde{f}_{L,r,T}^{*}\left(\omega\right) & =M_{S,T}^{-1}\sum_{j\in\mathbf{S}_{r}}f_{L,h,T}^{*}\left(j/T,\,\omega\right)
\end{align*}
 with $f_{L,h,T}^{*}\left(j/T,\,\omega\right)=f_{L,h,T}\left(j/T,\,\omega\right)$
$-\mathbb{E}(f_{L,h,T}(j/T,\,\omega))$. The  consistency result
follows from results in \citet{casini_hac}. The rate of convergence
of $\widehat{\sigma}_{L,r}^{2}\left(\omega\right)$ is $O(\sqrt{\widetilde{M}_{S,T}b_{1,T}}).$
Given Assumption \ref{Taper, nT, W and K1, b1}-(iv), $(\widetilde{M}_{S,T}b_{1,T})^{-1/2}M_{S,T}^{1/2}\log T\rightarrow0$
and so one can replace $\widehat{\sigma}_{L,r}\left(\omega\right)$
by $\sigma_{L,r}\left(\omega\right)$ in the definition of $\mathrm{S}_{\mathrm{max},T}\left(\omega\right)$
and $\mathrm{S}_{\mathrm{Dmax},T}$ throughout the proofs of Lemma
\ref{Lemma Smax - Smax_old 0 } and of Theorem \ref{Theorem Asymptotic H0 Distrbution Smax}-\ref{Theorem Asymptotic H0 Distrbution SDmax}. 
\begin{lem}
\label{Lemma Smax - Smax_old 0 }Let Assumption \ref{Assumption Locally Stationary}-\ref{Assumption Minimum f}
and Condition \ref{Condition n_T h_T} hold. Under $\mathcal{H}_{0}$,
$\sqrt{\log\left(M_{T}\right)}M_{S,T}^{1/2}(\mathrm{S}_{\mathrm{max},T}\left(\omega\right)-\mathrm{\widetilde{S}}_{\mathrm{max},T}\left(\omega\right))\overset{\mathbb{P}}{\rightarrow}0$
for any $\omega\in\left[-\pi,\,\pi\right]$ where 
\begin{align*}
\mathrm{\widetilde{S}}_{\mathrm{max},T}\left(\omega\right) & \triangleq\max_{r=1,\ldots,\,M_{T}-2}\left|\frac{\widetilde{f}_{r,T}\left(\omega\right)-\widetilde{f}_{r+1,T}\left(\omega\right)}{\sigma_{f,r}\left(\omega\right)}\right|.
\end{align*}
\end{lem}
\noindent\textit{Proof of Lemma }\ref{Lemma Smax - Smax_old 0 }\textit{.
}Note that for arbitrary sequences of numbers $\left(a_{i}\right)_{i=1,\ldots,\,N}$
and $\left(b_{i}\right)_{i=1,\ldots,\,N}$ with $N\geq1$, we have
for any $i,$ 
\begin{align}
\left|a_{i}\right| & \leq\left|a_{i}-b_{i}\right|+\left|b_{i}\right|\leq\max_{i=1,\ldots,N}\left|a_{i}-b_{i}\right|+\max_{i=1,\ldots,N}\left|b_{i}\right|.\label{Eq. Ineqaulity ai bi}
\end{align}
The inequality still holds if on the left-hand side we replace $\left|a_{i}\right|$
by $\max_{i=1,\ldots,N}\left|a_{i}\right|$. We then have 
\begin{align}
\mathrm{S}_{\mathrm{max},T} & \left(\omega\right)-\mathrm{\widetilde{S}}_{\mathrm{max},T}\left(\omega\right)\nonumber \\
 & =\max_{r=1,\ldots,\,M_{T}-2}\left|\frac{\widetilde{f}_{L,r,T}\left(\omega\right)-\widetilde{f}_{R,r+1,T}\left(\omega\right)}{\sigma_{L,r}\left(\omega\right)}\right|-\max_{r=1,\ldots,\,M_{T}-2}\left|\frac{\widetilde{f}_{r,T}\left(\omega\right)-\widetilde{f}_{r+1,T}\left(\omega\right)}{\sigma_{f,r}\left(\omega\right)}\right|.\label{Eq. Smax -Smax old}
\end{align}
Using \eqref{Eq. Ineqaulity ai bi} the right-hand side of \eqref{Eq. Smax -Smax old}
is less than or equal to
\begin{align*}
\max_{r=1,\ldots,\,M_{T}-2} & \left|\frac{\widetilde{f}_{L,r,T}\left(\omega\right)-\widetilde{f}_{R,r+1,T}\left(\omega\right)}{\sigma_{L,r}\left(\omega\right)}-\frac{\widetilde{f}_{r,T}\left(\omega\right)-\widetilde{f}_{r+1,T}\left(\omega\right)}{\sigma_{L,r}\left(\omega\right)}\right|\\
 & +\max_{r=1,\ldots,\,M_{T}-2}\left|\frac{\widetilde{f}_{r,T}\left(\omega\right)-\widetilde{f}_{r+1,T}\left(\omega\right)}{\sigma_{L,r}\left(\omega\right)}\right|-\max_{r=1,\ldots,\,M_{T}-2}\left|\frac{\widetilde{f}_{r,T}\left(\omega\right)-\widetilde{f}_{r+1,T}\left(\omega\right)}{\sigma_{f,r}\left(\omega\right)}\right|.
\end{align*}
 The second line converges to zero in probability given the uniform
asymptotic equivalence of $\sigma_{L,r}\left(\omega\right)$ and $\sigma_{f,r}\left(\omega\right)$
with an error $O\left(T^{-1}\right)$. Thus, it is sufficient to show
\begin{align*}
\max_{r=1,\ldots,\,M_{T}-2} & \sqrt{\log\left(M_{T}\right)}M_{S,T}^{1/2}\left|\frac{\widetilde{f}_{L,r,T}\left(\omega\right)-\widetilde{f}_{R,r+1,T}\left(\omega\right)-\left(\widetilde{f}_{r,T}\left(\omega\right)-\widetilde{f}_{r+1,T}\left(\omega\right)\right)}{\sigma_{L,r}\left(\omega\right)}\right|\overset{\mathbb{P}}{\rightarrow}0.
\end{align*}
We use the following decomposition,
\begin{align}
\sqrt{\log\left(M_{T}\right)}M_{S,T}^{1/2} & \left|\frac{\widetilde{f}_{L,r,T}\left(\omega\right)-\widetilde{f}_{R,r+1,T}\left(\omega\right)-\left(\widetilde{f}_{r,T}\left(\omega\right)-\widetilde{f}_{r+1,T}\left(\omega\right)\right)}{\sigma_{L,r}\left(\omega\right)}\right|\nonumber \\
 & \leq\sqrt{\log\left(M_{T}\right)}M_{S,T}^{1/2}\left|\frac{\widetilde{f}_{L,r,T}\left(\omega\right)-\widetilde{f}_{r,T}\left(\omega\right)}{\sigma_{L,r}\left(\omega\right)}\right|\label{Eq.fL -f R}\\
 & \quad+\sqrt{\log\left(M_{T}\right)}M_{S,T}^{1/2}\left|\frac{\left(\widetilde{f}_{R,r+1,T}\left(\omega\right)-\widetilde{f}_{r+1,T}\left(\omega\right)\right)}{\sigma_{L,r}\left(\omega\right)}\right|.\nonumber 
\end{align}
 Let us consider the first term on the right-hand side of \eqref{Eq.fL -f R}.
Note that $\widetilde{f}_{L,r,T}\left(\omega\right)$ and $\widetilde{f}_{r,T}\left(\omega\right)$
are weighted averages of  stochastic $\theta$-H\"older continuous
variables each standardized by $\sigma_{L,r}\left(\omega\right).$
These variables standardized by $\sigma_{L,r}\left(\omega\right)$
belong to the same block $r$ of the sample. Thus, their weighted
average is asymptotically approximated by a stochastic variable that
satisfies stochastic $\theta$-H\"older continuity as for $X_{t,T}$
in \eqref{Eq. Stochastic Lip cont}, where this property holds uniformly
in $r$. Thus, it follows that
\begin{align*}
\sigma_{L,r}^{-1}\left(\omega\right)\sqrt{\log\left(M_{T}\right)} & M_{S,T}^{1/2}\left|\widetilde{f}_{L,r,T}\left(\omega\right)-\widetilde{f}_{r,T}\left(\omega\right)\right|\\
 & =\sqrt{\log\left(M_{T}\right)}M_{S,T}^{1/2}\left(O\left(\left(m_{T}/T\right)^{\theta}\right)+O_{\mathbb{P}}\left(\left(n_{T}b_{W,T}\right)^{-1/2}\right)\right)\\
 & =o_{\mathbb{P}}\left(1\right),
\end{align*}
where the $O_{\mathbb{P}}((n_{T}b_{W,T})^{-1/2})$ rate follows from
Theorem \ref{Theorem 5.6.4 Brillinger} and the last equality uses
Condition \ref{Condition n_T h_T}. Using Markov's inequality, this
shows that for all $\epsilon>0$, 
\begin{align*}
\mathbb{P}\left(\max_{r=1,\ldots,\,M_{T}-2}\sqrt{\log\left(M_{T}\right)}M_{S,T}^{1/2}\left|\widetilde{f}_{L,r,T}\left(\omega\right)-\widetilde{f}_{r,T}\left(\omega\right)\right|>\epsilon\right) & \rightarrow0.
\end{align*}
 The argument for the second term of \eqref{Eq.fL -f R} is analogous.
$\square$
\begin{lem}
\label{Lemma Rmax - Rmax_old 0 }Let Assumption \ref{Assumption Locally Stationary}-\ref{Assumption Minimum f}
and Condition \ref{Condition n_T h_T} hold. Under $\mathcal{H}_{0}$,
we have $\sqrt{\log\left(M_{T}\right)}M_{S,T}^{1/2}(\mathrm{R}_{\mathrm{max},T}\left(\omega\right)-\mathrm{\widetilde{R}}_{\mathrm{max},T}\left(\omega\right))\overset{\mathbb{P}}{\rightarrow}0$
for any $\omega\in\left[-\pi,\,\pi\right]$, where 
\begin{align*}
\mathrm{\widetilde{R}}_{\mathrm{max},T}\left(\omega\right) & \triangleq\max_{r=1,\ldots,\,M_{T}-2}\left|\frac{\widetilde{f}_{r,T}\left(\omega\right)}{\widetilde{f}_{r+1,T}\left(\omega\right)}-1\right|.
\end{align*}
\end{lem}
\noindent\textit{Proof of Lemma }\ref{Lemma Rmax - Rmax_old 0 }\textit{.
}Using \eqref{Eq. Ineqaulity ai bi}, we have 
\begin{align}
\left|\mathrm{R}_{\mathrm{max},T}\left(\omega\right)-\mathrm{\widetilde{R}}_{\mathrm{max},T}\left(\omega\right)\right| & \leq\max_{r=1,\ldots,\,M_{T}-2}\left|\frac{\widetilde{f}_{L,r,T}\left(\omega\right)}{\widetilde{f}_{R,r+1,T}\left(\omega\right)}-1-\left(\frac{\widetilde{f}_{r,T}\left(\omega\right)}{\widetilde{f}_{r+1,T}\left(\omega\right)}-1\right)\right|\nonumber \\
 & \leq\max_{r=1,\ldots,\,M_{T}-2}\left|\widetilde{f}_{L,r,T}\left(\omega\right)\left(\frac{1}{\widetilde{f}_{R,r+1,T}\left(\omega\right)}-\frac{1}{\widetilde{f}_{r+1,T}\left(\omega\right)}\right)\right|\label{Eq. Rmax -Rmax old}\\
 & \quad+\max_{r=1,\ldots,\,M_{T}-2}\left|\frac{\widetilde{f}_{L,r,T}\left(\omega\right)-\widetilde{f}_{r,T}\left(\omega\right)}{\widetilde{f}_{r+1,T}\left(\omega\right)}\right|.\nonumber 
\end{align}
 Let us consider the second term on the right-hand side of \eqref{Eq. Rmax -Rmax old}.
Note that for all $\epsilon>0$ and all constants $C>0,$ we have
\begin{align}
\mathbb{P} & \left(\max_{r=1,\ldots,\,M_{T}-2}\left|\sqrt{\log\left(M_{T}\right)}M_{S,T}^{1/2}\frac{\widetilde{f}_{L,r,T}\left(\omega\right)-\widetilde{f}_{r,T}\left(\omega\right)}{\widetilde{f}_{r+1,T}\left(\omega\right)}\right|>\epsilon\right)\nonumber \\
 & \leq\mathbb{P}\left(\max_{r=1,\ldots,\,M_{T}-2}\sqrt{\log\left(M_{T}\right)}M_{S,T}^{1/2}\left|\widetilde{f}_{L,r,T}\left(\omega\right)-\widetilde{f}_{r,T}\left(\omega\right)\right|\cdot\max_{r=1,\ldots,\,M_{T}-2}\left|\frac{1}{\widetilde{f}_{r+1,T}\left(\omega\right)}\right|>\epsilon\right)\nonumber \\
 & \leq\mathbb{P}\left(\max_{r=1,\ldots,\,M_{T}-2}\sqrt{\log\left(M_{T}\right)}M_{S,T}^{1/2}\left|\widetilde{f}_{L,r,T}\left(\omega\right)-\widetilde{f}_{r,T}\left(\omega\right)\right|>\frac{\epsilon}{C}\right)\label{Eq. (69) BJV-1}\\
 & \quad+\leq\mathbb{P}\left(\max_{r=1,\ldots,\,M_{T}-2}\left|\frac{1}{\widetilde{f}_{r+1,T}\left(\omega\right)}\right|>C\right).\nonumber 
\end{align}
Theorem \ref{Theorem 5.6.1} implies that
\begin{align*}
\mathbb{E}\left(f_{h,T}\left(u,\,\omega\right)\right) & =f\left(u,\,\omega\right)+O\left(\left(n_{T}/T\right)^{-2}\right)+O\left(\log\left(n_{T}\right)n_{T}^{-1}\right)+o\left(b_{W,T}^{2}\right).
\end{align*}
 The same result holds for $f_{L,h,T}\left(u,\,\omega\right)$. By
Assumption \ref{Assumption Smothness of A} we have 
\begin{align}
f\left(\left(\left(r+1\right)m_{T}+j\right)/T,\,\omega\right)-f\left(\left(rm_{T}+j\right)/T,\,\omega\right) & =O\left(\left(m_{T}/T\right)^{\theta}\right),\qquad\mathrm{uniformly\,in}\,r\,\mathrm{and}\,j.\label{Eq. Uniform Lip f}
\end{align}
 Thus, using the bound for the variance of $f_{h,T}\left(u,\,\omega\right)$
in Theorem \ref{Theorem 5.6.4 Brillinger} and Condition \ref{Condition n_T h_T},
we have 

\begin{align*}
\max_{r=1,\ldots,\,M_{T}-2}\sqrt{\log\left(M_{T}\right)} & M_{S,T}^{1/2}\left|\widetilde{f}_{L,r,T}\left(\omega\right)-\widetilde{f}_{r,T}\left(\omega\right)\right|\\
 & =\sqrt{\log\left(M_{T}\right)}M_{S,T}^{1/2}\left(O\left(\left(m_{T}/T\right)^{\theta}\right)+O_{\mathbb{P}}\left(\left(n_{T}b_{W,T}\right)^{-1/2}\right)\right)\\
 & =o_{\mathbb{P}}\left(1\right).
\end{align*}
By using Markov's inequality, this shows that 
\begin{align}
\mathbb{P}\left(\max_{r=1,\ldots,\,M_{T}-2}\sqrt{\log\left(M_{T}\right)}M_{S,T}^{1/2}\left|\widetilde{f}_{L,r,T}\left(\omega\right)-\widetilde{f}_{r,T}\left(\omega\right)\right|>\frac{\epsilon}{C}\right) & \rightarrow0.\label{Eq. before (69)}
\end{align}
By Theorem \ref{Theorem 5.6.4 Brillinger}, $\widetilde{f}_{r+1,T}\left(\omega\right)=f\left(\left(r+1\right)m_{T},\,\omega\right)+o_{\mathbb{P}}\left(1\right)$.
Thus, the second term of \eqref{Eq. (69) BJV-1} also converges to
zero for example by choosing $C=3/f_{-}.$ Altogether we obtain that
the right-hand side of \eqref{Eq. (69) BJV-1} converges to zero.
Next, we consider the first term of \eqref{Eq. Rmax -Rmax old}. For
any $\epsilon>0$ and any $C>0,$ we have
\begin{align}
\mathbb{P} & \left(\max_{r=1,\ldots,\,M_{T}-2}\sqrt{\log\left(M_{T}\right)}M_{S,T}^{1/2}\left|\widetilde{f}_{L,r,T}\left(\omega\right)\left(\frac{1}{\widetilde{f}_{R,r+1,T}\left(\omega\right)}-\frac{1}{\widetilde{f}_{r+1,T}\left(\omega\right)}\right)\right|>\epsilon\right)\nonumber \\
 & \leq\mathbb{P}\left(\max_{r=1,\ldots,\,M_{T}-2}\sqrt{\log\left(M_{T}\right)}M_{S,T}^{1/2}\widetilde{f}_{L,r,T}\left(\omega\right)\left|\widetilde{f}_{r+1,T}\left(\omega\right)-\widetilde{f}_{R,r+1,T}\left(\omega\right)\right|>\frac{\epsilon}{C}\right)\label{Eq. (69) BJV-1-1}\\
 & \quad+\mathbb{P}\left(\max_{r=1,\ldots,\,M_{T}-2}\left|\frac{1}{\widetilde{f}_{r+1,T}\left(\omega\right)\widetilde{f}_{R,r+1,T}\left(\omega\right)}\right|>C\right).\nonumber 
\end{align}
The first term on the right-hand side above is less than or equal
to,
\begin{align*}
\mathbb{P} & \left(\max_{r=1,\ldots,\,M_{T}-2}\left|\widetilde{f}_{L,r,T}\left(\omega\right)\right|>C_{2}\right)\\
 & \quad+\mathbb{P}\left(\max_{r=1,\ldots,\,M_{T}-2}\sqrt{\log\left(M_{T}\right)}M_{S,T}^{1/2}\left|\widetilde{f}_{r+1,T}\left(\omega\right)-\widetilde{f}_{R,r+1,T}\left(\omega\right)\right|>\frac{\epsilon}{C\cdot C_{2}}\right),
\end{align*}
 for all $C_{2}>0.$ We can choose $C_{2}$ large enough such that
the first term above converges to zero. The second term above converges
to zero by the same argument as in \eqref{Eq. before (69)}. The second
term on the right-hand side of \eqref{Eq. (69) BJV-1-1} can be expanded
as follows, 
\begin{align*}
\mathbb{P} & \left(\max_{r=1,\ldots,\,M_{T}-2}\left|\frac{1}{\widetilde{f}_{r+1,T}\left(\omega\right)\widetilde{f}_{R,r+1,T}\left(\omega\right)}\right|>C\right)\\
 & \leq\mathbb{P}\left(\min_{r=1,\ldots,\,M_{T}-2}\left|\widetilde{f}_{R,r+1,T}\left(\omega\right)\right|<C^{-1/2}\right)+\mathbb{P}\left(\min_{r=1,\ldots,\,M_{T}-2}\left|\widetilde{f}_{r+1,T}\left(\omega\right)\right|<C^{-1/2}\right)\\
 & \leq\mathbb{P}\left(\min_{r=1,\ldots,\,M_{T}-2}\left|\widetilde{f}_{R,r+1,T}\left(\omega\right)\right|<C^{-1/2}\right)+\mathbb{P}\left(\min_{r=1,\ldots,\,M_{T}-2}\left|\widetilde{f}_{R,r+1,T}\left(\omega\right)\right|<2C^{-1/2}\right)\\
 & \quad+\mathbb{P}\left(\max_{r=1,\ldots,\,M_{T}-2}\left|\widetilde{f}_{R,r+1,T}\left(\omega\right)-\widetilde{f}_{r+1,T}\left(\omega\right)\right|>2C^{-1/2}\right).
\end{align*}
 The first two terms on the right-hand side have already been discussed
above. The third term has also been discussed above with the multiplicative
factor $\sqrt{\log\left(M_{T}\right)}M_{S,T}^{1/2}$. $\square$ 

\subsubsection{Proof of Lemma \ref{Lemma Theorem 4.1 in Shao and Wu (2007)}}

Let $J_{t}=\kappa_{\mathbf{X},t}\left(k_{1},\cdots,\,k_{r-1}\right)$
where $t\leq t+k_{1}\leq\ldots\leq t+k_{r-1}$. For $1\leq l\leq r-1$,
let $n_{l}=k_{l}-k_{l-1}$. Define the vector  $y_{0}=y_{0,l}=(k_{1}-k_{l-1},\cdots,\,k_{l-2}-k_{l-1},\,0).$
Let $\mathcal{G}_{n}=\left(\ldots,\,e'_{n-1},\,e'_{n}\right)$ where
$n\in\mathbb{N}$ and for $t\geq0$ define $X_{t,T}^{*}=H\left(t/T,\,\{\mathcal{G}_{0},\,e_{1},\cdots,\,e_{t}\}\right)$.
Following Proposition 2 in \citeReferencesSupp{wu/zhao:04} and by
the additivity of cumulants, 

\begin{align}
\sup_{t}\left|J_{t}\right| & =\sup_{k_{l-1},\,l\in\left\{ 1,\ldots,\,r-1\right\} }\left|\kappa_{\mathbf{X},-k_{l-1}}\left(y_{0},\,k_{l}-k_{l-1},\,k_{l+1}-k_{l-1},\cdots,\,k_{r-1}-k_{l-1}\right)\right|\label{Eq (Sup Jt)}\\
 & =\left|\kappa_{\mathbf{X},-k_{l_{*}-1}}\left(y_{0}^{*},\,k_{l_{*}}-k_{l_{*}-1},\,k_{l_{*}+1}-k_{l_{*}-1},\cdots,\,k_{r-1}-k_{l_{*}-1}\right)\right|\nonumber \\
 & \leq\sum_{j=0}^{r-l_{*}-1}\biggl|\mathrm{cum}\biggl(Y_{0},\,X_{k_{l_{*}}-k_{l_{*}-1},T}^{*},\cdots,\,X_{k_{l_{*}+j-1}-k_{l_{*}-1},T}^{*},\nonumber \\
 & \quad X_{k_{l_{*}+j}-k_{l_{*}-1},T}-X_{k_{l_{*}+j}-k_{l_{*}-1},T}^{*},\,X_{k_{l_{*}+j+1}-k_{l_{*}-1},T},\,\cdots,\,X_{k_{k-1}-k_{l_{*}-1},T}\biggr)\biggr|\nonumber \\
 & =\sum_{j=0}^{r-l_{*}-1}B_{j}\left(l_{*}\right),\nonumber 
\end{align}
 where $y_{0}^{*}=y_{0,l_{*}}^{*}=(k_{1}-k_{l_{*}-1},\cdots,\,k_{l_{*}-2}-k_{l_{*}-1},\,0)$,
and $Y_{0}=Y_{0,l_{*}}=(X_{k_{1}-k_{l_{*}-1},T},\cdots,\,X_{k_{l_{*}-2}-k_{l_{*}-1},T},\,X_{0,T})$.
Let $\zeta_{t}=||X_{t,T}-X_{j,T}^{*}||_{r}$ and $S_{j}=S_{j,r}=\sum_{i=j}^{\infty}\phi_{i,r}^{2}.$
By Proposition 2 in \citeReferencesSupp{wu/zhao:04} we have $|B_{j}\left(l_{*}\right)|\leq C_{1}\zeta_{k_{l_{*}+j}-k_{l_{*}-1}}$,
where $C_{1}$ only depends on $r$ and on $\sup_{t}\mathbb{E}|X_{t,T}|^{m}$
with $1\leq m\leq r$. Thus, using \eqref{Eq (Sup Jt)} and Proposition
2 of \citeReferencesSupp{wu:07}, we have
\begin{align*}
\sup_{t}|J_{t}| & \leq C_{1}\sum_{j=0}^{r-l_{*}-1}\zeta_{k_{l_{*}+j}-k_{j-1}}\leq C_{2}\sum_{j=0}^{r-l_{*}-1}S_{k_{l_{*}+j}-k_{j-1}}^{1/2}\leq C_{3}S_{n_{l_{*}}}^{1/2},
\end{align*}
where $C_{2}=18r^{3/2}\left(r-1\right)^{-1/2}C_{1}$ and $C_{3}=C_{2}r$.
Since $1\leq l_{*}\leq r-1$, we have $\sup_{t}|J_{t}|\leq C_{3}S_{n_{l_{*}}}^{1/2}$.
Then, 
\begin{align*}
\sum_{k_{1},\ldots,\,k_{r-1}=-\infty}^{\infty}\left(1+|k_{j}|^{l}\right)\sup_{t}\left|\kappa_{\mathbf{X},t}\left(k_{1},\cdots,\,k_{r-1}\right)\right| & \leq2r!\sum_{k_{1},\ldots,\,k_{r-1}=0}^{\infty}|k_{r-1}|^{l}\sup_{t}\left|\kappa_{\mathbf{X},t}\left(k_{1},\cdots,\,k_{r-1}\right)\right|\\
 & \leq2r!\sum_{k_{1},\ldots,\,k_{r-1}=0}^{\infty}|k_{r-1}|^{l}C_{3}S_{n_{l_{*}}}^{1/2}\\
 & \leq2C_{3}r!\,\sum_{n=0}^{\infty}n^{l+r-2}S_{n}^{1/2}\\
 & <\infty,
\end{align*}
where the last inequality follows from Assumption \ref{Assumption WZ 2011}.
$\square$

\subsubsection{Proof of Theorem \ref{Theorem Asymptotic H0 Distrbution Smax}}

From Lemma \ref{Lemma Smax - Smax_old 0 } it is sufficient to show
the result for $\mathrm{\widetilde{S}}_{\mathrm{max},T}\left(\omega\right)$
since the latter is asymptotically equivalent to $\mathrm{S}_{\mathrm{max},T}\left(\omega\right).$
Define $f_{h,T}^{*}\left(j/T,\,\omega\right)=f_{h,T}\left(j/T,\,\omega\right)-\mathbb{E}\left(f_{h,T}\left(j/T,\,\omega\right)\right)$.
For $\omega\in\left[-\pi,\,\pi\right]$ let $S_{r+1}\left(\omega\right)=\sum_{j\in\left\{ \mathbf{S}_{s},\,s=1,\ldots,\,r+1\right\} }f_{h,T}^{*}\left(j/T,\,\omega\right)$
and 
\begin{align*}
R_{r,T}\left(\omega\right) & =\frac{1}{M_{S,T}}\left(S_{r+1}\left(\omega\right)-\sum_{j\in\left\{ \mathbf{S}_{s},\,s=1,\ldots,\,r+1\right\} }\mathscr{W}_{j}\left(\omega\right)-\left(S_{r}\left(\omega\right)-\sum_{j\in\left\{ \mathbf{S}_{s},\,s=1,\ldots,\,r\right\} }\mathscr{W}_{j}\left(\omega\right)\right)\right),
\end{align*}
where $\mathscr{W}_{j}\left(\omega\right)=\sigma_{j}\left(\omega\right)Z_{j}$
with $Z_{j}\sim\mathrm{i.i.d}.\,\mathscr{N}\left(0,\,1\right)$. Write
\begin{align}
\widetilde{f}_{r,T}\left(\omega\right) & =M_{S,T}^{-1}\sum_{j\in\mathbf{S}_{r}}f_{h,T}\left(j/T,\,\omega\right)\label{Eq(IB)-1}\\
 & =M_{S,T}^{-1}\sum_{j\in\mathbf{S}_{r}}\left(f_{h,T}^{*}\left(j/T,\,\omega\right)+\mathbb{E}\left(f_{h,T}\left(j/T,\,\omega\right)\right)\right)\nonumber \\
 & =\frac{1}{M_{S,T}}\sum_{j\in\mathbf{S}_{r}}\mathscr{W}_{j}\left(\omega\right)+R_{r,T}+\frac{1}{M_{S,T}}\sum_{j\in\mathbf{S}_{r}}\mathbb{E}\left(f_{h,T}\left(j/T,\,\omega\right)\right).\nonumber 
\end{align}
Under Assumption \ref{Assumption WZ 2011}, Theorem 1 in \citeReferencesSupp{wu/zhou:11}
yields $\max_{0\leq r\leq M_{S,T}-1}\left|R_{r,T}\right|=O_{\mathbb{P}}\left(\tau_{T}/M_{S,T}\right)$.
By Theorem \ref{Theorem 5.6.1}, 
\begin{align*}
\mathbb{E}\left(f_{h,T}\left(j/T,\,\omega\right)\right) & =f\left(j/T,\,\omega\right)+O\left(\left(n_{T}/T\right)^{2}\right)+O\left(b_{W,T}^{2}\right)+O\left(\log\left(n_{T}\right)/n_{T}\right).
\end{align*}
 Using \eqref{Eq. Uniform Lip f} we obtain 
\begin{align}
\sqrt{M_{S,T}} & \left(\widetilde{f}_{r+1,T}\left(\omega\right)-\widetilde{f}_{r,T}\left(\omega\right)\right)\label{Eq(32) WZ (07)-1}\\
 & =\frac{1}{\sqrt{M_{S,T}}}\left(\sum_{j\in\mathbf{S}_{r+1}}\mathscr{W}_{j}\left(\omega\right)-\sum_{j\in\mathbf{S}_{r}}\mathscr{W}_{j}\left(\omega\right)\right)\nonumber \\
 & \quad+O\left(M_{S,T}^{1/2}m_{T}^{\theta}/T^{\theta}\right)+O_{\mathbb{P}}\left(\tau_{T}/M_{S,T}^{1/2}\right)+O_{\mathbb{P}}\left(M_{S,T}^{1/2}\left(n_{T}/T\right)^{2}+M_{S,T}^{1/2}b_{W,T}^{2}+M_{S,T}^{1/2}\log\left(n_{T}\right)/n_{T}\right)\nonumber \\
 & =\frac{1}{\sqrt{m_{T}}}\left(\sum_{j\in\mathbf{S}_{r+1}}\mathscr{W}_{j}\left(\omega\right)-\sum_{j\in\mathbf{S}_{r}}\mathscr{W}_{j}\left(\omega\right)\right)\nonumber \\
 & \quad+o_{\mathbb{P}}\left(\left(\log M_{T}\right)^{-1/2}\right).\nonumber 
\end{align}
The result then follows from Lemma 1 in \citeReferencesSupp{wu/zhao:07}.
$\square$ 

\subsubsection{Proof of Theorem \ref{Theorem Asymptotic H0 Distrbution SDmax}}
\begin{lem}
\label{Lemma Extrem Value =00003D> Local Periodogram}Let $\mathscr{V}\left(\omega\right)$
denote a random variable defined by $\mathbb{P}\left(\mathscr{V}\left(\omega\right)\leq v\right)=\exp(-\pi^{-1/2}\exp(-v))$
for $\omega\in\Pi$. Assume that for $\omega,\,\omega'\in\Pi$ the
variables $\mathscr{V}\left(\omega\right)$ and $\mathscr{V}\left(\omega'\right)$
are independent. Let $\mathscr{V}^{*}\triangleq\max_{\omega\in\Pi}\mathscr{V}\left(\omega\right)-\log\left(n_{\omega}\right)$.
Then, $\mathbb{P}\left(\mathscr{V}^{*}\leq v\right)=\exp(-\pi^{-1/2}\exp\left(-v\right))$.
\end{lem}
\noindent\textit{Proof. }Since $\mathscr{V}\left(\omega\right)$
is independent from any $\mathscr{V}\left(\omega'\right)$ with $\omega\neq\omega'$,
we have 
\begin{align*}
\log\mathbb{P}\left(\mathscr{V}^{*}\leq v\right) & =\sum_{j=1}^{n_{\omega}}\log\mathbb{P}\left(\mathscr{V}\left(\omega_{j}\right)\leq\left(\log\left(n_{\omega}\right)+v\right)\right)\\
 & =\sum_{j=1}^{n_{\omega}}\left(-\pi^{-1/2}\exp\left(\log\left(n_{\omega}^{-1}\right)\right)\exp\left(-v\right)\right)\\
 & =-\pi^{-1/2}\exp\left(-v\right).
\end{align*}
Thus, $\mathbb{P}\left(\mathscr{V}^{*}\leq v\right)=\exp(-\pi^{-1/2}\exp\left(-v\right))$.
$\square$ \medskip{}

\noindent\textit{Proof of Theorem }\ref{Theorem Asymptotic H0 Distrbution SDmax}\textit{.
}From Theorem \ref{Theorem 5.6.4 Brillinger}, it follows that $f_{h,T}\left(u,\,\omega_{j}\right)$
and $f_{h,T}\left(u,\,\omega_{k}\right)$ are asymptotically independent
if $\omega_{k}\pm\omega_{k}\not\equiv0\,(\mathrm{mod\,}2\pi),\,1\leq j<k\leq n_{\omega}$.
The result then follows from Lemma \ref{Lemma Smax - Smax_old 0 }
and \ref{Lemma Extrem Value =00003D> Local Periodogram}, and Theorem
\ref{Theorem Asymptotic H0 Distrbution Smax}. $\square$ 

\subsubsection{Proof of Theorem \ref{Theorem Asymptotic H0 Distrbution Rmax and RDmax}}

Due to the self-normalization nature of the test statistic, we can
use Lemma \ref{Lemma Rmax - Rmax_old 0 } and steps similar to Proposition
A1-A.3 in \citeReferencesSupp{bibinger/jirak/vetter:16} to show that
it is sufficient to consider the behavior of 
\begin{align*}
\mathrm{\widetilde{R}}^{*}\left(\omega\right) & =\max_{r=1,\ldots,\,M_{T}-2}\left|M_{S,T}^{-1}\sum_{j\in\mathbf{S}_{r}}g_{T}\left(j/T,\,\omega\right)-M_{S,T}^{-1}\sum_{j\in\mathbf{S}_{r+1}}g_{T}\left(j/T,\,\omega\right)\right|,
\end{align*}
where $g_{T}\left(j/T,\,\omega\right)$ are random variables with
mean $\mathbb{E}\left(f_{h,T}\left(j/T,\,\omega\right)\right)$, unit
variance and satisfying Assumption \ref{Assumption WZ 2011}. For
$\omega\in\left[-\pi,\,\pi\right]$ let $S_{r+1}\left(\omega\right)=\sum_{j\in\left\{ \mathbf{S}_{s},\,s=1,\ldots,\,r+1\right\} }\left(g_{T}\left(j/T,\,\omega\right)-\mathbb{E}\left(f_{h,T}\left(j/T,\,\omega\right)\right)\right)$
and 
\begin{align*}
R_{r,T}\left(\omega\right) & =\frac{1}{M_{S,T}}\left(S_{r+1}\left(\omega\right)-\sum_{j\in\left\{ \mathbf{S}_{s},\,s=1,\ldots,\,r+1\right\} }\mathscr{W}_{j}\left(\omega\right)-\left(S_{r}\left(\omega\right)-\sum_{j\in\left\{ \mathbf{S}_{s},\,s=1,\ldots,\,r\right\} }\mathscr{W}_{j}\left(\omega\right)\right)\right),
\end{align*}
where $\mathscr{W}_{j}\left(\omega\right)=Z_{j}$ with $Z_{j}\sim\mathrm{i.i.d.}\,\mathscr{N}\left(0,\,1\right)$.
Write
\begin{align}
M_{S,T}^{-1}\sum_{j\in\mathbf{S}_{r}}g_{T}\left(j/T,\,\omega\right) & =M_{S,T}^{-1}\sum_{j\in\mathbf{S}_{r}}\left(\left(g_{T}\left(j/T,\,\omega\right)-\mathbb{E}\left(g_{T}\left(j/T,\,\omega\right)\right)\right)+\mathbb{E}\left(g_{T}\left(j/T,\,\omega\right)\right)\right)\label{Eq(IB)-1-1}\\
 & =\frac{1}{M_{S,T}}\sum_{j\in\mathbf{S}_{r}}\mathscr{W}_{j}\left(\omega\right)+R_{r,T}+\frac{1}{M_{S,T}}\sum_{j\in\mathbf{S}_{r}}\mathbb{E}\left(g_{T}\left(j/T,\,\omega\right)\right).\nonumber 
\end{align}
As in the proof of Theorem \ref{Theorem Asymptotic H0 Distrbution Smax},
we have $\max_{0\leq r\leq M_{S,T}-1}\left|R_{r,T}\right|=O_{\mathbb{P}}\left(\tau_{T}/M_{S,T}\right)$.
By Theorem \ref{Theorem 5.6.1}, $\mathbb{E}\left(g_{T}\left(j/T,\,\omega\right)\right)=f\left(j/T,\,\omega\right)+O\left(\left(n_{T}/T\right)^{2}\right)+O\left(b_{W,T}^{2}\right)+O\left(\log\left(n_{T}\right)/n_{T}\right)$.
Using \eqref{Eq. Uniform Lip f}, we obtain 
\begin{align}
\sqrt{M_{S,T}} & \left(M_{S,T}^{-1}\sum_{j\in\mathbf{S}_{r+1}}g_{T}\left(j/T,\,\omega\right)-M_{S,T}^{-1}\sum_{j\in\mathbf{S}_{r}}g_{T}\left(j/T,\,\omega\right)\right)\label{Eq(32) WZ (07)-1-1}\\
 & =\frac{1}{\sqrt{M_{S,T}}}\left(\sum_{j\in\mathbf{S}_{r+1}}\mathscr{W}_{j}\left(\omega\right)-\sum_{j\in\mathbf{S}_{r}}\mathscr{W}_{j}\left(\omega\right)\right)\nonumber \\
 & \quad+O\left(M_{S,T}^{1/2}m_{T}^{\theta}/T^{\theta}\right)+O_{\mathbb{P}}\left(\tau_{T}/M_{S,T}^{1/2}\right)+O_{\mathbb{P}}\left(M_{S,T}^{1/2}\left(n_{T}/T\right)^{2}+M_{S,T}^{1/2}b_{W,T}^{2}+M_{S,T}^{1/2}\log\left(n_{T}\right)/n_{T}\right)\nonumber \\
 & =\frac{1}{\sqrt{m_{T}}}\left(\sum_{j\in\mathbf{S}_{r+1}}\mathscr{W}_{j}\left(\omega\right)-\sum_{j\in\mathbf{S}_{r}}\mathscr{W}_{j}\left(\omega\right)\right)\nonumber \\
 & \quad+o_{\mathbb{P}}\left(\left(\log M_{T}\right)^{-1/2}\right).\nonumber 
\end{align}
The result about $\mathrm{R}_{\mathrm{max},T}\left(\omega\right)$
follows from Lemma 1 in \citet{wu/zhao:07}. The result concerning
$\mathrm{R}_{\mathrm{Dmax},T}$ follows using the same argument as
in the proof of Theorem \ref{Theorem Asymptotic H0 Distrbution SDmax}.
$\square$ 

\subsection{\label{subsec:Proofs-of-the Consistency}Proofs of the Results in
Section \ref{Section Consistency-and-Minimax}}

For a sequence of random variables $\{\xi_{j}\},$ let $\mathbb{P}_{\{\xi_{j}\}}$
denote the law of the observations $\{\xi_{j}\}$. Let $||\mathbb{P}_{\{\xi_{j}\}}-\mathbb{P}_{\{\xi_{j}^{*}\}}||_{\mathrm{TV}}$
define the total variation distance between the probability measures
$\mathbb{P}_{\{\xi_{j}\}}$ and $\mathbb{P}_{\{\xi_{j}^{*}\}}.$ For
two random variables $Y$ and $X$ with distributions $\mathbb{P}_{Y}$
and $\mathbb{P}_{X}$, respectively, denote the Kullback-Leibler divergence
by $D_{\mathrm{KL}}\left(Y||\,X\right)=D_{\mathrm{KL}}\left(\mathbb{P}_{Y}||\,\mathbb{P}_{X}\right)=\int\log\left(d\mathbb{P}_{Y}/d\mathbb{P}_{X}\right)d\mathbb{P}_{Y}$. 

\subsubsection{Proof of Theorem \ref{Theorem 4.1 BJV}}

The proof is based on several steps of information-theoretic reductions
that allow us to show the asymptotic equivalence in the strong Le
Cam sense of our statistical problem to a special high-dimensional
signal detection problem. The minimax lower bound is then obtained
by using classical arguments as in \citeReferencesSupp{ingster/suslina:03}.
Information-theoretic reductions were also used by \citeReferencesSupp{bibinger/jirak/vetter:16}
to establish a minimax lower bound for change-point testing in volatility
in the context of high-frequency data. Our derivations differ from
theirs in several ways because we deal with serially correlated observations
while they had independent observations. Furthermore, our testing
problem is more complex because our observations have an unknown distribution
while their observations are squared  standard normal variables.

We first consider alternatives as in $\mathcal{H}_{1}^{\mathrm{B}}$.
Throughout the proof we set 
\begin{align}
m_{T} & =C_{T}\left(\sqrt{\log\left(M_{T}\right)}T^{\theta}/D\right)^{\frac{2}{2\theta+1}},\label{Eq. (85) BJV 17}
\end{align}
with a constant $C_{T}>0$. We begin by granting the experimenter
additional knowledge thereby focusing on a simpler sub-model. This
additional knowledge can only decrease the lower bound on minimax
distinguishability and therefore such lower bound carries over to
the original model. We restrict attention to a sub-class of $\boldsymbol{F}_{1,\lambda_{b}^{0},\omega_{0}}\left(\theta,\,b_{T},\,D\right)$
which is characterized by a break at time $\lambda_{b}^{0}\in\left(0,\,1\right)$
with $\left|f\left(\lambda_{b}^{0},\,\omega_{0}\right)-f\left(\lambda_{b}^{0}+,\,\omega_{0}\right)\right|\geq b_{T}$,
where $f\left(\lambda_{b}^{0}+,\,\omega\right)=\lim_{s\downarrow\lambda_{b}^{0}}f\left(s,\,\omega\right)$.
We further assume that the break point is an integer multiple of $m_{T}$,
i.e., $T\lambda_{b}^{0}m_{T}^{-1}\in\left\{ 1,\,2,\ldots,\,\left\lfloor T/m_{T}\right\rfloor -1\right\} $.

In order to simplify the proof, we consider a simplified version of
the problem following \citeReferencesSupp{bibinger/jirak/vetter:16}.
We set $f_{-}\left(\omega_{0}\right)=1$ and let 
\begin{align}
f\left(j/T,\,\omega_{0}\right) & =\begin{cases}
1+\left(m_{T}-j\,\mathrm{mod}\,m_{T}\right)^{\theta}T^{-\theta}, & T\lambda_{b}^{0}<j\leq T\lambda_{b}^{0}+m_{T}\\
1, & \mathrm{else}
\end{cases}.\label{Eq. (86) BJV}
\end{align}
We discuss the general case $f_{-}\left(\omega_{0}\right)\neq1$ at
the end of this proof. Eq. \eqref{Eq. (86) BJV} specifies that the
spectrum at frequency $\omega_{0}$ exhibits a break of order $b_{T}$
at $\lambda_{b}^{0}$ and then decays on the interval $(\lambda_{b}^{0},\,\lambda_{b}^{0}+T^{-1}m_{T}]$
smoothly with regularity $\mathfrak{\theta}$ and is constant elsewhere.
Name this sub-class $\boldsymbol{F}_{\lambda_{b}^{0},\omega_{0}}^{+}.$
Note that here the location of $\lambda_{b}^{0}$ is still unknown.
To establish the lower bound, it suffices to focus on the sub-class
of the above form. 

Next, we introduce a stepwise approximation to $f\left(j/T,\,\omega_{0}\right)$.
Define, for a given sequence $a_{T}$ with $a_{T}\rightarrow\infty$
and $a_{T}m_{T}^{-1}=o(1/\log\left(M_{T}\right)),$ 
\begin{align*}
\widetilde{f}\left(j/T,\,\omega_{0}\right) & =\begin{cases}
1+\left(m_{T}-la_{T}\right)^{\theta}T^{-\theta}, & T\lambda_{b}^{0}+\left(l-1\right)a_{T}<j\leq T\lambda_{b}^{0}+la_{T},\quad1\leq l\leq m_{T}/a_{T}\\
1, & \mathrm{else}
\end{cases}.
\end{align*}
We are given the observations $I_{L,h,T}\left(j/T,\,\omega\right)$
for $j=n_{T}+1,\ldots,\,T$ and $\omega\in\left[-\pi,\,\pi\right]$.
Assume without loss of generality that $\omega_{0}\neq\pm\pi,\,\pm3\pi,\ldots$.
By Theorem \ref{Theorem 5.2.8 Brillinger}-(ii), $I_{L,h,T}\left(j/T,\,\omega_{0}\right)$
is approximately $f\left(j/T,\,\omega_{0}\right)\chi_{2}^{2}/2$ for
$j/T\neq\lambda_{b}^{0}$. For $j/T=\lambda_{b}^{0}$, $I_{L,h,T}\left(j/T,\,\omega_{0}\right)$
is approximately $f\left(j/T,\,\omega_{0}\right)\chi_{2}^{2}/2$,
which also follows from Theorem \ref{Theorem 5.2.8 Brillinger}-(ii)
since by Assumption \ref{Assumption Smothness of A} $f\left(\cdot,\,\omega_{0}\right)$
is continuous from the left at $\lambda_{b}^{0}$. However, note that
$I_{L,h,T}\left(j/T,\,\omega_{0}\right)$ is not asymptotically independent
of $I_{L,h,T}\left(l/T,\,\omega_{0}\right)$ for $l=j-n_{T}+1,\ldots,\,j$.
Let $S_{J}=\left\{ n_{T}+1,\,n_{T}+1+m_{S,T},\ldots\right\} .$ Let
$\zeta_{j}=f\left(j/T,\,\omega_{0}\right)\chi_{2}^{2}/2$ and $\zeta_{j}^{*}=f\left(j/T,\,\omega_{0}\right)\chi_{2}^{2}/2$
where $\zeta_{j}^{*}$ are independent across $j$. Define $\widetilde{\zeta}_{j}^{*}=\widetilde{f}\left(j/T,\,\omega_{0}\right)\chi_{2}^{2}/2$
where $\widetilde{\zeta}_{j}^{*}$ are independent across $j$. 

We distinguish between two cases: (i) $\theta>1/2$ and (ii) $\theta\leq1/2$.

(i) Case $\theta>1/2$. Let us consider the following distinct experiments:

$\mathcal{E}_{1}:$ Observe $\{\zeta_{j}\}_{j=n_{T}+1}^{T}$ and information
$T\lambda_{b}^{0}m_{T}^{-1}\in\{1,\,2,\ldots,\,\left\lfloor T/m_{T}\right\rfloor -1\}$
is provided. 

$\mathcal{E}_{2}:$ Observe $\{\zeta_{j}^{*}\}_{j=n_{T}+1}^{T}$ and
information $T\lambda_{b}^{0}m_{T}^{-1}\in\{1,\,2,\ldots,\,\left\lfloor T/m_{T}\right\rfloor -1\}$
is provided. 

$\mathcal{E}_{3}:$ Observe $\{\widetilde{\zeta}_{j}^{*}\}_{j=n_{T}+1}^{T}$
and information $T\lambda_{b}^{0}m_{T}^{-1}\in\{1,\,2,\ldots,\,\left\lfloor T/m_{T}\right\rfloor -1\}$
is provided. 

$\mathcal{E}_{4}:$ Observe $\chi=((\widetilde{f}(jm_{T}/T,\,\omega_{0})\chi_{2m_{T},j}^{2})_{j\in\mathcal{I}_{1}},\,(\widetilde{f}(\lambda_{b}^{0}+\left(\left(j-1\right)a_{T}+1\right)/T,\,\omega_{0})\widetilde{\chi}_{2m_{T},j}^{2})_{j\in\mathcal{I}_{2}}),$
where $\mathcal{I}_{1}=\{1,\ldots,\,\lambda_{b}^{0}Tm_{T}^{-1},\,\lambda_{b}^{0}Tm_{T}^{-1}+2,\ldots,\,\left\lfloor T/m_{T}\right\rfloor \}$,
$\mathcal{I}_{2}=\{1,\,2,\ldots,\,m_{T}a_{T}^{-1}\}$, and $\{\chi_{2m_{T},j}^{2}\}_{j\in\mathcal{I}_{1}}$
and $\{\widetilde{\chi}_{2a_{T},j}^{2}\}_{j\in\mathcal{I}_{2}}$ are
i.i.d. sequences of chi-square random variables with $2m_{T}$ and
$2a_{T}$ degrees of freedom, respectively. Further, information $T\lambda_{b}^{0}m_{T}^{-1}\in\{1,\,2,\ldots,\,\left\lfloor T/m_{T}\right\rfloor -1\}$
is provided. 

$\mathcal{E}_{5}:$ Observe $\xi=((m_{T}^{1/2}\xi_{j}\widetilde{f}(jm_{T}/T,\,\omega_{0})+\widetilde{f}(jm_{T}/T,\,\omega_{0}))_{j\in\mathcal{I}_{1}},\,(a_{T}^{1/2}\widetilde{\xi}_{j}\widetilde{f}(\lambda_{b}^{0}+((j-1)a_{T}+1)/T,\,\omega_{0})+\widetilde{f}(\lambda_{b}^{0}+((j-1)a_{T}+1)/T,\,\omega_{0}))_{j\in\mathcal{I}_{2}}),$
where $\{\xi_{j}\}_{j\in\mathcal{I}_{1}}$ and $\{\widetilde{\xi}_{j}\}_{j\in\mathcal{I}_{2}}$
are i.i.d. standard normal random variables. Further, information
$T\lambda_{b}^{0}m_{T}^{-1}\in\{1,\,2,\ldots,\,\left\lfloor T/m_{T}\right\rfloor -1\}$
is provided. 

We assume that $\{\zeta_{j}\}$ and $\{\zeta_{j}^{*}\}$ are realized
on the same probability space which is rich enough to allow for both
sequences to be realized there. This is richer then the probability
space in which $\{\zeta_{j}\}$ is realized. Thus, the latter probability
space is extended in the usual way using product spaces. The symbol
$\approx$ denotes asymptotic equivalence while $\sim$ denotes strong
Le Cam equivalence. Our proof consists of showing the following strong
Le Cam equivalence of statistical experiments:
\begin{align}
\mathcal{E}_{1}\approx\mathcal{E}_{2}\thickapprox\mathcal{E}_{3}\sim\mathcal{E}_{4}\thickapprox\mathcal{E}_{5} & .\label{Eq. (87) BJV Strong Equivalence}
\end{align}
 Therefore, given the relation \eqref{Eq. (87) BJV Strong Equivalence},
the lower bound for $\mathcal{E}_{5}$ carries over to the less informative
experiment $\mathcal{E}_{1}$. We prove \eqref{Eq. (87) BJV Strong Equivalence}
in steps. 

Step 1: $\mathcal{E}_{1}\thickapprox\mathcal{E}_{2}$. Given $\zeta_{j}=f\left(j/T,\,\omega_{0}\right)\chi_{2}^{2}/2$
and the boundness of $f\left(\cdot,\,\cdot\right)$, Theorem 1 in
\citeReferencesSupp{berkes/philipp:79} implies that there exists
a sequence $\{\zeta_{j}^{*}\}_{j\in S_{J}}$ of independent random
variables such that $\zeta_{j}^{*}$ has the same distribution as
$\zeta_{j}$ and $\mathbb{P}(|\zeta_{j}-\zeta_{j}^{*}|\geq\nu_{j})\leq\nu_{j}$
with $\nu_{j}>0$. In view of Assumption \ref{Assumption WZ 2011},
we have $\sum_{j=1}^{\infty}\nu_{j}<\infty$, which in turn yields,
\begin{align}
\sum_{j=1}^{\infty}\left|\zeta_{j}-\zeta_{j}^{*}\right|<\infty, & \qquad\mathbb{P}-\mathrm{almost\,surely.}\label{Eq. (after 1.4) BP79}
\end{align}
Note that
\begin{align*}
\left|S_{J}\right|^{-1}\sum_{j\in S_{J}}\left|\zeta_{j}-\zeta_{j}^{*}\right| & =\left|S_{J}\right|^{-1}\sum_{j=n_{T}+1}^{J_{1}}\left|\zeta_{j}-\zeta_{j}^{*}\right|+\left|S_{J}\right|^{-1}\sum_{j\in S_{J},\,j>J_{1}}\left|\zeta_{j}-\zeta_{j}^{*}\right|.
\end{align*}
 Choose $J_{1}$ large enough such that $\sum_{j\in S_{J},\,j>J_{1}}\left|\zeta_{j}-\zeta_{j}^{*}\right|=o_{\mathrm{a.s}}$$\left(\left|S_{J}\right|\right).$
Thus, $\left|S_{J}\right|^{-1}\sum_{j\in S_{J}}\left|\zeta_{j}-\zeta_{j}^{*}\right|\rightarrow0$
$\mathbb{P}$-almost surely. This implies that $||\mathbb{P}_{\{\left|S_{J}\right|^{-1}\zeta_{j}\}}-\mathbb{P}_{\{\left|S_{J}\right|^{-1}\zeta_{j}^{*}\}}||_{\mathrm{TV}}\rightarrow0$.
The latter shows that $\mathcal{E}_{1}\approx\mathcal{E}_{2}$.

Step 2: $\mathcal{E}_{2}\thickapprox\mathcal{E}_{3}$. Note that $c\chi_{2}^{2}$
with $c>0$ is approximately distributed as $\Gamma\left(1,\,2c\right)$
where $\Gamma\left(a,\,b\right)$ is the Gamma distribution with parameters
$\left(a,\,b\right)$. The Kullback-Leibler divergence of $\Gamma\left(1,\,2c\right)$
from $\Gamma\left(1,\,2\widetilde{c}\right)$ is given by 
\begin{align*}
D_{\mathrm{KL}}\left(\mathbb{P}_{c}||\,\mathbb{P}_{\widetilde{c}}\right)= & \left(\log c-\log\widetilde{c}\right)+\frac{\widetilde{c}-c}{c}.
\end{align*}
 For $c=\widetilde{c}+\delta$ with $\delta\rightarrow0,$ we obtain
\begin{align}
D_{\mathrm{KL}}\left(\mathbb{P}_{c}||\,\mathbb{P}_{\widetilde{c}}\right) & =\log\left(\frac{\widetilde{c}+\delta}{\widetilde{c}}\right)+\frac{\widetilde{c}-\left(\widetilde{c}+\delta\right)}{\widetilde{c}+\delta}\nonumber \\
 & =-\frac{\delta^{2}}{2\widetilde{c}^{2}}+O\left(\delta^{2}\right)+O\left(\delta^{3}\right).\label{Eq. (88) BJV 17}
\end{align}
 By Pinsker's inequality,
\begin{align*}
\left\Vert \mathbb{P}_{\left\{ \zeta_{j}^{*}\right\} }-\mathbb{P}_{\left\{ \widetilde{\zeta}_{j}^{*}\right\} }\right\Vert _{\mathrm{TV}}^{2} & \leq\frac{1}{2}D_{\mathrm{KL}}\left(\mathbb{P}_{\zeta_{j}^{*}}||\,\mathbb{P}_{\widetilde{\zeta}_{j}^{*}}\right).
\end{align*}
Thus, using \eqref{Eq. (88) BJV 17} and the additivity of Kullback-Leibler
divergence for independent distributions, we have 
\begin{align*}
D_{\mathrm{KL}}\left(\mathbb{P}_{\zeta_{j}^{*}}||\,\mathbb{P}_{\widetilde{\zeta}_{j}^{*}}\right) & =C\sum_{s=1}^{m_{T}a_{T}^{-1}}\sum_{j=1}^{a_{T}}\left(jT^{-1}\right)^{2\theta}=CO\left(a_{T}T^{-1}\right)^{2\theta}m_{T}.
\end{align*}
This tends to zero in view of \eqref{Eq. (85) BJV 17} and $m_{T}^{-1}a_{T}\rightarrow0.$ 

Step 3: $\mathcal{E}_{3}\sim\mathcal{E}_{4}$. The vector of averages 

\begin{align*}
\left(\left(m_{T}^{-1}\sum_{s=1}^{m_{T}}\widetilde{\zeta}_{jm_{T}+s-1}^{*}\right)_{j\in\mathcal{I}_{1}},\,\left(a_{T}^{-1}\sum_{s=1}^{a_{T}}\widetilde{\zeta}_{T\lambda_{b}^{0}+\left(j-1\right)a_{T}+s-1}^{*}\right)_{j\in\mathcal{I}_{2}}\right) & ,
\end{align*}
forms a sufficient statistic for $\left\{ \widetilde{f}\left(j/T,\,\omega_{0}\right)\right\} _{\left(j/T\right)\in\left[0,\,1\right]}$.
Hence, by Lemma 3.2 of \citeReferencesSupp{brown/low:96} this yields
the strong Le Cam equivalence. 

Step 4: $\mathcal{E}_{4}\thickapprox\mathcal{E}_{5}$. Let 
\begin{align*}
\chi^{*} & =(m_{T}^{-1/2}(\widetilde{f}(jm_{T}/T,\,\omega_{0})(\chi_{2m_{T},j}^{2}-2m_{T}))_{j\in\mathcal{I}_{1}},\\
 & \quad a_{T}^{-1/2}(\widetilde{f}(\lambda_{b}^{0}+((j-1)a_{T}+1)/T,\,\omega_{0})(\widetilde{\chi}_{2m_{T},j}^{2}-2a_{T}))_{j\in\mathcal{I}_{2}})\\
\xi^{*} & =((\xi_{j}\widetilde{f}(jm_{T}/T,\,\omega_{0}))_{j\in\mathcal{I}_{1}},\,(\widetilde{\xi}_{j}\widetilde{f}(\lambda_{b}^{0}+((j-1)a_{T}+1)/T,\,\omega_{0}))_{j\in\mathcal{I}_{2}}).
\end{align*}
Note that $\left\Vert \mathbb{P}_{\chi}-\mathbb{P}_{\xi}\right\Vert _{\mathrm{TV}}^{2}=\left\Vert \mathbb{P}_{\chi^{*}}-\mathbb{P}_{\xi^{*}}\right\Vert _{\mathrm{TV}}^{2}.$
By Pinsker's inequality and independence, 
\begin{align*}
|| & \mathbb{P}_{\chi^{*}}-\mathbb{P}_{\xi^{*}}||_{\mathrm{TV}}^{2}\\
 & \leq2^{-1}D_{\mathrm{KL}}\left(\mathbb{P}_{\chi^{*}}||\,\mathbb{P}_{\xi^{*}}\right)\\
 & \leq2^{-1}\sum_{j\in\mathcal{I}_{1}}D_{\mathrm{KL}}\left(m_{T}^{-1/2}\left(\widetilde{f}\left(jm_{T}/T,\,\omega_{0}\right)\left(\chi_{2m_{T},j}^{2}-2m_{T}\right)\right)||\,\xi_{j}\widetilde{f}\left(jm_{T}/T,\,\omega_{0}\right)\right)\\
 & \quad+2^{-1}\sum_{j\in\mathcal{I}_{2}}D_{\mathrm{KL}}(a_{T}^{-1/2}(\widetilde{f}(\lambda_{b}^{0}+((j-1)a_{T}+1)/T,\,\omega_{0})(\chi_{2a_{T},j}^{2}-2a_{T}))\\
 & \quad||\,\xi_{j}\widetilde{f}(\lambda_{b}^{0}+((j-1)a_{T}+1)/T,\,\omega_{0})).
\end{align*}
 We now apply Theorem 1.1 in \citeReferencesSupp{bobkov/Chi/Gotze:13}
with $c_{1}=12^{-1}\kappa_{3}^{2}$ in their eq. (1.3), where $\kappa_{3}$
is the third-order cumulant of the variable in question. This gives
the following bounds,
\begin{align*}
D_{\mathrm{KL}}\left(\left(m_{T}\right)^{-1/2}\left(\widetilde{f}\left(jm_{T}/T,\,\omega_{0}\right)\left(\chi_{2m_{T},j}^{2}-2m_{T}\right)\right)||\,\xi_{j}\widetilde{f}\left(jm_{T}/T,\,\omega_{0}\right)\right) & =\frac{1}{12}\left(\frac{8}{2m_{T}}\right)+o\left(\frac{1}{m_{T}\log m_{T}}\right),
\end{align*}
 and 
\begin{align*}
D_{\mathrm{KL}} & (a_{T}^{-1/2}(\widetilde{f}(\lambda_{b}^{0}+((j-1)a_{T}+1)/T,\,\omega_{0})(\chi_{2a_{T},j}^{2}-2a_{T}))\\
 & ||\,\xi_{j}\widetilde{f}(\lambda_{b}^{0}+((j-1)a_{T}+1)/T,\,\omega_{0}))=\frac{1}{12}\left(\frac{8}{2a_{T}}\right)+o\left(\frac{1}{a_{T}\log a_{T}}\right).
\end{align*}
Hence, $||\mathbb{P}_{\chi^{*}}-\mathbb{P}_{\xi^{*}}||_{\mathrm{TV}}^{2}=O(Tm_{T}^{-2})+O(m_{T}a_{T}^{-2})$.
Since $\theta>1/2$, we have $Tm_{T}^{-2}\rightarrow0$. Finally,
since $m_{T}^{-1}a_{T}\rightarrow0$, we can choose $a_{T}$ increasing
sufficiently fast such that $m_{T}a_{T}^{-2}\rightarrow0$. Thus,
we have  $||\mathbb{P}_{\chi}-\mathbb{P}_{\xi}||_{\mathrm{TV}}\rightarrow0.$ 

By step 1-4, it is sufficient to establish the minimax lower bound
for experiment $\mathcal{E}_{5}$. After adding an additional drift
$\xi$, which gives an equivalent problem, we cast the problem as
a high dimensional location signal detection problem {[}cf. \citeReferencesSupp{ingster/suslina:03}{]}
from which the bound can be derived using classical arguments. Consider
the observations 
\begin{align*}
\xi^{*} & =((m_{T}^{-1/2}\xi_{j}\widetilde{f}(jm_{T}/T,\,\omega_{0})+\widetilde{f}(jm_{T}/T,\,\omega_{0})-1)_{j\in\mathcal{I}_{1}},\\
 & \quad(a_{T}^{-1/2}\widetilde{\xi}_{j}\widetilde{f}(\lambda_{b}^{0}+((j-1)a_{T}+1)/T,\,\omega_{0})+\widetilde{f}(\lambda_{b}^{0}+((j-1)a_{T}+1)/T,\,\omega_{0})-1)_{j\in\mathcal{I}_{2}}),
\end{align*}
 and the hypothesis 
\begin{align}
\mathcal{H}_{0}:\,\sup_{j}\left(\widetilde{f}\left(j/T,\,\omega_{0}\right)-1\right) & =0\qquad\mathrm{versus}\qquad\mathcal{H}_{1}:\,\sup_{j}\left(\widetilde{f}\left(j/T,\,\omega_{0}\right)-1\right)\geq b_{T}.\label{eq. (91) BJV 17}
\end{align}
The goal is to find the maximal value $b_{T}\rightarrow0$ such that
the hypotheses $\mathcal{H}_{0}$ and $\mathcal{H}_{1}$ are non-distinguishable
in the minimax sense or $\lim_{T\rightarrow\infty}\inf_{\psi}\gamma_{\psi}\left(\theta,\,b_{T}\right)=1.$
Here, the detection rate is $b_{T}\propto\left(T^{-1}m_{T}\right)^{\theta}\propto T^{-\frac{\theta}{2\theta+1}}.$
Consider the product measures $\mathbb{P}_{\mathcal{H}_{0}}=\mathbb{P}_{\xi^{*}}\times\mathbb{P}_{0}$
and $\mathbb{P}_{\mathcal{H}_{1}}=\mathbb{P}_{\xi^{*}}\times\mathbb{P}_{\lambda_{b}^{0},1}$
where $\mathbb{P}_{\xi^{*}}$ is the probability law of $\xi^{*}$
and $\mathbb{P}_{0}$ is the measure for the no break case. Thus,
$\mathbb{P}_{\mathcal{H}_{0}}$ is the probability measure under $\mathcal{H}_{0}$
while $\mathbb{P}_{\mathcal{H}_{1}}$ is the probability measure under
$\mathcal{H}_{1}$ which draws a break at time $\lambda_{b}^{0}$
with $T\lambda_{b}^{0}m_{T}^{-1}\in\left\{ 1,\,2,\ldots,\,\left\lfloor T/m_{T}\right\rfloor -1\right\} $
uniformly from this set. From similar derivations that yield eq. (2.20)-(2.22)
in \citeReferencesSupp{ingster/suslina:03}, it follows that 
\begin{align*}
\inf_{\psi}\gamma_{\psi}\left(\theta,\,b_{T}\right) & \geq1-\frac{1}{2}\left\Vert \mathbb{P}_{\mathcal{H}_{1}}-\mathbb{P}_{\mathcal{H}_{0}}\right\Vert _{\mathrm{TV}}\geq1-\frac{1}{2}\left|\mathbb{E}_{\mathbb{P}_{\mathcal{H}_{0}}}\left(\mathscr{L}_{0,1}^{2}-1\right)\right|^{1/2},
\end{align*}
 where $\mathscr{L}_{0,1}=d\mathbb{P}_{\mathcal{H}_{1}}/d\mathbb{P}_{\mathcal{H}_{0}}$
is the likelihood ratio between $\mathbb{P}_{\mathcal{H}_{1}}$ and
$\mathbb{P}_{\mathcal{H}_{0}}.$ By the above inequality, it is sufficient
to show $\mathbb{E}_{\mathbb{P}_{\mathcal{H}_{0}}}(\mathscr{L}_{0,1}^{2})\rightarrow1$.
The proof of the latter result follows similar arguments as in \citeReferencesSupp{bibinger/jirak/vetter:16}.

It remains to consider the case $\theta\leq1/2$. In a different setting,
\citeReferencesSupp{bibinger/jirak/vetter:16} considered separately
the case where their regularity exponent $\mathfrak{a}$ satisfies
$\mathfrak{a}\leq1/2$ to obtain the minimax lower bound. The same
arguments can be applied in our context which lead to the same result
as for the case $\theta>1/2$. 

The general case with $f_{-}\left(\omega_{0}\right)>0$ rather than
with $f_{-}\left(\omega_{0}\right)=1$ as discussed above follows
from the same arguments after we rescale the equations in \eqref{eq. (91) BJV 17}.
The only difference is the form of the detection rate which is now
$b_{T}\leq f_{-}(\omega_{0})D(T^{-1}m_{T})^{\theta}$.

The proof of the lower bound for the alternative $\mathcal{H}_{1}^{\mathrm{S}}$
is similar to the proof discussed above. The minor differences in
the proof outlined by \citeReferencesSupp{bibinger/jirak/vetter:16}
also apply here. $\square$

\subsubsection{Proof of Theorem \ref{Theorem 4.2 BJV 17}}

We present the proof for the statistic $\mathrm{S}_{\mathrm{max},T}$.
The proof for the other test statistics discussed in Section \ref{Section Change-Point-Tests}
is similar and omitted. From the same reasoning as in the proofs of
the results of Section \ref{Section Change-Point-Tests}, we can replace
$\widehat{\sigma}_{L,r}\left(\omega\right)$ by $\sigma_{L,r}\left(\omega\right)$
throughout the proof. Without loss of generality, we assume that $\omega_{0}\neq\pm\pi.$
Let $M_{S,T}^{*}=m_{T}^{*}/m_{S,T}^{*}$ and $m_{S,T}/m_{T}^{*}\rightarrow[0,\,\infty)$.
If $\left\lfloor T\lambda_{b}^{0}\right\rfloor \notin\{\{\mathbf{S}_{r}\}\cup\{\mathbf{S}_{r+1}\}\}$
or if $\omega\neq\omega_{0}$, then 
\begin{align*}
\biggl| & \frac{\widetilde{f}_{L,r,T}\left(\omega\right)-\widetilde{f}_{R,r+1,T}\left(\omega\right)}{\sigma_{L,r}\left(\omega\right)}\biggr|\\
 & =\left|\frac{\left(M_{S,T}^{*}\right)^{-1}\sum_{j\in\mathbf{S}_{r}}\left(f_{L,h,T}^{*}\left(j/T,\,\omega\right)+\mathbb{E}\left(f_{L,h,T}\left(j/T,\,\omega\right)\right)\right)}{\sigma_{L,r}\left(\omega\right)}\right.\\
 & -\left.\frac{\left(M_{S,T}^{*}\right)^{-1}\sum_{j\in\mathbf{S}_{r+1}}\left(f_{R,h,T}^{*}\left(j/T,\,\omega\right)+\mathbb{E}\left(f_{R,h,T}\left(j/T,\,\omega\right)\right)\right)}{\sigma_{L,r}\left(\omega\right)}\right|\\
 & =\left|\frac{\left(M_{S,T}^{*}\right)^{-1}\sum_{j\in\mathbf{S}_{r}}f_{L,h,T}^{*}\left(j/T,\,\omega\right)-\left(M_{S,T}^{*}\right)^{-1}\sum_{j\in\mathbf{S}_{r+1}}f_{R,h,T}^{*}\left(j/T,\,\omega\right)}{\sigma_{L,r}\left(\omega\right)}\right|\\
 & \quad+O\left(\left(m_{T}^{*}/T\right)^{\theta}\right)+O_{\mathbb{P}}\left(\left(n_{T}/T\right)^{2}+\log\left(n_{T}\right)/n_{T}\right)+o\left(b_{W,T}^{2}\right)\\
 & \triangleq\mathring{f}_{r,T}\left(\omega\right)+O\left(\left(m_{T}^{*}/T\right)^{\theta}\right)+O_{\mathbb{P}}\left(\left(n_{T}/T\right)^{2}+\log\left(n_{T}\right)/n_{T}\right)+o\left(b_{W,T}^{2}\right)\\
 & =\mathring{f}_{r,T}\left(\omega\right)+o_{\mathbb{P}}\left(\left(\sqrt{m_{T}^{*}}\right)^{-1}\right),
\end{align*}
where the last inequality follows from \eqref{Eq. Cond for Sec. 5}.
As in the proof of Theorem \ref{Theorem Asymptotic H0 Distrbution Smax},
we have $\sqrt{M_{S,T}^{*}}\mathring{f}_{r,T}\left(\omega\right)=O_{\mathbb{P}}\left(1\right)$
for $1\leq r\leq M_{T}^{*}-2$. This can be used to obtain the following
inequality, if $\left\lfloor T\lambda_{b}^{0}\right\rfloor \in\{\{\mathbf{S}_{r}\}\cup\{\mathbf{S}_{r+1}\}\}$
and $\omega=\omega_{0}$, 
\begin{align}
\mathrm{S}_{\mathrm{max},T} & \left(\omega_{0}\right)\geq-\mathring{f}_{r,T}\left(\omega_{0}\right)\nonumber \\
 & \quad+\frac{T}{m_{T}^{*}}\left|\int_{\left(rm_{T}^{*}-m_{T}^{*}/2+n_{T}/2+1\right)/T}^{\lambda_{b}^{0}}f\left(u,\,\omega_{0}\right)du-\int_{\lambda_{b}^{0}}^{(\left(r+1\right)m_{T}^{*}+n_{T}/2+m_{S,T}M_{S,T}/2)/T}f\left(u,\,\omega_{0}\right)du\right|\nonumber \\
 & \quad\times\frac{\left(1-o_{\mathbb{P}}\left(1\right)\right)}{\sup_{u}f(u,\,\omega_{0})}\nonumber \\
 & \geq-O_{\mathbb{P}}\left(\left(m_{T}^{*}\right)^{-1/2}\right)\label{Eq. (97) BJV 17}\\
 & \quad+\frac{T}{m_{T}^{*}}\left|\int_{\left(rm_{T}^{*}-m_{T}^{*}/2+n_{T}/2+1\right)/T}^{\lambda_{b}^{0}}f\left(u,\,\omega_{0}\right)du-\int_{\lambda_{b}^{0}}^{(\left(r+1\right)m_{T}^{*}+n_{T}/2+m_{S,T}M_{S,T}/2)/T}f\left(u,\,\omega_{0}\right)du\right|\nonumber \\
 & \quad\times\frac{\left(1-o_{\mathbb{P}}\left(1\right)\right)}{\sup_{u}f(u,\,\omega_{0})}.\nonumber 
\end{align}
Note that $\gamma_{\psi^{*}}\left(\theta,\,b_{T}^{*}\right)\rightarrow0$
follows from 
\begin{align}
\mathbb{P}\left(\mathrm{S}_{\mathrm{max},T}\left(\omega\right)<2D^{*}\sqrt{\log\left(M_{T}^{*}\right)/m_{T}^{*}}\right)\rightarrow1, & \quad\mathrm{for}\,\mathrm{all}\,\omega\in\left[-\pi,\,\pi\right],\qquad\mathrm{under}\,\mathcal{H}_{0}\label{eq. (99) BJV 17}\\
\mathbb{P}\left(\mathrm{S}_{\mathrm{max},T}\left(\omega\right)\geq2D^{*}\sqrt{\log\left(M_{T}^{*}\right)/m_{T}^{*}}\right)\rightarrow1, & \quad\mathrm{for}\,\mathrm{some}\,\omega\in\left[-\pi,\,\pi\right],\quad\mathrm{under}\,\mathcal{H}_{1}^{\mathrm{B}}\,\mathrm{or}\,\mathcal{H}_{1}^{\mathrm{S}}.\label{E. (98) BJV 17}
\end{align}
We first show \eqref{eq. (99) BJV 17}. Note that 
\begin{align*}
2D^{*}\sqrt{\log\left(M_{T}^{*}\right)/m_{T}^{*}} & \geq2\sqrt{\log\left(M_{T}^{*}\right)/m_{T}^{*}}+D\left(m_{T}^{*}/T\right)^{\theta}.
\end{align*}
Under $\mathcal{H}_{0}$, since $\theta'<\theta$ we have for all
$\omega\in\left[-\pi,\,\pi\right]$,
\begin{align*}
\mathrm{S}_{\mathrm{max},T}\left(\omega\right) & \leq\max_{1\leq r\leq M_{T}^{*}-2}\mathring{f}_{r,T}\left(\omega\right)+D\left(m_{T}^{*}/T\right)^{\theta'}+O_{\mathbb{P}}\left(\left(n_{T}/T\right)^{2}+\log\left(n_{T}\right)/n_{T}+o\left(b_{W,T}^{2}\right)\right).
\end{align*}
 Given \eqref{Eq. Cond for Sec. 5}, to conclude the proof, we have
to show 
\[
\mathbb{P}\left(\max_{1\leq r\leq M_{T}^{*}-2}\mathring{f}_{r,T}\left(\omega_{0}\right)\leq\sqrt{\log\left(M_{T}^{*}\right)/m_{T}^{*}}\right)\rightarrow1.
\]
 The latter result follows from $\sqrt{\log\left(M_{T}^{*}\right)/m_{T}^{*}}\leq\sqrt{\log\left(M_{T}^{*}\right)/M_{S,T}^{*}}$
which is implied by Theorem \ref{Theorem Asymptotic H0 Distrbution Smax}.

We now prove \eqref{E. (98) BJV 17} under $\mathcal{H}_{1}^{\mathrm{B}}$.
We have to show that the second term on the right hand side of \eqref{E. (98) BJV 17}
is greater than or equal to $2D^{*}\sqrt{\log\left(M_{T}^{*}\right)/m_{T}^{*}}.$
The term in question is larger than $b_{T}^{*}-2D\left(m_{T}^{*}/T\right)^{\theta}.$
In view of \eqref{Eq. (34) BJV 17} with $\theta'=0$ the result
follows. 

We now prove \eqref{E. (98) BJV 17} under $\mathcal{H}_{1}^{\mathrm{S}}$.
For $h\leq2m_{T}^{*}/T$, we have $f\left(\lambda_{b}^{0}+h,\,\omega_{0}\right)\geq f\left(\lambda_{b}^{0},\,\omega_{0}\right)+b_{T}^{*}h^{\theta'}$
or $f\left(\lambda_{b}^{0}+h,\,\omega_{0}\right)\leq f\left(\lambda_{b}^{0},\,\omega_{0}\right)-b_{T}^{*}h^{\theta'}$.
Thus, 
\begin{align*}
\frac{T}{m_{T}^{*}}\left|\int_{\lambda_{b}^{0}+m_{T}^{*}/T}^{\lambda_{b}^{0}+2m_{T}^{*}/T}\left(f\left(u,\,\omega_{0}\right)-f\left(u-m_{T}^{*}/T,\,\omega_{0}\right)\right)du\right| & \geq b_{T}\left(m_{T}^{*}/T\right)^{\theta'}\\
 & \geq2D^{*}\sqrt{\log\left(M_{T}^{*}\right)/m_{T}^{*}},
\end{align*}
where the second equality follows from \eqref{Eq. (34) BJV 17}.
$\square$ 

\subsection{\label{Subsec: Proofs Change-point estimator}Proofs of the Results
of Section \ref{Section Estimation-of-the Breaks}}

From the same reasoning as in the proofs of the results of Section
\ref{Section Change-Point-Tests}, we can replace $\widehat{\sigma}_{L,r}\left(\omega\right)$
by $\sigma_{L,r}\left(\omega\right)$ throughout the proofs of this
section.

\subsubsection{Proof of Proposition \ref{Prop 4.5 BJV Consitency Change-point}}

The following lemma is simple to verify.
\begin{lem}
\label{Lemma D1 BJV}Let $C\left(u\right)$ and $d\left(u\right)$
be functions on $\left[0,\,\lambda_{b}^{0}\right]$ such that $d\left(u\right)$
is increasing. As long as $d\left(\lambda_{b}^{0}\right)-d\left(\lambda_{b}^{0}-\kappa\right)\geq\sup_{0\leq u\leq\lambda_{b}^{0}}\left|C\left(u\right)\right|$
for some $\kappa\in\left[0,\,\lambda_{b}^{0}\right]$ we have that,
\begin{align}
\underset{0\leq u\leq\lambda_{b}^{0}}{\mathrm{argmax}}\left(d\left(u\right)+C\left(u\right)\right) & \geq\lambda_{b}^{0}-\kappa.\label{Eq. (BJV) Lemma}
\end{align}
An analogous results holds if $C\left(u\right)$ and $d\left(u\right)$
are functions on $\left[\lambda_{b}^{0},\,1\right]$ and $d\left(u\right)$
is decreasing.
\end{lem}
\noindent\textit{Proof of Proposition} \ref{Prop 4.5 BJV Consitency Change-point}.
For $\lambda_{b}^{0}\in\left(0,\,1\right)$ define $\overline{r}_{b}=\left\lceil T\lambda_{b}^{0}+1\right\rceil $,
i.e., the smallest integer such that $\overline{r}_{b}/T$ is larger
than or equal to $\lambda_{b}^{0}+1/T$. Denote by $\{\widetilde{f}\left(u,\,\omega_{0}\right)\}_{u\in\left[0,\,1\right]}$
the path of the spectrum $f\left(\cdot,\,\omega_{0}\right)$ without
the break: $f\left(r/T,\,\omega\right)=\widetilde{f}\left(r/T,\,\omega\right)+\delta_{T}\mathbf{1}\left\{ r\geq\overline{r}_{b}\right\} .$
Without loss of generality, we assume $\delta_{T}>0$. Define $d\left(r/T,\,\omega\right)=0$
for $\omega\neq\omega_{0}$ and
\begin{align*}
d\left(r/T,\,\omega_{0}\right) & =\begin{cases}
0 & \textrm{if }r+m_{T}<\overline{r}_{b},\\
\left(r+m_{T}-\overline{r}_{b}\right)m_{S,T}^{-1}M_{S,T}^{-1/2}\delta_{T} & \textrm{if }r=\overline{r}_{b}-m_{T},\,\overline{r}_{b}-m_{T}+m_{S,T}\ldots,\\
M_{S,T}^{1/2}\delta_{T} & \textrm{if }r>\overline{r}_{b},
\end{cases}\overline{r}_{b},
\end{align*}
 and $\left\{ d\left(u,\,\omega_{0}\right)\right\} _{u\in\left[0,\,1\right]}$
is the associated piecewise constant increasing step function. By
Lemma \ref{Lemma Smax - Smax_old 0 } it is sufficient to consider
\begin{align}
\mathrm{D}'_{r,T}\left(\omega\right)=M_{S,T}^{-1/2}\left|\sum_{j\in\mathbf{S}_{L,r}}f_{h,T}\left(j/T,\,\omega\right)-\sum_{j\in\mathbf{S}_{R,r}}f_{h,T}\left(j/T,\,\omega\right)\right|, & \qquad\omega\in\left[-\pi,\,\pi\right].\label{Eq. D'}
\end{align}
For $r=m_{T},\,2m_{T},\ldots,$ write 
\begin{align*}
\sum_{j\in\mathbf{S}_{L,r}} & f_{h,T}\left(j/T,\,\omega_{0}\right)-\sum_{j\in\mathbf{S}_{R,r}}f_{h,T}\left(j/T,\,\omega_{0}\right)\\
 & =\sum_{j\in\mathbf{S}_{L,r}}\left(f_{h,T}\left(j/T,\,\omega_{0}\right)-\mathbb{E}\left(f_{h,T}\left(j/T,\,\omega_{0}\right)\right)\right)-\sum_{j\in\mathbf{S}_{R,r}}\left(f_{h,T}\left(j/T,\,\omega_{0}\right)-\mathbb{E}\left(f_{h,T}\left(j/T,\,\omega_{0}\right)\right)\right)\\
 & \quad+\sum_{j\in\mathbf{S}_{L,r}}\left(\mathbb{E}\left(f_{h,T}\left(j/T,\,\omega_{0}\right)\right)-\widetilde{f}\left(j/T,\,\omega_{0}\right)\right)-\sum_{j\in\mathbf{S}_{R,r}}\left(\mathbb{E}\left(f_{h,T}\left(j/T,\,\omega_{0}\right)\right)-f\left(j/T,\,\omega_{0}\right)\right)\\
 & \quad+\sum_{j\in\mathbf{S}_{L,r}}\widetilde{f}\left(j/T,\,\omega_{0}\right)-\sum_{j\in\mathbf{S}_{R,r}}\widetilde{f}\left(j/T,\,\omega_{0}\right)-\sum_{j\in\mathbf{S}_{R,r}}\left(f\left(j/T,\,\omega_{0}\right)-\widetilde{f}\left(j/T,\,\omega_{0}\right)\right).
\end{align*}
 For $r=2m_{T},\ldots,\,\overline{r}_{b}$ let $C\left(r/T,\,\omega\right)=\mathrm{D}'_{r,T}\left(\omega\right)$
for $\omega\neq\omega_{0}$ and 
\begin{align*}
C\left(r/T,\,\omega_{0}\right) & =M_{S,T}^{-1/2}\left(\sum_{j\in\mathbf{S}_{L,r}}f_{h,T}\left(j/T,\,\omega_{0}\right)-\sum_{j\in\mathbf{S}_{R,r}}f_{h,T}\left(j/T,\,\omega_{0}\right)\right.\\
 & \quad\left.+\sum_{j\in\mathbf{S}_{R,r},\,j>\overline{r}_{b}}\left(f\left(j/T,\,\omega_{0}\right)-\widetilde{f}\left(j/T,\,\omega_{0}\right)\right)\right),
\end{align*}
for $\omega=\omega_{0}$. Note that $C\left(s/T,\,\omega\right)$
does not involve any break for any $\omega.$  Thus, we can proceed
similarly as in the proofs of Section \ref{Section Change-Point-Tests}.
That is, we exploit the smoothness of $f\left(\cdot,\,\cdot\right)$
under $\mathcal{H}_{0}$ to yield $\sup_{u\in\left[0,\,\lambda_{b}^{0}\right]}\sup_{\omega\in\left[-\pi,\,\pi\right]}\left|C\left(u,\,\omega\right)\right|=O_{\mathbb{P}}(\sqrt{\log\left(T\right)})$.
This combined with the definition of $d\left(r/T,\,\omega_{0}\right)$
implies that for each $r=\overline{r}_{b}-\left\lfloor m_{T}/B\right\rfloor ,\ldots,\,\overline{r}_{b}$,
where $B$ is any finite integer with $B>1,$
\[
\left|d\left(r/T,\,\omega_{0}\right)\right|>\max_{\omega\in\left[-\pi,\,\pi\right]}\left(\left|C\left(r/T,\,\omega\right)\right|\right)>0,
\]
with probability approaching one and  
\begin{align*}
\mathrm{D}_{r,T}\left(\omega\right)=\left|d\left(r/T,\,\omega\right)+C\left(r/T,\,\omega\right)\right| & =d\left(r/T,\,\omega\right)+\mathrm{sign}\left(C\left(r/T,\,\omega\right)\right)\left|C\left(r/T,\,\omega\right)\right|.
\end{align*}
By the definition of $d\left(\cdot,\,\omega_{0}\right)$, for $\kappa_{T}\in[0,\,m_{T}/\left(BT\right)]$,
\begin{align*}
d\left(\overline{r}_{b}/T,\,\omega_{0}\right)-d\left(\overline{r}_{b}/T-\kappa_{T},\,\omega_{0}\right) & =\left\lfloor \kappa_{T}T\right\rfloor m_{S,T}^{-1}\delta_{T}M_{S,T}^{-1/2}.
\end{align*}
 In order to apply Lemma \ref{Lemma D1 BJV}, we need to choose $\kappa_{T}$
such that $\left\lfloor \kappa_{T}T\right\rfloor m_{S,T}^{-1}\delta_{T}M_{S,T}^{-1/2}/\sqrt{\log\left(T\right)}\geq1$
or $\sqrt{M_{S,T}\log\left(T\right)}m_{S,T}/\left(\delta_{T}T\right)=o\left(\kappa_{T}\right).$
Lemma \ref{Lemma D1 BJV} then yields 
\begin{align*}
\frac{\overline{r}_{b}}{T}\geq & \underset{r=2m_{T},\,3m_{T}\ldots;\,r<\overline{r}_{b}}{\mathrm{argmax}}\max_{\omega\in\left[-\pi,\,\pi\right]}T^{-1}\mathrm{D}_{r,T}\left(\omega\right)=\underset{r=m_{T},\ldots,\,\overline{r}_{b}}{\mathrm{argmax}}T^{-1}\mathrm{D}_{r,T}\left(\omega_{0}\right)\geq\frac{\overline{r}_{b}}{T}-\kappa_{T}.
\end{align*}
The case $r>\overline{r}_{b}$ can be treated similarly by symmetry.
It results in 
\begin{align*}
\frac{\overline{r}_{b}}{T}\leq & \underset{r=\overline{r}_{b},\ldots,\,T-m_{T}}{\mathrm{argmax}}\max_{\omega\in\left[-\pi,\,\pi\right]}T^{-1}\mathrm{D}_{r,T}\left(\omega\right)=\underset{r=\overline{r}_{b},\ldots,\,T-m_{T}}{\mathrm{argmax}}T^{-1}\mathrm{D}_{r,T}\left(\omega_{0}\right)\leq\frac{\overline{r}_{b}}{T}+\kappa_{T}.
\end{align*}
 Therefore, we conclude that $|\widehat{\lambda}_{b}-\overline{r}_{b}/T|=O_{\mathbb{P}}\left(\kappa_{T}\right)\rightarrow0$.
$\square$

\subsubsection{Proof of Proposition \ref{Prop 4.9 BJV Consitency Change-point}}

Set $\widehat{\mathcal{I}}=\left\{ 2m_{T},\,3m_{T},\ldots,\,\left(M_{T}-1\right)m_{T}-n_{T}\right\} \backslash\left\{ 2m_{T}\right\} $
and $\widehat{\mathcal{T}}=\emptyset$. Under $\mathcal{H}_{1,\mathrm{M}}$,
the arguments in the proof of Theorem \ref{Theorem Asymptotic H0 Distrbution SDmax}
yields,
\begin{align*}
\underset{r\in\widehat{\mathcal{I}}}{\mathrm{max}}\max_{\omega\neq\omega_{1},\ldots,\,\omega_{m_{0}}}\mathrm{D}_{r,T}\left(\omega\right)=O_{\mathbb{P}}\left(\sqrt{\log\left(T\right)}\right) & .
\end{align*}
Let $r_{L,l},\,r_{R,l}\in\widehat{\mathcal{I}}$ $(l=1,\ldots,\,m_{0})$
such that $r_{R,l}=r_{L,l}+m_{T}$ and $r_{L,l}\leq T_{l}^{0}<r_{R,l}$.
For any $\omega$, we have 
\begin{align*}
\underset{r\in\widehat{\mathcal{I}}\backslash\left\{ r_{L,1},\,r_{R,1},\ldots,\,r_{L,m_{0}},\,r_{R,m_{0}}\right\} }{\mathrm{max}}\mathrm{D}_{r,T}\left(\omega\right)=O_{\mathbb{P}}\left(\sqrt{\log\left(T\right)}\right) & .
\end{align*}
For each $r\in\widehat{\mathcal{I}}$, we draw $K$ points $r_{k,r}^{\diamond}$
with $k=1,\ldots,\,K$ uniformly (without replacement) from $\mathbf{I}\left(r\right)$.
Consider the following events,
\begin{align*}
\mathbf{D}_{1}= & \left\{ \forall r\in\widehat{\mathcal{I}}\,\mathrm{and}\,\forall\,k=1,\ldots,\,K,\,\left(\exists!\,1\leq l\leq m_{0}\right)\vee\left(\nexists\,1\leq l\leq m_{0}\right)\,\mathrm{s.t.}\,T_{l}^{0}\in\left[r_{k,r}^{\diamond}-m_{T},\,r_{k,r}^{\diamond}+m_{T}\right]\right\} \\
\mathbf{D}_{2}= & \left\{ \forall l=1,\ldots,\,m_{0}\,\exists\,r\in\widehat{\mathcal{I}}\,\mathrm{s.t.}\,\exists\,k=1,\ldots,\,K,\,\mathrm{s.t.}\,\left|T_{l}^{0}-r_{k,r}^{\diamond}\right|=Cm_{T}\,\mathrm{for}\,\mathrm{some}\,C\in[0,\,1)\right\} .
\end{align*}
Let $\mathbf{A}^{c}$ denote the complement of a set $\mathbf{A}.$
Note that $\mathbb{P}\left(\left(\mathbf{D}_{1}\cap\mathbf{D}_{2}\right)^{c}\right)=\mathbb{P}\left(\left(\mathbf{D}_{2}\right)^{c}\right)$
by Assumption \ref{Assumption small shift Multiple breaks} and that
$\mathbb{P}\left(\left(\mathbf{D}_{2}\right)^{c}\right)=0$ if there
are still undetected breaks. 

The remaining arguments will be valid on the set $\mathbf{D}_{1}\cap\mathbf{D}_{2}$
as long as there are undetected breaks. Let $r_{l},\,r_{l+1}\in\widehat{\mathcal{I}}$
be such that $T_{l}^{0}\in[r_{l},\,r_{l+1})$. As in the proof of
Proposition \ref{Prop 4.5 BJV Consitency Change-point}, 
\begin{align*}
\mathrm{D}_{r_{l},T}\left(\omega_{l}\right) & =|O_{\mathbb{P}}\left(M_{S,T}^{-1/2}\delta_{l,T}\left(M_{S,T}-\left(r_{l}-T_{l}^{0}\right)\mathbf{1}\left\{ r_{l-1}<T_{l}^{0}\leq r_{l}\right\} +\left(r_{l+1}-T_{l}^{0}\right)\mathbf{1}\left\{ r_{l}<T_{l}^{0}<r_{l+1}\right\} \right)\right)|.
\end{align*}
Note that if $\mathrm{D}_{r_{l},T}\left(\omega_{l}\right)/\left(\delta_{l,T}\sqrt{M_{S,T}}\right)\overset{\mathbb{P}}{\rightarrow}0$
then we must have $\mathrm{D}_{r_{l+1},T}\left(\omega_{l}\right)=O_{\mathbb{P}}\left(\delta_{l,T}\sqrt{M_{S,T}}\right)$.
Using a similar argument as in Lemma \ref{Lemma Smax - Smax_old 0 }
one can show that $\mathrm{S}_{\mathrm{Dmax},T}\left(\widehat{\mathcal{I}}\right)$
is asymptotically equivalent to $\underset{r\in\widehat{\mathcal{I}}}{\mathrm{max}}\,\underset{k\in K}{\mathrm{max}}\underset{\omega\in\left[-\pi,\,\pi\right]}{\mathrm{max}}\mathrm{D}_{r_{k,r}^{\diamond},T}\left(\omega\right).$
Thus, in step (2) $\psi(\left\{ X_{t,T}\right\} ,\,\widehat{\mathcal{I}})=1$
because for large enough $T$,
\begin{align*}
\underset{r\in\widehat{\mathcal{I}}}{\mathrm{max}}\,\underset{k\in K}{\mathrm{max}}\underset{\omega\in\left[-\pi,\,\pi\right]}{\mathrm{max}}\mathrm{D}_{r_{k,r}^{\diamond},T}\left(\omega_{0}\right) & \geq\underset{r\in\widehat{\mathcal{I}}}{\mathrm{max}}\underset{\omega\in\left[-\pi,\,\pi\right]}{\mathrm{max}}\mathrm{D}_{r,T}\left(\omega\right)\\
 & =|\delta_{l,T}O_{\mathbb{P}}\left(\sqrt{M_{S,T}}\right)|\\
 & \geq\inf_{1\leq l\leq m_{0}}|\delta_{l,T}O_{\mathbb{P}}\left(\sqrt{M_{S,T}}\right)|\\
 & =2D^{*}\left(\log\left(T\right)\right)^{2/3}\\
 & >2D^{*}\sqrt{\log\left(M_{T}^{*}\right)},
\end{align*}
 where the last equality follows from Assumption \ref{Assumption small shift Multiple breaks}.
We now move to step (3). By the arguments in the proof of Proposition
\ref{Prop 4.5 BJV Consitency Change-point}, there exists $1\leq l\leq m_{0}$
such that $|\lambda_{l}^{0}-\widehat{\lambda}_{T}(\widehat{\mathcal{I}})|\leq m_{T}/T$.
Since $\inf_{1\leq l\leq m_{0}-1}|\lambda_{l+1}^{0}-\lambda_{l}^{0}|\geq\nu_{T}^{-1}$
and $m_{T}/v_{T}\rightarrow0$ there can exist exactly one $l$ that
satisfies $|\lambda_{l}^{0}-\widehat{\lambda}_{T}(\widehat{\mathcal{I}})|\leq m_{T}/T$.
For such a $\lambda_{l}^{0}$ define $\overline{r}_{l,b}=\left\lceil T\lambda_{l}^{0}+1\right\rceil $,
the smallest integer such that $\overline{r}_{l,b}/T$ is larger than
or equal to $\lambda_{l}^{0}+1/T$. Denote by $\{\widetilde{f}\left(u,\,\omega\right)\}_{u\in\left[0,\,1\right]}$
the path of the spectrum $f\left(\cdot,\,\omega\right)$ without the
break $\delta_{l,T}$: 
\begin{align*}
f\left(r/T,\,\omega_{l}\right) & =\widetilde{f}\left(r/T,\,\omega_{l}\right)+\delta_{l,T}\mathbf{1}\left\{ r\geq\overline{r}_{l,b}\right\} .
\end{align*}
Without loss of generality, we assume $\delta_{l,T}>0$. Define $d_{l}\left(r/T,\,\omega\right)=0$
for $\omega\neq\omega_{l}$ and
\begin{align*}
d_{l}\left(r/T,\,\omega_{l}\right) & =\begin{cases}
0 & \textrm{if }r+m_{T}<\overline{r}_{l,b},\\
\left(r+m_{T}-\overline{r}_{l,b}\right)m_{S,T}^{-1}M_{S,T}^{-1/2}\delta_{l,T} & \textrm{if }r=\overline{r}_{l,b}-m_{T},\,\overline{r}_{l,b}-m_{T}+m_{S,T},\ldots,\\
M_{S,T}^{1/2}\delta_{l,T} & \textrm{if }r>\overline{r}_{l,b},
\end{cases}\overline{r}_{l,b},
\end{align*}
for $\omega=\omega_{l}$. Let $\left\{ d\left(u\right)\right\} _{u\in\left[0,\,1\right]}$
be the associated piecewise constant increasing step function. For
any $r\in\widehat{\mathcal{I}}$, write 
\begin{align*}
\sum_{j\in\mathbf{S}_{L,r}} & f_{h,T}\left(j/T,\,\omega\right)-\sum_{j\in\mathbf{S}_{R,r}}f_{h,T}\left(j/T,\,\omega\right)\\
 & =\sum_{j\in\mathbf{S}_{L,r}}\left(f_{h,T}\left(j/T,\,\omega\right)-\mathbb{E}\left(f_{h,T}\left(j/T,\,\omega\right)\right)\right)-\sum_{j\in\mathbf{S}_{R,r}}\left(f_{h,T}\left(j/T,\,\omega\right)-\mathbb{E}\left(f_{h,T}\left(j/T,\,\omega\right)\right)\right)\\
 & \quad+\sum_{j\in\mathbf{S}_{L,r}}\left(\mathbb{E}\left(f_{h,T}\left(j/T,\,\omega\right)\right)-\widetilde{f}\left(j/T,\,\omega\right)\right)-\sum_{j\in\mathbf{S}_{R,r}}\left(\mathbb{E}\left(f_{h,T}\left(j/T,\,\omega\right)\right)-f\left(j/T,\,\omega\right)\right)\\
 & \quad+\sum_{j\in\mathbf{S}_{L,r}}\widetilde{f}\left(j/T,\,\omega\right)-\sum_{j\in\mathbf{S}_{R,r}}\widetilde{f}\left(j/T,\,\omega\right)-\sum_{j\in\mathbf{S}_{R,r}}\left(f\left(j/T,\,\omega\right)-\widetilde{f}\left(j/T,\,\omega\right)\right).
\end{align*}
 For $r=2m_{T},\,3m_{T}\ldots,\,\overline{r}_{b},$ let $C_{l}\left(r/T,\,\omega\right)=\mathrm{D'}_{r,T}\left(\omega\right)$
for $\omega\neq\omega_{l}$, where $\mathrm{D'}_{r,T}\left(\omega\right)$
is given in \eqref{Eq. D'} and 

\begin{align*}
C_{l}\left(r/T,\,\omega_{l}\right) & =M_{S,T}^{-1/2}\left(\sum_{j\in\mathbf{S}_{L,r}}f_{h,T}\left(j/T,\,\omega_{l}\right)-\sum_{j\in\mathbf{S}_{R,r}}f_{h,T}\left(j/T,\,\omega_{l}\right)\right.\\
 & \quad\left.+\sum_{j\in\mathbf{S}_{R,r},\,j>\overline{r}_{l,b}}\left(f\left(j/T,\,\omega_{l}\right)-\widetilde{f}\left(j/T,\,\omega_{l}\right)\right)\right),
\end{align*}
for $\omega=\omega_{l}$.   We proceed as in the proof of Proposition
\ref{Prop 4.5 BJV Consitency Change-point}. We have 
\begin{align*}
d\left(r/T,\,\omega_{l}\right)\geq\max_{\omega\in\{\left[-\pi,\,\pi\right]/\{\omega_{1},\ldots,\,\omega_{m_{0}}\}\}}\left|d\left(r/T,\,\omega\right)\right|>0 & ,
\end{align*}
 with probability approaching one. Exploiting the smoothness on $(\lambda_{l-1}^{0},\,\lambda_{l}^{0}]$,
we have
\[
\sup_{u\in(\lambda_{l-1}^{0},\,\lambda_{l}^{0}]}\sup_{\omega\in\left[-\pi,\,\pi\right]}\left|C_{l}\left(u,\,\omega\right)\right|=O_{\mathbb{P}}\left(\sqrt{\log\left(T\right)}\right).
\]
This implies 
\begin{align*}
\mathrm{D}_{r,T}\left(\omega\right)=\left|d_{l}\left(r/T,\,\omega\right)+C_{l}\left(r/T,\,\omega\right)\right| & =\left(d_{l}\left(r/T,\,\omega\right)+\mathrm{sign}\left(C_{l}\left(r/T,\,\omega\right)\right)\left|C_{l}\left(r/T,\,\omega\right)\right|\right),
\end{align*}
for each $r=\overline{r}_{l,b}-\left\lfloor m_{T}/B\right\rfloor ,\ldots,\,\overline{r}_{l,b},$
where $B$ is any integer with $1<B<\infty.$ By the definition of
$d_{l}\left(\cdot,\,\omega_{l}\right),$ for $\kappa_{T}\in\left[0,\,m_{T}/\left(BT\right)\right]$
we have 
\begin{align*}
d_{l}\left(\overline{r}_{l,b}/T,\,\omega_{l}\right)-d_{l}\left(\overline{r}_{l,b}/T-\kappa_{T},\,\omega_{l}\right) & =\left\lfloor \kappa_{T}T\right\rfloor m_{S,T}^{-1}\delta_{l,T}M_{S,T}^{-1/2}.
\end{align*}
In order to apply Lemma \ref{Lemma D1 BJV}, we need to choose $\kappa_{T}$
such that $\left\lfloor \kappa_{T}T\right\rfloor \delta_{l,T}M_{S,T}^{-1/2}/m_{S,T}\sqrt{\log\left(T\right)}\geq1$
or $m_{S,T}\sqrt{M_{S,T}\log\left(T\right)}/\delta_{l,T}T=o\left(\kappa_{T}\right).$
Lemma \ref{Lemma D1 BJV} then yields 
\begin{align*}
\frac{\overline{r}_{l,b}}{T}\geq & \underset{r\in\left(\widehat{\mathcal{I}}\backslash\left\{ r:\,r>\overline{r}_{l,b}\right\} \right)}{\mathrm{argmax}}\max_{\omega\in\left[-\pi,\,\pi\right]}T^{-1}\mathrm{D}_{r,T}\left(\omega\right)=\underset{r\in\left(\widehat{\mathcal{I}}\backslash\left\{ r:\,r>\overline{r}_{l,b}\right\} \right)}{\mathrm{argmax}}T^{-1}\mathrm{D}_{r,T}\left(\omega_{l}\right)\geq\frac{\overline{r}_{l,b}}{T}-\kappa_{T}.
\end{align*}
The case $r>\overline{r}_{l,b}$ can be treated similarly by symmetry.
It results in 
\begin{align*}
\frac{\overline{r}_{l,b}}{T}\leq & \underset{r\in\left(\widehat{\mathcal{I}}\backslash\left\{ r:\,r<\overline{r}_{l,b}\right\} \right)}{\mathrm{argmax}}\max_{\omega\in\left[-\pi,\,\pi\right]}T^{-1}\mathrm{D}_{r,T}\left(\omega\right)=\underset{r\in\left(\widehat{\mathcal{I}}\backslash\left\{ r:\,r<\overline{r}_{l,b}\right\} \right)}{\mathrm{argmax}}T^{-1}\mathrm{D}_{r,T}\left(\omega_{l}\right)\leq\frac{\overline{r}_{l,b}}{T}+\kappa_{T}.
\end{align*}
Therefore, we conclude $|\widehat{\lambda}_{T}-\overline{r}_{l,b}/T|=O_{\mathbb{P}}\left(\kappa_{T}\right)\rightarrow0$.
Now set $\widehat{\mathcal{I}}=\widehat{\mathcal{I}}\backslash\{T\widehat{\lambda}_{T}(\mathcal{\widehat{I}})-v_{T},\ldots,\,T\widehat{\lambda}_{T}(\widehat{\mathcal{I}})+v_{T}\}$
and $\widehat{\mathcal{T}}=\widehat{\mathcal{T}}\cup\{T\widehat{\lambda}_{T}(\mathcal{\widehat{I}})\}$.
Since $\mathbb{P}((\mathbf{D}_{2})^{c})=0$ if there are still undetected
breaks, we can repeat the above steps (1)-(4). The final results are
$\mathbb{P}(|\widehat{\mathcal{T}}-m_{0}|>\epsilon_{2})\rightarrow0$
for any $\epsilon_{2}>0$ and, after ordering the elements of $\widehat{\mathcal{T}}$
in chronological order, $\sup_{1\leq l\leq m_{0}}|\widehat{\lambda}_{l,T}-\lambda_{l}^{0}|=O_{\mathbb{P}}(m_{S,T}\sqrt{M_{S,T}\log\left(T\right)}/(T\inf_{1\leq l\leq m_{0}}\delta_{l,T}))$.

Assume without loss of generality that $\delta_{1,T}\geq\delta_{2,T}\geq\cdots\geq\delta_{m_{0},T}$.
Let $\widehat{\lambda}_{T}^{\left(q\right)}$ $\left(q=1,\ldots,\,m_{0}\right)$
denote the $q$th break detected by the procedure. It remains to prove
that if $K\rightarrow\infty$ then $\widehat{\lambda}_{T}^{\left(q\right)}$
is consistent for $\lambda_{q}^{0}$ $\left(q=1,\ldots,\,m_{0}\right)$.
Consider the first break $\lambda_{1}^{0}$. In order for the algorithm
to return $\widehat{\lambda}_{T}^{\left(1\right)}$ such that $|\widehat{\lambda}_{T}^{\left(1\right)}-\lambda_{1}^{0}|\overset{\mathbb{P}}{\rightarrow}0$
we need the following event to occur with sufficiently high probability,
$\mathbf{W}=\{\mathrm{for\,}l=1\,\exists r\in\widehat{\mathcal{I}}\,\mathrm{and}\,k=1,\ldots,\,K\,\mathrm{s.t.}\,r_{r,k}^{\diamond}=T_{1}^{0}\}.$
Note that 
\begin{align*}
\mathbf{W}^{c} & =\left\{ T_{1}^{0}\,\mathrm{not\,sampled\,in\,}K\,\mathrm{draws}\,\mathrm{from\,}T_{1}^{0}-m_{T}+1,\ldots,\,T_{1}^{0}\,\mathrm{without\,replacement}\right\} .
\end{align*}
 Thus, 

\begin{align*}
1-\mathbb{P}\left(\mathbf{W}^{c}\right) & =1-\frac{m_{T}-1}{m_{T}}\times\frac{m_{T}-2}{m_{T}-1}\times\cdots\times\frac{m_{T}-K}{m_{T}-K+1}\\
 & =1-\frac{m_{T}-K}{m_{T}}\\
 & \rightarrow1,
\end{align*}
 only if $K=O\left(a_{T}m_{T}\right)$ with $a_{T}\in(0,\,1]$ such
that $a_{T}\rightarrow1.$ Note that $K\leq m_{T}$ by construction.
The same argument can be repeated for $l=2,\ldots,\,m_{0}$. $\square$ 

\section{\label{Section Sensitivity-Analysis}Sensitivity Analyses}

In this section, we conduct Monte Carlo simulations to assess how
the finite-sample performance of the test statistics and change-point
estimators change when we implement them with different choices for
the tuning parameters. Recall that our recommended choices  are $m_{T}=T^{0.66},\,n_{T}=T^{5/8}$,
$b_{W,T}=n_{T}^{-1/6}$ and $n_{\omega}=7.$ Since our choice for
$m_{T}$ and $n_{T}$ corresponds to the upper bound allowed by Condition
\ref{Condition n_T h_T}, here we consider smaller values of $m_{T}$
and $n_{T}$. We consider Model M1 and $T=250,\,500$ and 1000. Table
\ref{Table SM1 Size Robustness} shows that the null rejection rates
are not much affected by the change in the choice of the tuning parameters.
The null rejection rates become less accurate only for substantially
smaller of $m_{T}$, $n_{T}$ and $b_{W,T}$ are chosen to be too
small, the null rejection rates become less accurate. For example,
for $m_{T}=T^{0.58}$ and $n_{T}=T^{0.56}$ or for $b_{W,T}=n_{T}^{-0.25}$,
some of the tests show some over-rejection. The choice of the number
of frequencies $n_{\omega}$ and of their locations do not matter
much for the finite-sample performance of $\mathrm{S}_{\mathrm{Dmax}}$
and $\mathrm{R}_{\mathrm{Dmax}}$. Table \ref{Table SM1 Power Robustness}
shows the results about the power. Any tuning parameter choice results
in good monotonic power for all tests. Overall, the results suggest
that reasonable changes in the choice of the tuning parameters yield
little changes in the finite-sample performance of the tests. The
change in the results become larger as the smoothing bandwidths $m_{T}$,
$n_{T}$ and $b_{W,T}$ are set too small. 

We move to the results about the change-point estimator. We consider
Model M6 with $T=1000.$ Table \ref{Table SM1 Estimation Robustness}
shows that small changes in the tuning parameters result in little
changes in the precision of the change-point estimator and of the
estimator of the number of change-points.

\begin{table}[H]
\caption{\label{Table SM1 Size Robustness}Empirical small-sample size for
model M1}

\begin{centering}
{\small{}}%
\begin{tabular}{lccccccccc}
\hline 
{\footnotesize{}} & \multicolumn{3}{c}{{\footnotesize{}$m_{T}=T^{0.63},\,n_{T}=T^{0.6}$}} & \multicolumn{3}{c}{{\footnotesize{}$m_{T}=T^{0.6},\,n_{T}=T^{0.58}$}} & \multicolumn{3}{c}{{\footnotesize{}$m_{T}=T^{0.58},\,n_{T}=T^{0.56}$}}\tabularnewline
{\footnotesize{}$\alpha=0.05$} & {\footnotesize{}$T=250$} & {\footnotesize{}$T=500$} & {\footnotesize{}$T=1000$} & {\footnotesize{}$T=250$} & {\footnotesize{}$T=500$} & {\footnotesize{}$T=1000$} & {\footnotesize{}$T=250$} & {\footnotesize{}$T=500$} & {\footnotesize{}$T=1000$}\tabularnewline
\hline 
\hline 
{\footnotesize{}$\mathrm{S}_{\max,T}\left(0\right)$} & {\footnotesize{}0.079} & {\footnotesize{}0.061} & {\footnotesize{}0.056} & {\footnotesize{}0.059} & {\footnotesize{}0.054} & {\footnotesize{}0.059} & {\footnotesize{}0.109} & {\footnotesize{}0.088} & {\footnotesize{}0.076}\tabularnewline
{\footnotesize{}$\mathrm{S}_{\mathrm{Dmax},T}$} & {\footnotesize{}0.045} & {\footnotesize{}0.067} & {\footnotesize{}0.063} & {\footnotesize{}0.058} & {\footnotesize{}0.070} & {\footnotesize{}0.058} & {\footnotesize{}0.132} & {\footnotesize{}0.069} & {\footnotesize{}0.062}\tabularnewline
{\footnotesize{}$\mathrm{R}_{\mathrm{max},T}\left(0\right)$} & {\footnotesize{}0.088} & {\footnotesize{}0.085} & {\footnotesize{}0.077} & {\footnotesize{}0.109} & {\footnotesize{}0.096} & {\footnotesize{}0.082} & {\footnotesize{}0.169} & {\footnotesize{}0114} & {\footnotesize{}0.102}\tabularnewline
{\footnotesize{}$\mathrm{R}_{\mathrm{Dmax},T}$} & {\footnotesize{}0.039} & {\footnotesize{}0.050} & {\footnotesize{}0.034} & {\footnotesize{}0.093} & {\footnotesize{}0.086} & {\footnotesize{}0.072} & {\footnotesize{}0.101} & {\footnotesize{}0.092} & {\footnotesize{}0.074}\tabularnewline
{\footnotesize{}} & \multicolumn{3}{c}{{\footnotesize{}$b_{W,T}=n_{T}^{-0.15}$}} & \multicolumn{3}{c}{{\footnotesize{}$b_{W,T}=n_{T}^{-0.2}$}} & \multicolumn{3}{c}{{\footnotesize{}$b_{W,T}=n_{T}^{-0.25}$}}\tabularnewline
{\footnotesize{}$\alpha=0.05$} & {\footnotesize{}$T=250$} & {\footnotesize{}$T=500$} & {\footnotesize{}$T=1000$} & {\footnotesize{}$T=250$} & {\footnotesize{}$T=500$} & {\footnotesize{}$T=1000$} & {\footnotesize{}$T=250$} & {\footnotesize{}$T=500$} & {\footnotesize{}$T=1000$}\tabularnewline
\hline 
\hline 
{\footnotesize{}$\mathrm{S}_{\max,T}\left(0\right)$} & {\footnotesize{}0.025} & {\footnotesize{}0.038} & {\footnotesize{}0.037} & {\footnotesize{}0.072} & {\footnotesize{}0.068} & {\footnotesize{}0.063} & {\footnotesize{}0.084} & {\footnotesize{}0.103} & {\footnotesize{}0.092}\tabularnewline
{\footnotesize{}$\mathrm{S}_{\mathrm{Dmax},T}$} & {\footnotesize{}0.032} & {\footnotesize{}0.041} & {\footnotesize{}0.044} & {\footnotesize{}0.047} & {\footnotesize{}0.061} & {\footnotesize{}0.057} & {\footnotesize{}0.021} & {\footnotesize{}0.067} & {\footnotesize{}0.096}\tabularnewline
{\footnotesize{}$\mathrm{R}_{\mathrm{max},T}\left(0\right)$} & {\footnotesize{}0.032} & {\footnotesize{}0.031} & {\footnotesize{}0.045} & {\footnotesize{}0.063} & {\footnotesize{}0.061} & {\footnotesize{}0.058} & {\footnotesize{}0.094} & {\footnotesize{}0.089} & {\footnotesize{}0.084}\tabularnewline
{\footnotesize{}$\mathrm{R}_{\mathrm{Dmax},T}$} & {\footnotesize{}0.029} & {\footnotesize{}0.031} & {\footnotesize{}0.032} & {\footnotesize{}0.027} & {\footnotesize{}0.038} & {\footnotesize{}0.042} & {\footnotesize{}0.091} & {\footnotesize{}0.085} & {\footnotesize{}0.082}\tabularnewline
{\footnotesize{}} & \multicolumn{3}{c}{{\footnotesize{}$n_{\omega}=15$}} & \multicolumn{3}{c}{{\footnotesize{}$n_{\omega}=11$}} & \multicolumn{3}{c}{{\footnotesize{}$n_{\omega}=5$}}\tabularnewline
{\footnotesize{}$\alpha=0.05$} & {\footnotesize{}$T=250$} & {\footnotesize{}$T=500$} & {\footnotesize{}$T=1000$} & {\footnotesize{}$T=250$} & {\footnotesize{}$T=500$} & {\footnotesize{}$T=1000$} & {\footnotesize{}$T=250$} & {\footnotesize{}$T=500$} & {\footnotesize{}$T=1000$}\tabularnewline
\hline 
\hline 
{\footnotesize{}$\mathrm{S}_{\mathrm{Dmax},T}$} & {\footnotesize{}0.011} & {\footnotesize{}0.037} & {\footnotesize{}0.043} & {\footnotesize{}0.016} & {\footnotesize{}0.037} & {\footnotesize{}0.039} & {\footnotesize{}0.023} & {\footnotesize{}0.070} & {\footnotesize{}0.064}\tabularnewline
{\footnotesize{}$\mathrm{R}_{\mathrm{Dmax},T}$} & {\footnotesize{}0.022} & {\footnotesize{}0.024} & {\footnotesize{}0.030} & {\footnotesize{}0.024} & {\footnotesize{}0.027} & {\footnotesize{}0.033} & {\footnotesize{}0.026} & {\footnotesize{}0.032} & {\footnotesize{}0.035}\tabularnewline
\hline 
\end{tabular}{\small\par}
\par\end{centering}
\end{table}

\begin{table}[H]
\caption{\label{Table SM1 Power Robustness}Empirical small-sample power for
model M1}

\begin{centering}
{\footnotesize{}}%
\begin{tabular}{lccccccccc}
\hline 
 & \multicolumn{3}{c}{{\footnotesize{}$m_{T}=T^{0.63},\,n_{T}=T^{0.6}$}} & \multicolumn{3}{c}{{\footnotesize{}$m_{T}=T^{0.6},\,n_{T}=T^{0.58}$}} & \multicolumn{3}{c}{{\footnotesize{}$m_{T}=T^{0.58},\,n_{T}=T^{0.56}$}}\tabularnewline
{\footnotesize{}$\alpha=0.05$} & {\footnotesize{}$T=250$} & {\footnotesize{}$T=500$} & {\footnotesize{}$T=1000$} & {\footnotesize{}$T=250$} & {\footnotesize{}$T=500$} & {\footnotesize{}$T=1000$} & {\footnotesize{}$T=250$} & {\footnotesize{}$T=500$} & {\footnotesize{}$T=1000$}\tabularnewline
\hline 
\hline 
{\footnotesize{}$\mathrm{S}_{\max,T}\left(0\right)$} & {\footnotesize{}0.686} & {\footnotesize{}0.791} & {\footnotesize{}0.913} & {\footnotesize{}0.605} & {\footnotesize{}0.748} & {\footnotesize{}0.886} & {\footnotesize{}0.682} & {\footnotesize{}0.783} & {\footnotesize{}0.892}\tabularnewline
{\footnotesize{}$\mathrm{S}_{\mathrm{Dmax},T}$} & {\footnotesize{}0.771} & {\footnotesize{}0.846} & {\footnotesize{}0.936} & {\footnotesize{}0.747} & {\footnotesize{}0.828} & {\footnotesize{}0.936} & {\footnotesize{}0.822} & {\footnotesize{}0.914} & {\footnotesize{}0.942}\tabularnewline
{\footnotesize{}$\mathrm{R}_{\mathrm{max},T}\left(0\right)$} & {\footnotesize{}0.813} & {\footnotesize{}0.899} & {\footnotesize{}0.956} & {\footnotesize{}0.821} & {\footnotesize{}0.915} & {\footnotesize{}0.942} & {\footnotesize{}0.891} & {\footnotesize{}0.877} & {\footnotesize{}0.932}\tabularnewline
{\footnotesize{}$\mathrm{R}_{\mathrm{Dmax},T}$} & {\footnotesize{}0.582} & {\footnotesize{}0.745} & {\footnotesize{}0.871} & {\footnotesize{}0.649} & {\footnotesize{}0.827} & {\footnotesize{}0.864} & {\footnotesize{}0.784} & {\footnotesize{}0.791} & {\footnotesize{}0.863}\tabularnewline
 & \multicolumn{3}{c}{{\footnotesize{}$b_{W,T}=n_{T}^{-0.15}$}} & \multicolumn{3}{c}{{\footnotesize{}$b_{W,T}=n_{T}^{-0.2}$}} & \multicolumn{3}{c}{{\footnotesize{}$b_{W,T}=n_{T}^{-0.25}$}}\tabularnewline
{\footnotesize{}$\alpha=0.05$} & {\footnotesize{}$T=250$} & {\footnotesize{}$T=500$} & {\footnotesize{}$T=1000$} & {\footnotesize{}$T=250$} & {\footnotesize{}$T=500$} & {\footnotesize{}$T=1000$} & {\footnotesize{}$T=250$} & {\footnotesize{}$T=500$} & {\footnotesize{}$T=1000$}\tabularnewline
\hline 
\hline 
{\footnotesize{}$\mathrm{S}_{\max,T}\left(0\right)$} & {\footnotesize{}0.607} & {\footnotesize{}0.815} & {\footnotesize{}0.886} & {\footnotesize{}0.786} & {\footnotesize{}0.921} & {\footnotesize{}0.996} & {\footnotesize{}0.779} & {\footnotesize{}0.926} & {\footnotesize{}0.969}\tabularnewline
{\footnotesize{}$\mathrm{S}_{\mathrm{Dmax},T}$} & {\footnotesize{}0.628} & {\footnotesize{}0.899} & {\footnotesize{}0.907} & {\footnotesize{}0.754} & {\footnotesize{}0.921} & {\footnotesize{}0.995} & {\footnotesize{}0.759} & {\footnotesize{}0.915} & {\footnotesize{}0.953}\tabularnewline
{\footnotesize{}$\mathrm{R}_{\mathrm{max},T}\left(0\right)$} & {\footnotesize{}0.771} & {\footnotesize{}0.902} & {\footnotesize{}0.940} & {\footnotesize{}0.821} & {\footnotesize{}0.952} & {\footnotesize{}0.996} & {\footnotesize{}0.926} & {\footnotesize{}0.982} & {\footnotesize{}0.996}\tabularnewline
{\footnotesize{}$\mathrm{R}_{\mathrm{Dmax},T}$} & {\footnotesize{}0.458} & {\footnotesize{}0.641} & {\footnotesize{}0.704} & {\footnotesize{}0.604} & {\footnotesize{}0.821} & {\footnotesize{}0.952} & {\footnotesize{}0.852} & {\footnotesize{}0.928} & {\footnotesize{}0.952}\tabularnewline
 & \multicolumn{3}{c}{{\footnotesize{}$n_{\omega}=15$}} & \multicolumn{3}{c}{{\footnotesize{}$n_{\omega}=11$}} & \multicolumn{3}{c}{{\footnotesize{}$n_{\omega}=3$}}\tabularnewline
{\footnotesize{}$\alpha=0.05$} & {\footnotesize{}$T=250$} & {\footnotesize{}$T=500$} & {\footnotesize{}$T=1000$} & {\footnotesize{}$T=250$} & {\footnotesize{}$T=500$} & {\footnotesize{}$T=1000$} & {\footnotesize{}$T=250$} & {\footnotesize{}$T=500$} & {\footnotesize{}$T=1000$}\tabularnewline
\hline 
\hline 
{\footnotesize{}$\mathrm{S}_{\mathrm{Dmax},T}$} & {\footnotesize{}0.769} & {\footnotesize{}0.977} & {\footnotesize{}0.907} & {\footnotesize{}0.731} & {\footnotesize{}0.966} & {\footnotesize{}0.962} & {\footnotesize{}0.645} & {\footnotesize{}0.756} & {\footnotesize{}0.864}\tabularnewline
{\footnotesize{}$\mathrm{R}_{\mathrm{Dmax},T}$} & {\footnotesize{}0.401} & {\footnotesize{}0.635} & {\footnotesize{}0.708} & {\footnotesize{}0.455} & {\footnotesize{}0.634} & {\footnotesize{}0.701} & {\footnotesize{}0.742} & {\footnotesize{}0.756} & {\footnotesize{}0.834}\tabularnewline
\hline 
\end{tabular}{\footnotesize\par}
\par\end{centering}
\end{table}

\begin{table}[H]
\caption{\label{Table SM1 Estimation Robustness}Empirical distribution of
$\widehat{m}-m_{0}$ for model M6}

\begin{centering}
\begin{tabular}{cccccccccc}
\hline 
\multicolumn{5}{c}{$m_{T}=T^{0.63},\,n_{T}=T^{0.6}$} & \multicolumn{5}{c}{$m_{T}=T^{0.60},\,n_{T}=T^{0.58}$}\tabularnewline
\hline 
Percent time $\widehat{m}=m_{0}$ &  & $Q_{0.25}$ & Median & $Q_{0.75}$ & Percent time $\widehat{m}=m_{0}$ &  & $Q_{0.25}$ & Median & $Q_{0.75}$\tabularnewline
\hline 
\hline 
80.30 & $\widehat{T}_{1}$ & 297 & 328 & 348 & 78.60 & $\widehat{T}_{1}$ & 300 & 328 & 346\tabularnewline
 & $\widehat{T}_{2}$ & 605 & 650 & 680 &  & $\widehat{T}_{2}$ & 614 & 645 & 675\tabularnewline
\hline 
 &  &  &  &  &  &  &  &  & \tabularnewline
\multicolumn{5}{c}{$n_{\omega}=15$} & \multicolumn{5}{c}{$n_{\omega}=11$}\tabularnewline
Percent time $\widehat{m}=m_{0}$ &  & $Q_{0.25}$ & Median & $Q_{0.75}$ & Percent time $\widehat{m}=m_{0}$ &  & $Q_{0.25}$ & Median & $Q_{0.75}$\tabularnewline
 & $\widehat{T}_{1}$ &  &  &  &  &  &  &  & \tabularnewline
 & $\widehat{T}_{2}$ &  &  &  &  &  &  &  & \tabularnewline
\multicolumn{10}{c}{}\tabularnewline
\multicolumn{5}{c}{$n_{\omega}=11$} & \multicolumn{5}{c}{$n_{\omega}=3$ }\tabularnewline
Percent time $\widehat{m}=m_{0}$ &  & $Q_{0.25}$ & Median & $Q_{0.75}$ & Percent time $\widehat{m}=m_{0}$ &  & $Q_{0.25}$ & Median & $Q_{0.75}$\tabularnewline
 & $\widehat{T}_{1}$ &  &  &  & 0.814 & $\widehat{T}_{1}$ & 299 & 333 & 352\tabularnewline
 & $\widehat{T}_{2}$ &  &  &  &  & $\widehat{T}_{2}$ & 617 & 657 & 694\tabularnewline
\end{tabular}
\par\end{centering}
\end{table}

\section{\label{Section Additional-Monte-Carlo}Additional Monte Carlo Results}

In this section, we report simulations results for Model M6 where
the errors are drawn from the $t_{\nu}$ distribution. In particular,
in Model M6 $e_{t}\sim\mathrm{i.i.d.}\,t_{\nu}$ with $\nu=5,\,10$.
Table \ref{Table SM1-M2 t-distr}-\ref{Table Power M3-M5 t-dstr}
shows that the proposed test statistics have accurate null rejection
rates and good monotonic power similar to the case of Gaussian errors.
Note that the statistic $\widehat{D}$ continue to establish large
size. 

\begin{table}[H]
\caption{\label{Table SM1-M2 t-distr}Empirical small-sample size for model
M1 with $t$-distributed errors}

\begin{centering}
\begin{tabular}{lcccccc}
\hline 
 & \multicolumn{3}{c}{$t_{\nu},\,\nu=5$} & \multicolumn{3}{c}{$t_{\nu},\,\nu=10$}\tabularnewline
$\alpha=0.05$ & $T=250$ & $T=500$ & $T=1000$ & $T=250$ & $T=500$ & $T=1000$\tabularnewline
\hline 
\hline 
$\mathrm{S}_{\max,T}\left(0\right)$ & 0.065 & 0.106 & 0.097 & 0.039 & 0.054 & 0.035\tabularnewline
$\mathrm{S}_{\mathrm{Dmax},T}$ & 0.054 & 0.097 & 0.112 & 0.028 & 0.063 & 0.051\tabularnewline
$\mathrm{R}_{\mathrm{max},T}\left(0\right)$ & 0.071 & 0.056 & 0.049 & 0.038 & 0.065 & 0.039\tabularnewline
$\mathrm{R}_{\mathrm{Dmax},T}$ & 0.043 & 0.014 & 0.007 & 0.008 & 0.026 & 0.011\tabularnewline
$\widehat{D}$ statistic & 0.661 & 0.515 & 0.083 & 0.598 & 0.454 & 0.051\tabularnewline
\hline 
\end{tabular}
\par\end{centering}
\end{table}

\begin{table}[H]
\caption{\label{Table Power M3-M5 t-dstr}Empirical small-sample power for
model M6 with $t$-distributederrors}

\begin{centering}
\begin{tabular}{lcccccc}
\hline 
 & \multicolumn{3}{c}{$t_{\nu},\,\nu=5$} & \multicolumn{3}{c}{$t_{\nu},\,\nu=10$}\tabularnewline
$\alpha=0.05$ & $T=250$ & $T=500$ & $T=1000$ & $T=250$ & $T=500$ & $T=1000$\tabularnewline
\hline 
\hline 
$\mathrm{S}_{\max,T}\left(0\right)$ & 0.777 & 0.896 & 0.900 & 0.738 & 0.875 & 0.915\tabularnewline
$\mathrm{S}_{\mathrm{Dmax},T}$ & 0.788 & 0.894 & 0.908 & 0.717 & 0.859 & 0.890\tabularnewline
$\mathrm{R}_{\mathrm{max},T}\left(0\right)$ & 0.782 & 0.928 & 0.934 & 0.734 & 0.920 & 0.967\tabularnewline
$\mathrm{R}_{\mathrm{Dmax},T}$ & 0.582 & 0.803 & 0.800 & 0.522 & 0.796 & 0.895\tabularnewline
$\widehat{D}$ statistic & 0.989 & 0.996 & 0.992 & 0.982 & 0.979 & 0.856\tabularnewline
\hline 
\end{tabular}
\par\end{centering}
\end{table}

\newpage{}

\bibliographystyleReferencesSupp{elsarticle-harv}  
\bibliographyReferencesSupp{References_Supp}

\clearpage{}

\end{singlespace}

\end{document}